\documentclass[a4paper, 10pt]{article}

\usepackage[top=2.5cm, right = 4cm, bottom =2.5cm, left = 4cm]{geometry}
\usepackage{lmodern}
\usepackage{amsmath, amssymb, amsthm}
\usepackage{bm} 

\newtheorem{lem}{Lemma}[section]
\newtheorem{thm}[lem]{Theorem}
\newtheorem{prop}[lem]{Proposition}

\theoremstyle{plain}
\newtheorem{define}[lem]{Definition}
\newtheorem{ex}[lem]{Example}
\theoremstyle{remark}
\newtheorem{rem}[lem]{Remark}
\numberwithin{equation}{section}

\usepackage{enumitem}
\usepackage{xcolor}
\usepackage{pgfplots}
\pgfplotsset{compat=newest}
\usetikzlibrary{plotmarks}
\usetikzlibrary{arrows.meta}
\usepgfplotslibrary{patchplots}
\usepackage{grffile}
\usepackage{graphicx}
\graphicspath{{./Figs/}}
\usepackage{subcaption}

\usepackage[]{csquotes}
\newcommand{\R}{\mathbb{R}}
\newcommand{\Rd}[1][d]{%
\R^{#1}%
}

\newcommand\vc[1]{\bm{#1}}

\newcommand{\dual}[1]{{#1}^{\prime}}

\newcommand{\dualp}[3]{
    \langle {#1}, {#2} \rangle_{\dual{#3},{#3}}
}

\newcommand{\Lp}[2][\Omega]{\ensuremath{L^{#2}\left(#1\right)}}
\newcommand{\hone}[1][\Omega]{\ensuremath{H^1\left(#1\right)}}
\newcommand{\hzero}[1][\Omega]{\ensuremath{H^1_0\left(#1\right)}}
\newcommand{\sob}[3][\Omega]{\ensuremath{W^{#2,#3}\left(#1\right)}}
\newcommand{\sobzero}[3][\Omega]{\ensuremath{W^{#2,#3}_0\left(#1\right)}}

\newcommand\dm{\mathop{}\!\mathrm{d}}

\newcommand{\argmin}[1]{
	\mathrm{arg}\, \min_{#1}
}                                                   
\newcommand{\sign}{\mathrm{sgn}}                    

\begin{document}
\title{Eliminating Gibbs Phenomena: A Non-linear Petrov-Galerkin Method for the Convection-Diffusion-Reaction Equation}
\author{Paul Houston$^{a,}$\thanks{Paul.Houston@nottingham.ac.uk}, Sarah Roggendorf$^{a,}$\thanks{Corresponding author: Sarah.Roggendorf@nottingham.ac.uk} and Kristoffer G. van der Zee$^{a,}$\thanks{KG.vanderZee@nottingham.ac.uk} \bigskip\\
{\small $^a$School of Mathematical Sciences, The University of Nottingham},\\ {\small University Park, NG7 2RD, United Kingdom}}
%
%
%
%
%

\maketitle
\begin{abstract}
  In this article we consider the numerical approximation of the con\-vec\-tion-dif\-fusion-reaction equation. One
  of the main challenges of designing a numerical method for this problem
  is that boundary layers occurring in the convection-dominated case can lead to non-physical oscillations in the numerical
  approximation, often referred to as Gibbs phenomena.
  The idea of this article is to consider the approximation problem as a residual
  minimization in dual norms in $L^q$-type Sobolev spaces, with $1<q<\infty$. We
  then apply a non-standard, non-linear Petrov-Galerkin discretization, that is applicable to reflexive Banach spaces such that the
  space itself and its dual are strictly convex.
  Similar to discontinuous Petrov-Galerkin methods, this method is
  based on  minimizing the residual in a dual norm. Replacing the  intractable dual
  norm  by a suitable discrete dual norm gives rise to a non-linear inexact
  mixed method. This generalizes the
  Petrov-Galerkin framework developed in the context of discontinuous Petrov-Galerkin
  methods to more general Banach spaces. For the convection-diffusion-reaction equation,
  this yields a generalization of a similar approach from
  the $L^2$-setting to the $L^q$-setting.
  A key advantage of considering a more general Banach space setting is that, in
  certain cases, the oscillations in the numerical approximation vanish as $q$ tends to $1$, as we will demonstrate using a few simple numerical examples.
  \bigskip

  \bigskip

\end{abstract}
\section{Introduction}
The term Gibbs phenomenon originally refers to the effect that the partial sums
of the Fourier series approximating a function with a jump discontinuity exhibit
over- and undershoots near the discontinuity. The phenomenon is named after
Willard Gibbs who described it in 1899 \cite{Gibbs1899}, though it had already been discovered
earlier by Henry Wilbraham in 1848 \cite{Wilbraham1848}.  It also occurs in the best approximation by spline
functions in the $L^2$-metric \cite{Richards1991} and is one of the main challenges in the numerical approximation of partial differential equations
(PDEs) whose solutions contain sharp features such as shocks or thin layers.
In \cite{Saff1999} it is shown in one dimension that the best approximation of
jump discontinuities by polygonal lines on a uniform grid in one dimension does not lead to Gibbs phenomenon in
$L^1$. More precisely, it is shown that the overshoot of the best approximation
in $L^q$, $1<q<\infty$, is an increasing function of $q$ that tends to $0$ as $q$ tends to $1$.  We show in  \cite{Houston2019a} that this result is still true for the best approximation in $L^q$ by piecewise linear functions on certain meshes in two dimensions and certain non-uniform meshes in one dimension. However, there exist meshes both in one and two dimensions such that the (maximal) overshoot tends to some $\alpha >0$ as $q \rightarrow 1$. Nonetheless, it is suggested in \cite{Houston2019a} that if the location of the discontinuity is known, a mesh can be constructed in such a way that the overshoot vanishes.
In \cite{Moskona1995}, it is shown that Gibbs phenomena occur in the best approximation of a function containing a jump discontinuity by a trigonometric polynomial in the $L^1$-metric. This confirms that the choice of the finite dimensional approximation space crucially determines whether Gibbs phenomena can be eliminated by considering best approximation in $L^1$.

The motivation of this article is to exploit  the Gibbs-reducing property of $L^q$-type spaces
in the numerical approximation of PDEs using finite element methods.
A similar idea has already been pursued for numerical methods in $L^q$ by Guermond
\cite{Guermond2004}.
In \cite{Guermond2004}, Guermond points out
that there are only very few attempts to approximate PDEs directly in $L^1$ despite the fact that
first-order PDEs and their non-linear generalizations have been extensively studied
in $L^1$. The existing numerical methods to achieve this include the ones outlined in the articles by Lavery \cite{Lavery1988,Lavery1989,Lavery1991},
the reweighted least-squares method of Jiang \cite{Jiang1993,Jiang1998} and the methods outlined
in the series of articles by Guermond et al.\ \cite{Guermond2004,Guermond2007,Guermond2008,Guermond2008/09,Guermond2009}.

More recently, a novel approach to designing finite element methods in
a very general Banach space setting has been introduced in \cite{Muga2017}.
This approach is rooted in the so-called Discontinuous Petrov-Galerkin (DPG) methods  \cite{Demkowicz2014} and
extends the concept of optimal test norms and functions from Hilbert spaces to
more general Banach spaces yielding a scheme that can be interpreted as a \emph{non-linear} Petrov-Galerkin method. At least in an abstract sense, this approach outlines
how to design a numerical method that leads to a quasi-best approximation of the
solution in a space of choice, provided the continuous problem is well-posed in a suitable sense.
In this article we  apply this abstract approach to
the convection-diffusion-reaction equation in a $\sobzero{1}{q}$-$\sobzero{1}{q'}$ setting where $1/q+1/q' = 1$ and study the effect this has on the numerical approximation of boundary layers. We will see that as $q\rightarrow 1$, the Gibbs phenomenon can be eliminated entirely on some meshes. We will also consider certain choices of meshes for which this is not the case and consider how this can be fixed.

\subsection{Notation}
Throughout this article, we denote by $\Lp{q}$, $1\le q<\infty$, the Lebesgue space of $q$-integrable
functions on a bounded Lipschitz domain $\Omega \subset \Rd$, $d\in\{1,2,3,\ldots\}$;
$L^\infty(\Omega)$ is the Lebesgue space of functions on~$\Omega$ with finite essential supremum;
and $\sob{1}{q}$, $1\le q\le \infty$, is the Sobolev space of functions that are in $\Lp{q}$
such that their gradient is in $\Lp{q}^d$. Furthermore, $\sobzero{1}{q}\subset
\sob{1}{q}$ is the subspace of all functions with zero trace  on the boundary
$\partial\Omega$. The corresponding norms are denoted by $\|\cdot \|_{\Lp{q}}$
and $\|\cdot\|_{\sob{1}{q}}$, respectively, and the Sobolev-semi norm on $\sob{1}{q}$ is given
by $|\cdot|_{\sob{1}{q}}$. For $q=2$, we furthermore use the usual
notation $H^1(\Omega) := \sob{1}{2}$ and $\hzero:= \sobzero{1}{2}$.
For $1\leq q\leq \infty$ we write $q'$ to denote the dual exponent such that
$1/q+1/q' = 1$.  For any Banach space $V$, its norm is denoted by $\|\cdot \|_{V}$ and
its dual space by $V'$; the dual space of $\sobzero{1}{q}$ is given by $\sob{-1}{q'}$ and
$H^{-1}(\Omega):=\sob{-1}{2}$. For $v \in V$ and $\varphi \in V'$, we have
the duality pairing
\begin{align*}
  \langle \varphi, v \rangle_{V', V} : = \varphi(v).
\end{align*}
For any $\varphi \in V'$, its norm in the dual space $V'$ is given by
\begin{align*}
  \|\varphi\|_{V'} := \sup_{v \in V, v \neq 0} \frac{\langle \varphi, v \rangle_{V', V}}{\|v\|_V}.
\end{align*}
If $V$ is a Hilbert space, we denote the Riesz map on $V$ by $R_V$.
\subsection{Problem Statement}
The approach to generalizing the DPG framework to general Banach spaces introduced
in \cite{Muga2017} is based on residual minimization. To this end, the authors consider the abstract problem: find $u \in U$ such that
\begin{align*}
  Bu = f \text{ in } V',
\end{align*}
where $U$ and $V$ are Banach spaces, $B:U\rightarrow V'$ is a continuous, bounded-below
linear operator, and the right hand side $f$ is in the dual space $V'$. The associated
residual minimization problem for a given finite dimension subspace $U_n\subset U$, e.g., a finite element space, of dimension $n$ is defined as follows:
\begin{align*}
  u_n = \argmin{w_n \in U_n} \| f- Bw_n\|_{V'}.
\end{align*}
It can be shown that the solution $u_n$ is a quasi-best approximation of the analytical solution $u$, in the sense that
\begin{align*}
  \|u-u_n\|_U \leq \frac{M}{\gamma}\inf_{w_n \in U_n} \|u-w_n\|_U,
\end{align*}
where $M$ is the continuity constant of $B$ and $\gamma$ the bounded-below constant of $B$.

In \cite{Muga2017} it is shown that the minimization problem can be reformulated yielding the following non-linear mixed system: find $(r, u_n) \in V \times U_n$ such that
\begin{align}
  \begin{aligned}
 \langle J_V(r), v \rangle_{V', V} + \langle B u_n, v \rangle_{V', V} &= \langle f, v\rangle_{V', V}, && \forall v \in V, \\
 \langle B w_n, r \rangle_{V', V} & = 0, && \forall w_n \in U_n,
\end{aligned}\label{eq:npg_exact_mixed_method}
\end{align}
where $J_V$ is a so-called duality mapping. Duality mappings are a generalization of the Riesz map $R_V$ to general Banach spaces with the main difference that they are non-linear mappings unless $V$ is a Hilbert space. To turn the above mixed system into a practical method, $V$ is additionally replaced by a finite dimensional subspace $V_m$ of dimension $m$. As a result the minimization problem is no longer solved exactly. However, if  the spaces $U_n$ and $V_m$ are chosen in a suitable way, one can obtain a well-posed fully discrete mixed system that retains the quasi-best approximation property of $u_n$ with a modified constant.

The aim of this article is to apply this abstract framework to the convection-diffusion-reaction equation
\begin{align}
  \begin{aligned}
-\varepsilon \Delta u + \vc{b}\cdot\nabla u +cu&= f \qquad \text{ in }\Omega,\\
u&= 0 \qquad\text{ on } \Gamma = \partial \Omega,
\end{aligned}\label{eq:intro_conv_diff}
\end{align}
where $\varepsilon$, $b:\Omega \rightarrow \Rd$, and $c:\Omega \rightarrow \R$ are the
(positive) diffusion parameter, convection field and reaction coefficient, respectively,
and $f:\Omega \rightarrow \R$ is a given source term.
Therefore, we will consider the linear operator associated with the following bilinear form:
\begin{align*}
\mathcal{B}_{\varepsilon}(u,v) &= \varepsilon\int_{\Omega} \nabla u \cdot \nabla v \dm \vc{x} - \int_{\Omega}u \nabla \cdot (\vc{b}v) \dm \vc x + \int_{\Omega} c u v \dm \vc x.
\end{align*}

The resulting method is a generalization to the $\sobzero{1}{q}$-$\sobzero{1}{q'}$-setting of the method introduced for the $\hzero$-setting in \cite{Chan2014}, where the DPG framework is applied to the convection-diffusion-reaction equation without introducing broken test spaces.
Related Petrov-Galerkin formulations are studied in \cite{Broersen2014,Broersen2015,Chan2014a,Cohen2012,Demkowicz2013}.
\subsection{Summary of Results}

In this article we study how the mixed method \eqref{eq:npg_exact_mixed_method} and its fully discrete counterpart can be applied to the convection-diffusion-reaction equation \eqref{eq:intro_conv_diff}.
In particular, we study the choice of the space $V$, the corresponding norm $\|\cdot\|_V$, and  the discrete space $V_m$. The norm $\|\cdot\|_V$ on the space $V$ determines the operator $J_V$ in \eqref{eq:npg_exact_mixed_method} and thus crucially determines the resulting method. As in the special case $q=2$ investigated in \cite{Chan2014}, we introduce weakly imposed boundary conditions on the inflow boundary in the space $V_m$ to address robustness issues. We demonstrate in one dimension, how this choice of boundary conditions affects the approximation and compare this approach with the case when a weighted norm on $V$ is employed.

The main focus of our numerical investigation is eliminating Gibbs phenomena in the finite element approximation. Indeed, we will see that the numerical approximation generated by our method qualitatively behaves like the $\Lp{q}$-best approximation of the analytical solution. Thus, the Gibbs phenomenon can be eliminated by taking the limit $q \rightarrow 1$ provided that the $\Lp{1}$-best approximation does not exhibit Gibbs phenomena.
The results in  \cite{Houston2019a} show that this depends on the mesh that is chosen. In one dimension only certain non-uniform grids have the property that the $\Lp{1}$-best approximation contains over- and undershoots, whereas in higher dimensions this can even occur on structured, uniform meshes. We will demonstrate for one two-dimensional example that it is possible to use the insights from \cite{Houston2019a} to modify the mesh in order to eliminate Gibbs phenomena as $q \rightarrow 1$.

\subsection{Outline of the Paper}
This article is structured as follows: in Section \ref{sec:duality_map} we introduce duality mappings, which are essential to the abstract framework in \cite{Muga2017},  cf., \eqref{eq:npg_exact_mixed_method}. In Section \ref{sec:minres} we give a brief overview of the abstract framework in \cite{Muga2017}.
Then, in Section \ref{sec:confusion} we then apply this framework to the convection-diffusion-reaction equation and illustrate how this yields a practical method that can be implemented, e.g., in FEniCS \cite{AlnaesBlechta2015a,LoggMardalEtAl2012a}.
In Section \ref{sec:numerics} we investigate the practical performance of the proposed method for a range of different test problems in one and two dimensions. Finally, in  \ref{sec:conclusion}, we summarize the work presented in this article and highlight the potential and challenges of the proposed numerical method.

\section{Duality Mappings}\label{sec:duality_map}
One of the key ingredients used to extend concepts from Hilbert spaces to more general
Banach spaces both in \cite{Muga2015} and \cite{Houston2019} is replacing the Riesz
map with so-called duality mappings. In general, duality mappings are non-linear maps
that share certain properties with the Riesz map. If the underlying space is a Hilbert
space, any duality mapping is linear and the Riesz map is one possible choice for the duality
mapping. In this section we will give a precise definition of duality mappings, summarize
a range of very useful properties and list some relevant examples.

\begin{define}
  \leavevmode
\begin{enumerate}
\item A weight function  is a continuous and strictly increasing function
$\varphi: \R_+ \rightarrow \R_+$ such that $\varphi(0) = 0$ and
 $\lim_{t\rightarrow \infty} \varphi(t) = + \infty$.
\item
Let $V$ be a Banach space and $\varphi$ a weight function. Denote
by $\mathcal{P}(\dual{V})$ the power set of $\dual{V}$. Then the multivalued
map $\mathcal{J}_V^{\varphi}:V\rightarrow \mathcal{P}(\dual{V})$, defined by
\begin{align}
\mathcal{J}_{V}^{\varphi}(v):=\left\{\dual{v} \in \dual{V} \,:\,
\dualp{\dual{v}}{v}{V}=\|v\|_V\|\dual{v}\|_{\dual{V}},
\|\dual{v}\|_{\dual{V}}= \varphi(\|v\|_V)\right\}\label{eq:dualitymap}
\end{align}
is called a \emph{duality mapping of weight $\varphi$}.
\end{enumerate}
\end{define}
Due to a corollary of the Hahn-Banach theorem (cf., e.g., \cite[Corollary~1.3]{Brezis2011})
 the set $\mathcal{J}^{\varphi}_{V}(v)$ is non-empty. If we choose
 $\varphi(t) = t$ as a weight function we obtain the so-called \emph{normalized}
duality map. In many cases this choice is the simplest and most convenient one.
 In fact, some books only introduce the duality mapping for this choice of
 $\varphi$, cf., e.g., \cite{Brezis2011,Deimling1985,Zeidler1985}.
 As we will see later, however, $\varphi(t) = t^{q-1}$ is in some cases a
  more suitable choice for $L^{q}$-type spaces. If ${q}=2$, this choice
  coincides with the normalized duality map. Furthermore, note that we have
  $\mathcal{J}_{V}(v)=\{R_V(v)\}$ if $\mathcal{J}_{V}$ is the normalized duality
   map and $V$ is a Hilbert space due to the Riesz Representation theorem.

\subsection{Properties of Duality Mappings}
In the following proposition we summarize a few  properties of the duality map
 for Banach spaces with particular structures.
\begin{prop}\label{prop:duality_map_convex_spaces} Let $V$ be a Banach space
and denote by $\mathcal{J}_V^{\varphi}:V\rightarrow \mathcal{P}(\dual{V})$
the duality map of weight $\varphi$ on $V$. Then the following statements are
true:
\begin{enumerate}
\item  $\dual{V}$ is strictly convex\footnote{A Banach space $V$ is \emph{strictly convex} if for all $v_1, v_2 \in V$ such that $v_1 \neq v_2$ and $\|v_1\|_V = \|v_2\|_V =1$ it holds that $
  \|\vartheta v_1 + (1-\vartheta)v_2\|_V < 1 \quad \forall \, \vartheta \in (0,1).
$}
 if and only if $\mathcal{J}_V^{\varphi}$ is
  single valued, cf., \cite[Prop.~12.3]{Deimling1985}. In this case we define
  the duality map $J_V^{\varphi}:V\rightarrow \dual{V}$ such that
  $\mathcal{J}_V^{\varphi}(v) = \{J_V^{\varphi}(v)\}$ for all $v \in V$.
\item If $V$ is strictly convex, then $\mathcal{J}_V^{\varphi}(v)\cap
 \mathcal{J}_V^{\varphi}(w)=\emptyset$ for all $w\neq v$. In particular,
 $\mathcal{J}_V^{\varphi}$ is injective.
\item $V$ is reflexive if and only if $\mathcal{J}_V^{\varphi}$ is surjective
in the sense that for every $\dual{v}\in \dual{V}$ there is a $v\in V$ such that
$\dual{v} \in \mathcal{J}_V^{\varphi}(v)$, cf., \cite[Theorem~3.4, Chapter~II]{Cioranescu1990}.
\item If $V$ is a reflexive Banach space and $\mathcal{J}^{\varphi}_V$ is a
duality mapping of weight $\varphi$, then  $(\mathcal{J}^{\varphi}_V)^{-1}$ is
a duality mapping on $V'$ of weight $\varphi^{-1}$, cf.,
\cite[Cor.~3.5, Ch.~II]{Cioranescu1990}.
\end{enumerate}
\end{prop}
The main implication of the above proposition is that in a strictly convex and reflexive space $V$, the duality mapping is bijective and its inverse can be identified with a duality mapping on the dual space $V'$ with the inverse weight
$\varphi^{-1}$.

The following theorem is a special case of Theorem 4.4 in \cite[Chapter~I]{Cioranescu1990}
and states that the duality map on $V$ can be characterized using the subdifferential
of the norm on $V$. This is a key property of the duality map that will allow us
to derive the duality map for some specific Banach spaces in the special case
that the subdifferential is essentially the G\^ateaux or Fr\'echet derivative
of the norm.
\begin{thm}[Asplund, cf., {\cite[Ch.~I, Theorem~4.4]{Cioranescu1990}}]\label{thm:Asplund}
Let $V$ be a Banach space and define $F_V^{\varphi}:V\rightarrow \R$
by $F_V^{\varphi}(\cdot):= \psi(\|\cdot\|_V)$, where
$\psi(s) : = \int_{0}^s \varphi(t)\dm t$ and $\varphi$ is a weight function.
 Then for any $v\in V$, we have
\begin{align}
\mathcal{J}_V^{\varphi}(v) = \partial F_V^{\varphi}(v),
\end{align}
where $\partial F_V^{\varphi}(v)$ denotes the subdifferential of $F_V^{\varphi}$
at $v$.
\end{thm}
The result that finally allows us to compute duality maps is the following.
\begin{prop}[cf., {\cite[Proposition~47.19]{Zeidler1985}}]
  Let $V$ be a Banach space and define $F_V^{\varphi}:V\rightarrow \R$
  by $F_V^{\varphi}(\cdot):= \psi(\|\cdot\|_V)$, where
  $\psi(s) : = \int_{0}^s \varphi(t)\dm t$ and $\varphi$ is a weight function.
\begin{enumerate}
\item If $\dual{V}$ is strictly convex, then $\nabla F_V^{\varphi}(v)$ exists as a
G\^ateaux derivative and $\nabla F_V^{\varphi}(v)=J_V^{\varphi}(v)$ for all $v\in V$.
\item If $\dual{V}$ is uniformly convex, then $\nabla F_V^{\varphi}(v)$ exists as a
Fr\'echet derivative and $\nabla F_V^{\varphi}(v)=J_V^{\varphi}(v)$ for all $v\in V$.
\end{enumerate}
\end{prop}
Note that  uniform convexity
of $\dual{V}$ implies strict convexity of $\dual{V}$ (cf., \cite{Chidume2009}) and due to the
Milman-Pettis Theorem, $\dual{V}$ and hence also
$V$ are reflexive in this case.  Furthermore, if $F_V^{\varphi}$ is
G\^ateaux differentiable, Theorem \ref{thm:Asplund} guarantees that
the duality map is single valued which implies strict convexity of
the dual space $\dual{V}$ due to Proposition \ref{prop:duality_map_convex_spaces}.

\subsection{Some Examples of Duality Mappings on Sobolev Spaces} \label{sec:duality_map_example}
First, consider $V= \Lp{{q}}$ with the norm
\begin{align}
\|v\|_{\Lp{{q}}} = \left(\int_{\Omega} |v|^{q} \dm \vc x\right)^{1/{q}}.
\end{align}
Let us denote the duality mapping on $\Lp{{q}}$ with weight $\varphi(t) = t^{{q}-1}$
by $J_{q}$. In this case $\psi(s) = \int_0^s \varphi(t) \dm t = \frac{1}{{q}}s^{q}$
 and thus we can compute
\begin{align}
\begin{aligned}
\langle J_{q}(v), w \rangle_{\Lp{{q'}}, \Lp{{q}}} &= \left.\frac{\dm}{\dm t}
\left(\frac{1}{{q}}\|v+tw\|_{\Lp{{q}}}^{q}\right)\right|_{t=0}\\
&=\left.\frac{\dm}{\dm t}\left(\frac{1}{{q}}\int_{\Omega}|v+tw|^{q} \dm \vc x
\right)\right|_{t=0} \\
&= \int_{\Omega} |v|^{{q}-1} \sign (v) w \dm \vc x.
\end{aligned}\label{eq:dmapLp}
\end{align}
Note that $\varphi^{-1}(t) = t^{{q'}-1}$ for $q$ such that $1=1/q+1/q'$ and
therefore $J_{q'}^{-1} = J_{q}$. Moreover, we have for all $v \in \Lp{q}$
\begin{align}
\|J_{q}(v)\|_{\Lp{{q'}}} = \|v\|_{\Lp{{q}}}^{{q}-1}, \qquad \langle J_{q}(v), v
\rangle_{\Lp{{q'}}, \Lp{{q}}} = \|v\|_{\Lp{{q}}}^{q} .
\end{align}
Similarly, we can compute the normalized duality map $\tilde{J}_{q}$ on
$\Lp{{q}}$, i.e., the duality map with weight $\varphi(t)=t$, and obtain
\begin{align}
\langle \tilde{J}_{q}(v), w \rangle_{\Lp{{q'}}, \Lp{{q}}}
&= \|v\|^{2-{q}}_{\Lp{{q}}}\int_{\Omega} |v|^{{q}-1} \sign (v) w \dm \vc x.
\end{align}
Comparing the expressions for $J_{q}$ and $\tilde{J}_{q}$, we can see that it
may be useful, in particular for the implementation, to use $J_{q}$ instead of
the normalized duality mapping in order to avoid the additional scaling
$\|v\|_{\Lp{{q}}}^{2-{q}}$. Indeed, if $\tilde{J}_q$ is a non-linear operator
in a variational problem that is approximated using a finite element method, we
would need to evaluate the derivative of $\tilde{J}_q$ to determine the Jacobi matrix
used in Newton's method. The derivative is given by
\begin{align*}
  \langle \tilde{J}_q'(v)(z),w\rangle_{\Lp{q'}, \Lp{q}} &= (q-1)\|v\|_{\Lp{q}}^{2-q}\int_{\Omega}|v|^{q-2}zw\dm \bm x \\%
&\phantom{=}+(2-q)\|v\|^{2-2q}_{\Lp{q}}\int_{\Omega}|v|^{q-2}v z \dm \bm x \int_{\Omega}|v|^{q-2}v w \dm \bm x.
\end{align*}
If $w$ and $z$ are finite element functions with local support,
the first term is only non-zero if both $w$ and $z$ are non-zero and
thus this term would yield a sparse matrix for typical finite element spaces.
The second term on the other hand is always non-zero if $v$ is non-zero on the
whole domain and thus may lead to a dense Jacobi matrix. If we consider $J_q$ instead,
the derivative consists of only the first term without the scaling $\|v\|_{\Lp{q}}^{2-q}$
and therefore we obtain a sparse Jacobi matrix.

As a second example consider the space $\sobzero{1}{{q}}$ with the (semi-)norm
$|v|_{\sob{1}{q}}$. In the same way as before we can compute the
duality map of weight $\varphi(t)=t^{{q}-1}$ and obtain
\begin{align}
\langle J_{\sobzero{1}{{q}}}(v), w \rangle_{\Lp{{q'}}, \Lp{{q}}}
= \int_{\Omega} \sum_{i=1}^d|\partial_i v|^{q-1}\sign (\partial_i v) \partial_i w\dm \vc x.
\label{eq:dmapW1p}
\end{align}

\begin{rem}[q-Laplacian]
  Let $1\leq q<r<\infty$ and $x\in\Rd$. In this case we have the following norm
  equivalence on the finite dimensional space $\Rd$:
  \begin{subequations}
  \begin{align}
    \|x\|_{l^q} &= \left(\sum_{i=1}^d |x_i|^q \right)^{\frac{1}{q}} \leq d^{\frac{1}{q}-\frac{1}{r}}
    \left(\sum_{i=1}^d |x_i|^r \right)^{\frac{1}{r}}=d^{\frac{1}{q}-\frac{1}{r}}\|x\|_{l^r},\\
    \|x\|_{l^r} &\leq \|x\|_{l^q}.
  \end{align}
\end{subequations}
As a result
\begin{align}
  \left(\int_{\Omega}\|\nabla v\|_{l^2}^q \dm \vc x\right)^{1/q}
\end{align}
defines a norm on $\sobzero{1}{q}$ that is equivalent to $|v|_{\sob{1}{q}}$ and the
corresponding duality mapping for $q>2$ is given by
\begin{align}
  \int_{\Omega} \|\nabla v\|_{l^2}^{q-2}\nabla v \cdot \nabla w \dm \vc x,
\end{align}
which is also the weak form of the $q$-Laplacian.
\end{rem}
\section{Minimum Residual Methods in Banach Spaces}\label{sec:minres}

As mentioned before, the method we are considering is a generalization of the scheme
discussed in \cite{Chan2014}. Conceptually, this method is closely related to
DPG methods with the main difference that only continuous finite element spaces are
considered in \cite{Chan2014}. The key observation that makes it possible to extend
the methodology to general Banach spaces is the equivalence with a residual minimization
problem, cf.,  \cite{Chan2014} and \cite{Demkowicz2015} for the DPG method.
  In this section, we will see that, in the context of
Banach spaces, the Riesz map can be replaced by  a duality mapping, but due to the
non-linearity of the duality mapping we lose the concept of an optimal test space which
is used in the context of DPG methods, cf., \cite{Demkowicz2011}. The extension to Banach spaces was
first introduced in \cite{Muga2015} and we will repeat the main concepts in this section.

Let $U$ and $V$ be two Banach spaces and $b:U\times V\rightarrow \R$ a continuous bilinear form
that satisfies the following inf-sup conditions:
\begin{subequations}
\begin{align}
\inf_{w \in U} \sup_{v\in V} \frac{b(w,v)}{\|w\|_U \|v\|_V} = \gamma >0, \label{eq:infsup}\\
\left\{v\in V \, : \, b(w,v) = 0, \forall w \in U\right\} = \{0\}.\label{eq:infsup2}
\end{align}
\end{subequations}
For a given right-hand side $\ell \in V'$, we consider the problem: find $u \in U$ such that
\begin{align}
b(u,v) = \ell(v) \quad \forall v \in V.
\end{align}
Introducing a finite dimensional subspace $U_n \subset U$, we can formulate the
following residual minimization problem: find $u_n \in U_n$ such that
 \begin{align}
   \begin{split}
 u_n &= \argmin{w_n \in U_n} H_V(w_n), \\
  H_V(w_n)&:= \psi(\|Bw_n-\ell\|_{\dual{V}}),
\end{split}\label{eq:res_min_Banach}
 \end{align}
 where $\psi(0) = 0$, $\varphi :=\psi'$ is a weight function as defined in Section \ref{sec:duality_map} and $B:U\rightarrow V'$ denotes the linear operator associated with the bilinear form $b$.
 A typical choice for $\psi$ would  be $\psi(t) = t^{q}/q$, $1 < q < \infty$. If $\dual{V}$
 is strictly convex, then the duality mapping $J_{\dual{V}}^{\varphi}$ is single valued
 and $J_{\dual{V}}^{\varphi} = \nabla (v' \mapsto \psi(\|v'\|_{\dual{V}}))$.

 Note that
 this best approximation problem can be formulated in any Banach space. However, we
 need strict convexity for uniqueness of minimizers, cf., e.g., \cite{Stakgold2011,Muga2015}.
 The closed range theorem is still applicable in general Banach spaces
  and provides conditions for well-posedness of linear problems. Its reformulation
  in terms of inf-sup conditions only requires reflexivity of the test space to
  describe the kernel of $B'$ by \eqref{eq:infsup2}. Here, $B': V \rightarrow U'$ denotes the adjoint operator to $B$.

\subsection{Saddle Point Formulation}
If $\dual{V}$ is strictly convex, it can be shown that the residual minimization problem is equivalent to $\nabla H_V(u_n) = J_{V'}^{\varphi}(\ell-Bu_n) = 0$. If $V$ is additionally reflexive, as well as being
 strictly convex, we obtain the following non-linear Petrov Galerkin formulation:
 find $u_n \in U_n$ such that
 \begin{align}
   \begin{split}
   \langle J_{\dual{V}}^{\varphi}(\ell-B u_n), &{B w_n}\rangle_{V'',V'}\\
   &=\left\langle{B w_n},{\left(J_{{V}}^{\varphi^{-1}}\right)^{-1}(\ell-B u_n)}\right\rangle_{V',V}  =0 \qquad \forall w_n \in U_n. \label{eq:NPGform}
\end{split}
 \end{align}
Here, we used $J^{\varphi}_{V'} = \left(J_{{V}}^{\varphi^{-1}}\right)^{-1}$ by means of canonical identification, cf., Proposition \ref{prop:duality_map_convex_spaces}. Introducing an auxiliary variable $r = \left(J_{V}^{\varphi^{-1}}\right)^{-1}(\ell - Bu_n)$,
 this can be reformulated as a mixed method: find $u_n \in U_n$ and $r \in V$ such that
 \begin{subequations}
 \begin{align}
   \langle J_V^{\varphi^{-1}}(r), v \rangle_{\dual{V}, V} + \langle B u_n, v \rangle_{\dual{V}, V}
   &= \langle \ell, v \rangle_{\dual{V}, V} && \forall v \in V,\\
   \langle B w_n , r \rangle_{\dual{V}, V} & = 0 && \forall w_n \in U_n.
 \end{align}\label{eq:sys:mixed_banach}
 \end{subequations}
 Note that $J_V^{\varphi^{-1}}$ is a non-linear map unless $V$ is a Hilbert space.
 If $V$ is a Hilbert space and we select $\psi(t) = t^2/2$, then $J_V^{\varphi^{-1}} = R_V$
 and we recover the framework in \cite{Chan2014}.
\subsection{The Optimal Test Norm}
In applications, we are often more interested in minimising the error in the approximation with respect to the norm on $U$ rather than minimising the residual. In other words, we usually want to choose the
norm on $V$ in such a way that we ultimately control the error $\|u-u_n\|_U$.
With this in mind, employing the inf-sup condition \eqref{eq:infsup} and continuity of $B$, we deduce that
\begin{align}
\|u-u_n\|_U \leq \frac{1}{\gamma} \sup_{v\in V}\frac{b(u-u_n, v)}{\|v\|_V}
 = \frac{1}{\gamma}\|\ell - Bu_n\|_{\dual{V}} \leq \frac{M}{\gamma}\|u-u_n\|_U,
\end{align}
where $M$ is the continuity constant of the bilinear form $b$.
Hence, in order to control $\|u-u_n\|_U$, we require $\gamma$ and $M$
as close to one as possible and  independent of certain problem specific
parameters (the scaling of the diffusion term, for example, in case of
the convection-diffusion-reaction equation). The \emph{optimal} test norm is a concept
introduced in the context of DPG methods \cite{Zitelli2011} but unlike the concept of
optimal test functions and spaces it can easily be extended to Banach spaces.
The optimal test norm is defined as the norm on $V$ such that $M = \gamma = 1$ and is given by
\begin{align}
\|v\|_{\mathrm{opt}} := \sup_{u \in U} \frac{b(u,v)}{\|u\|_U} = \|B' v\|_{\dual{U}}.
\end{align}
Indeed,
\begin{align}
\sup_{v \in V} \frac{b(w,v)}{\|v\|_{\mathrm{opt}}} = \sup_{v \in V}\frac{(\dual{B}v)(w)}{\|\dual{B} v\|_{\dual{U}}} = \sup_{g \in \dual{U}} \frac{g(w)}{\|g\|_{\dual{U}}} = \|w\|_U.
\end{align}
Conversely, we have the optimal constants $M=\gamma=1$ for any given test norm,
if we endow $U$ with the so-called energy norm,
\begin{align}
  \|u\|_E := \sup_{v \in V} \frac{b(u,v)}{\|v\|_V}.
\end{align}

Ideally, we want to select the norm on $U$ in order to approximate the solution
in the desired norm and then work with the optimal test-norm. However, the optimal
test-norm is a dual norm which is, in general,  not computable. Thus, we have to replace
the optimal test-norm with an equivalent norm that is computable. The difficulty
here is to obtain equivalence constants, i.e., constants $\gamma$ and $M$, that are independent of problem parameters.

In the context of the convection-diffusion-reaction equation, Broersen and Stevenson
work directly with the optimal test norm in \cite{Broersen2014},
whereas the analysis presented in \cite{Demkowicz2013} and \cite{Chan2014a} relies on
robust estimates for $\gamma$ and $M$ by looking at the adjoint problem.
Extending any of the  robust estimates from Hilbert subspaces of $\Lp{2}$ to Banach spaces
$\Lp{q}$, $1<q<\infty$, is highly non-trivial and to the best of our knowledge remains an
open problem.

\subsection{The Inexact Method}\label{sec:inexact}
So far, \eqref{eq:sys:mixed_banach} does not define a method that can be implemented because the mixed method is still an infinite
dimensional problem that relies on the whole test space $V$.
Replacing the space $V$ in \eqref{eq:sys:mixed_banach} by a finite dimensional $V_m \subset V$,
however, allows us to approximate the solution to the best-approximation problem by
the solution of a finite-dimensional saddle point problem. More precisely, we obtain the
following fully discrete mixed problem: find $(r_m, u_n) \in V_m \times U_n$
such that
\begin{subequations}
\begin{align}
  \langle J_V^{\varphi^{-1}}(r_m), v_m \rangle_{\dual{V}, V} + \langle B u_n, v_m \rangle_{\dual{V}, V}
  &= \langle \ell, v_m \rangle_{\dual{V}, V} && \forall v_m \in V_m,\\
  \langle B w_n , r_m \rangle_{\dual{V}, V} & = 0 && \forall w_n \in U_n,
\end{align}\label{eq:sys:mixed_banach_inexact}
\end{subequations}
where $J_V^{\varphi^{-1}}=R_V$ if $V$ is a Hilbert space.

Replacing $V$ by a finite-dimensional subspace we obviously lose the best approximation
properties.  The question  of quantifying this `loss' was addressed
 in \cite{Gopalakrishnan2014} for Hilbert spaces by introducing the concept of a
Fortin-operator and was later extended to Banach spaces in \cite{Muga2015}. For well-posedness
of the inexact method we require the existence of a bounded projection operator
$\Pi:V \rightarrow V_m$ such that
\begin{align}
  \langle B w_n, v- \Pi v\rangle_{V', V} = 0 \qquad \forall w_n \in U_n, v \in V. \label{eq:Fortin_operator}
\end{align}
The norm of $\Pi$ then enters into the a priori estimate and quantifies the `loss' due to the inexactness, i.e.,
we can obtain an a priori estimate of the form
\begin{align}
  \|u-u_n\| \leq C\inf_{w_n \in U_n} \|u-w_n\|,\qquad C = \frac{(\|\Pi\|+\|I-\Pi\|)M}{\gamma}.
\end{align}
\begin{rem}
The constant $C$ can be improved both in the Hilbert setting, where we can
obtain $C = \|\Pi\|M/\gamma$, cf., \cite{Gopalakrishnan2013}, and in the Banach setting, where we would have to
introduce geometric constants for the improved estimate, cf., \cite{Muga2015}.
\end{rem}

\section{The Convection-Diffusion Equation}\label{sec:confusion}
The main focus of this article is to apply the abstract approach described in the previous section
to the scalar convection-diffusion-reaction equation in a $\sob{1}{q}$-$\sob{1}{q'}$
setting. To this end,
consider the following model problem:
\begin{subequations}
\begin{align}
-\varepsilon \Delta u + \vc{b}\cdot\nabla u +cu&= f \qquad \text{ in }\Omega,\\
u&= 0 \qquad\text{ on } \Gamma = \partial \Omega,
\end{align}\label{eq:confusion}
\end{subequations}
where $\varepsilon$, $\bm b:\Omega \rightarrow \Rd$ and $c:\Omega \rightarrow \R$ are the
(positive) diffusion parameter, convection field and reaction coefficient, respectively,
and $f:\Omega \rightarrow \R$ is a given source term.
Multiplying by a test function $v \in C_c^{\infty}(\Omega)$ and integrating by parts yields the bilinear form
\begin{align}
\mathcal{B}_{\varepsilon}(u,v) &= \varepsilon\int_{\Omega} \nabla u \cdot\nabla v \dm \vc{x} - \int_{\Omega}u \nabla \cdot (\vc{b}v) \dm \vc x + \int_{\Omega} c u v \dm \vc x.\label{eq:conv_diff}
\end{align}
This allows us to state the following variational problem: find $u\in U:=\sobzero{1}{{q}}$ such that
\begin{align}
	\mathcal{B}_{\varepsilon}(u,v) = \langle f, v\rangle_{V', V}\qquad \forall v \in V:=\sobzero{1}{{q'}}.\label{eq:Beps}
\end{align}
In \cite{Houston2019}, a proof of an inf-sup condition is given, where the norms on $U$ and $V$ are weighted versions of the standard norms $\sob{1}{q}$ and $\sob{1}{q'}$, respectively. Furthermore, this proof requires certain regularity assumptions on the solution to the Poisson problem and that the convection field $\bm b$ and the reaction coefficient $c$ satisfy $c- \nabla \cdot \frac{1}{q}\bm b \geq c_0>0$.
The continuity constant and the inf-sup constant that are established for the bilinear $\mathcal{B}_{\varepsilon}$ in \cite{Houston2019} depend on the problem specific parameters. It should be noted that the estimates
in \cite{Houston2019} can be expected to be sub-optimal since sharper bounds are
known for $q = 2$.
\subsection{Choices for the Test Norm}\label{sec:test_norms}
The choice of the norm for the space $V$ crucially defines the method described in
Section \ref{sec:inexact}. To this end, we endow $V$ with the following norm
\begin{align}
\|v\|_V^{q'} &:= \alpha\|v\|_{\Lp{{q'}}}^{q'}+\varepsilon \|\nabla v\|_{\Lp{q'}}^{q'}+\frac{|\Omega|^{1/2}}{\|\bm b\|_{\Lp{\infty}}}\|(\omega(x))^{1/q'}\bm b \cdot \nabla v\|_{\Lp{q'}}^{q'}, \label{eq:norm_eps1q}
\end{align}
where $\omega(x)$ is a positive and smooth weighting function.  It is well
known that for $q=2$ and $\omega(x) \equiv \mathrm{const} \geq 0$  the above choice of norm on $V$ does not yield a robust
formulation, i.e., $M/\gamma$ depends on the problem parameters. Here, $M$ denotes the continuity constant of the bilinear form
and $\gamma$ its inf-sup constant.
The inf-sup constant for $1<q<\infty$ obtained in \cite{Houston2019} corresponds to the choice $\omega(x) \equiv 0$ and $\alpha = c_0$, where $c_0$ is the constant in the Friedrich's positivity assumption, $c- \frac{1}{q}\nabla \cdot  \bm b \geq c_0$. As mentioned above, the continuity and inf-sup constants are not parameter independent but can be expected to be suboptimal.

One approach to address the robustness
issue is to introduce a non-constant weighting function, cf., \cite{Demkowicz2013,Chan2014}.
For example, we may require $\omega(x)$ to be of magnitude $\mathcal{O}(\varepsilon)$ near the inflow
boundary but $\mathcal{O}(1)$ elsewhere. We will show the effect
of this by considering a simple one dimensional example on $\Omega = (0,1)$ with the
inflow boundary at $0$. For this example $\omega(x) = x+ \varepsilon$ has the desired
properties; for comparison, we will also consider $\omega(x)\equiv 0$ and
$\omega(x) \equiv 1$.

\subsection{Weak Boundary Conditions on $r$}\label{sec:weak_inflow_bc}
Alternatively, the robustness issue can be addressed by changing the boundary conditions  on $r$; this has been considered
in different ways in \cite{Cohen2012,Chan2014}. In \cite{Chan2014a}, the boundary conditions on $u$ were modified instead. We will now present the approach in \cite{Chan2014} and extend it for $1<q<2$.
The idea is to relax the boundary condition on the test space on the inflow part of the boundary.
The reasoning behind this is that in the mixed method we essentially approximate the
adjoint equation in the test space in order to approximate the residual or the optimal
test functions. Under resolved layers at the inflow boundary ---which is the
outflow boundary for the adjoint equation--- then pollute the solution to the
primal problem in the inexact method.

Instead of $\sobzero{1}{q'}$, we consider the modified test space $V=W^{1,q'}_{0, \Gamma_+}(\Omega)$, where $q' = q/(q-1)$,
 and
\begin{align}
W^{1,q'}_{0, \Gamma_+}(\Omega) &:= \left\{v \in \sob{1}{q'} \,: \, v\big|_{\Gamma_+} = 0\right\},
\end{align}
Here,
\begin{align*}
  \Gamma_- := \{ x \in \partial \Omega \,:\, \bm b \cdot \bm n(x) \leq 0\}, && \Gamma_+ := \partial \Omega \setminus \Gamma_-,
\end{align*}
where $\bm n(x)$ denotes the unit outward normal at a point on the boundary $\partial \Omega$.
The modified bilinear form (cf., \cite{Chan2014}) is given by
\begin{align}
\tilde{\mathcal{B}}_{\varepsilon}(u,v) := \varepsilon \int_{\Omega} \nabla u \cdot \nabla v \dm \vc x
 + \int_{\Omega} (\vc b \cdot \nabla u) v \dm \vc x
  - \varepsilon \int_{\Gamma_-} \frac{\partial u}{\partial \vc n} v \dm s.\label{eq:weak_inflow_bc}
\end{align}
The boundary term is merely the term that is picked up from the integration by parts if $v$ is non-zero on the inflow boundary.
Note, however, that the term $\varepsilon \int_{\Gamma_-} \frac{\partial u}{\partial \vc n} v \dm s$ is
a variational crime if we assume $u \in \sobzero{1}{q}$. As noted in \cite{Chan2014}, the
correct way of including this term would be to introduce it as an additional unknown on the boundary
and the inclusion of $\varepsilon \int_{\Gamma_-} \frac{\partial u}{\partial \vc n} v \dm s$ can
be viewed as a discrete elimination of this unknown, cf., \cite{Chan2014,Broersen2014}.

\subsection{The Inexact Method for the Con\-vec\-tion-Dif\-fusion-Re\-ac\-tion Equation}\label{sec:method_confusion}
The key step for implementing the inexact method \eqref{eq:sys:mixed_banach_inexact}
for any specific problem is determining the duality mapping $J_V^{\varphi^{-1}}$.
We choose the weight function $\varphi^{-1}(t) = t^{q'-1}$ and compute the duality
mapping, similar to Section \ref{sec:duality_map_example}; thereby, we get
\begin{align}
\begin{aligned}
\langle J_V^{\varphi^{-1}} (v), w\rangle_{\dual{V},V} = &\phantom{+}\phantom{\varepsilon}\int_{\Omega}|v|^{q'-1} \mathrm{sgn}(v) w \dm \vc x\\&+\varepsilon\int_{\Omega}\sum_{i=1}^d|\partial_i v|^{q'-1} \mathrm{sgn}(\partial_i v)  \partial_i w \dm \vc x\\
&+\int_{\Omega}\omega(x)|\vc b \cdot \nabla v|^{q'-1} \mathrm{sgn}(\vc b \cdot \nabla v) \vc b \cdot \nabla  w\dm \vc x.
\end{aligned}
\end{align}

For a given  right hand side $\ell \in \dual{V}$, we solve the following
non-linear system: find $(u_n, r_m) \in U_n \times V_m \subset U\times V$ such that
\begin{subequations}
\begin{align}
\dualp{J_V^{\varphi^{-1}}(r_m)}{v_m}{V} + {\mathcal{B}}_{\varepsilon}(u_n, v_m) &=\ell(v) && \text{ for all } v_m \in V_m,\\
{\mathcal{B}}_{\varepsilon}(w_n, r_m) & = 0 && \text{ for all }w_n \in U_n.
\end{align}\label{eq:inexact_confusion}
\end{subequations}
 One can easily implement \eqref{eq:inexact_confusion} in, e.g.,   FEniCS \cite{AlnaesBlechta2015a,LoggMardalEtAl2012a}
using standard $\hone$-con\-for\-ming Lagrange finite elements for both $U_n$ and $V_m$. The spaces $U_n$
and $V_m$ are chosen over a common mesh. For $U_n$ a global polynomial degree $
p_n$ is chosen and for $V_m$ we choose an enriched finite element space with
polynomial degree $p_m = p_n + \varDelta p$, $\varDelta p \geq 1$. To consider the weak boundary conditions
on $r$ introduced in the previous section, we simply replace $\mathcal{B}_{\varepsilon}$
with $\tilde{\mathcal{B}}_{\varepsilon}$ and adjust the space $V_m$ to only
satisfy Dirichlet boundary conditions on the outflow boundary.

\subsection{The Limit Case $\varepsilon = 0$}
Since we are interested in the convection-dominated case, i.e., $\varepsilon \ll \|\bm b \|_{\Lp{\infty}}$, it makes sense to consider the limit $\varepsilon \rightarrow 0$. Simply setting $\varepsilon = 0$ in \eqref{eq:conv_diff} does not yield a well-posed problem unless we only consider boundary conditions on the inflow boundary $\Gamma_-$. However, even the ill-posed
problem with Dirichlet boundary conditions on the whole boundary $\partial \Omega$ is of interest, since
this is essentially the problem that is approximated numerically on coarse meshes when $\varepsilon \ll \|\bm b \|_{\Lp{\infty}}$, cf., the discussion in \cite{Guermond2004}. On a discrete level we can think of imposing boundary conditions for the approximation $u_n$ on $\Gamma_+$ simply as considering an approximation problem in a smaller subspace of $U$.

When $\varepsilon = 0$, there are two different weak formulations of the \emph{convection-reaction} equation  depending on whether $\int_{\varOmega}(\bm b \cdot \nabla u)v \dm \bm x$ is integrated by parts or not. The two cases differ in the regularity of the trial and test spaces which is reflected in the norms chosen on $U$ and $V$.  Thus, in each case the residual is measured in a different norm and (quasi-)best approximation of the analytical solution is achieved in a different space. Both weak formulations can be extended to  weak formulations of the convection-diffusion-equation by adding the diffusion term and accounting for different boundary conditions and regularity requirements. If $\varepsilon >0$, we require $u \in \sobzero{1}{q}$ in both cases and the two weak formulations are formally equivalent.
 However, the differences observed in the limit case $\varepsilon = 0$ can be reflected in the choice of the norm on $\sobzero{1}{q}$ by choosing the weighting of the terms $\| w\|_{\Lp{q}}$, $\|\nabla w\|_{\Lp{q}}$ and $\|\nabla \cdot \bm b w \|_{\Lp{q}}$ accordingly.

 We show that for the first choice the exact mixed method is equivalent to the formulation used in \cite{Guermond2004} and that in the second case a quasi-best approximation in $\Lp{q}$ can be computed by determining the optimal test norm. We use this to interpret certain choices of the test norm for the convection-diffusion-reaction equation.

\subsubsection{Residual Minimisation in $\Lp{q}$}
Both weak formulations for the convection-reaction equation are obtained by multiplying the convection-diffusion-reaction equation \eqref{eq:confusion} with $\varepsilon = 0$ by a smooth test function and integrating over the domain $\Omega$. This immediately yields the first possible choice for the weak formulation: find $u \in W_{0, \Gamma_-}^q(\bm b, \varOmega)$ such that
\begin{align}
\hat{\mathcal{B}}_0(u,v) = \int_{\Omega}(\bm b\cdot \nabla u)v \dm \bm x + \int_{\Omega} cuv \dm \bm x &&\forall v \in \Lp{q'}.
\end{align}
Here, the bilinear form $\hat{\mathcal{B}}_0$ is well defined for $v \in \Lp{q'}$ and  $u$ in the graph space
\begin{align}
W^q_{0, \Gamma_-} (\bm b, \Omega):= \{ w \in \Lp{q} \,:\, \bm b \cdot \nabla u \in \Lp{q} \text{ and } u = 0 \text{ on } \Gamma_-\}
\end{align}
endowed with the norm
\begin{align}
\| u\|_{W^q_{0, \Gamma_-} (\bm b, \Omega)}^q = \|u\|_{\Lp{q}}^q+\|\bm b \cdot \nabla u\|_{\Lp{q}}^q.
\end{align}
The associated residual minimisation problem is
\begin{align}
u_n = \argmin{w_n \in U_n} \|\hat{\mathcal{B}}_0(w_n , \cdot) - f\|_{\Lp{q}};
\end{align}
this was considered in \cite{Guermond2004}.
In general,  method \eqref{eq:sys:mixed_banach_inexact} yields a formulation for
solving the residual minimisation problem inexactly. However, since in this case
$J^{\varphi}_{V'} = J_q$, we can avoid the inexactness by directly
implementing \eqref{eq:NPGform}; this is exactly the approach considered in \cite{Guermond2004} where
a regularisation is introduced to compute \eqref{eq:NPGform}.
The norm on  $V = \Lp{q'}$  corresponds to the choice of the test-norm for the convection-diffusion-reaction equation given in \eqref{eq:norm_eps1q} with $\omega(x) \equiv 0$ and $\alpha =1$ since this yields the $\Lp{q'}$-norm if $\varepsilon$ is set to zero. This suggests, that this choice of $\omega(x)$ and $\alpha$ corresponds to the natural extension of this bilinear form to the case $\varepsilon > 0$.
Moreover, the equivalence of the exact minimum residual method and the approach in \cite{Guermond2004} for $\varepsilon = 0$  together with $\varepsilon$ weighting of the norm of the gradient suggest that for this choice of the test norm, our scheme  closely resembles the approach in \cite{Guermond2004} for $0 < \varepsilon \ll 1$.

\subsubsection{Best-$L^q$ Approximation}
 Next, we can integrate the convection term by parts in order to obtain the second weak formulation. To this end, we note that
\begin{align}
\int_{\Omega} (\bm b \cdot \nabla u) v\dm \bm x = -\int_{\Omega} \nabla \cdot (\bm b v) u \dm \bm x+\int_{\partial \Omega} (\bm b \cdot \bm n) uv \dm \bm s.
\end{align}
The boundary term on $\Gamma_-$ vanishes since $u = 0$ on $\Gamma_-$ (for non-zero boundary conditions this term would be absorbed into the right hand side by inserting the boundary condition). For $v$ in the graph space
\begin{align}
W^{q'}_{0, \Gamma_+} (\bm b, \Omega):= \{ w \in \Lp{q'} \,:\, \bm b \cdot \nabla w \in \Lp{q'} \text{ and } w = 0 \text{ on } \Gamma_+\},
\end{align}
the boundary term also vanishes on $\Gamma_+$ and we obtain the bilinear form
\begin{align}
\mathcal{B}_0(u,v) = - \int_{\Omega}u \nabla \cdot (\vc{b}v) \dm \vc x + \int_{\Omega} c u v \dm \vc x,
\end{align}
which is well-defined for $u \in \Lp{q}$. In this case we can compute the optimal test norm
\begin{align}
\|v\|_{\mathrm{opt}} := \sup_{u \in \Lp{q}}\frac{\mathcal{B}_0(u,v)}{\|u\|_{\Lp{q}}} = \|\mathcal{B}_0(\cdot, v)\|_{\left(\Lp{q}\right)'} = \|-\nabla\cdot(\bm b v)+cv\|_{\Lp{q'}}.
\end{align}
Choosing the optimal test norm, allows us to obtain the $L^q$-best approximation in a given finite
dimensional space $U_n$. It is easy to see that the optimal test norm is bounded from above up to a constant by the graph norm
\begin{align*}
  \|w\|_{W^{q'}_{0,\Gamma_+}}^{q'}:= \| w\|_{\Lp{q'}}^{q'} + \|\bm b \cdot \nabla w\|_{\Lp{q'}}^{q'}.
\end{align*}
In \cite{MugTylZee2018} this formulation of the convection-reaction equation is analysed in detail and an inf-sup condition is established assuming $c - \frac{1}{q}\nabla \cdot \bm b \geq c_0 >0$. This implies that the optimal test norm is up to a constant also bounded from below by the graph norm. In other words, the graph norm is equivalent to the optimal test norm. The graph norm can be obtained from \eqref{eq:norm_eps1q} by choosing $\omega(x) \equiv 1$, $\alpha = 1$ and setting $\varepsilon$ to zero. This suggests that this choice of $\alpha$ and $\omega(x)$ yields the natural extension of this formulation to the convection-diffusion-reaction equation.
The dependence of the inf-sup constant in \cite{MugTylZee2018} on $c_0$ suggests
that the equivalence constants can be improved by choosing $\alpha = c_0$.

\section{Numerical Examples}\label{sec:numerics}
In this section we consider a range of numerical test cases to illustrate the performance of our proposed numerical scheme. To this end, in Section \ref{sec:convergence}  we demonstrate that in the diffusion-dominated regime optimal convergence rates are achieved and moreover that in the convection-dominated regime the convergence rate is as expected. We also study the effect of different choices of the discrete test space $V_m$ on the convergence rates. In Section \ref{sec:1D_numerics}, we then compare different choices for the test norm on $V$ and the boundary conditions as described in Sections \ref{sec:weak_inflow_bc}  and \ref{sec:test_norms} for a one-dimensional example.
 Furthermore, for the same example we study for which choices of $V$, $\|\cdot\|_{V}$ and $V_m$ we can observe vanishing over- and undershoots as $q \rightarrow 1$ and whether the method is robust in $\varepsilon$.
In Section \ref{sec:numerics_2D} we consider three two-dimensional examples, where for certain selected meshes the overshoots disappear as $q \rightarrow 1$; this behaviour can be predicted by considering $L^1$-best approximations of discontinuities.

The four examples we are using in this section are given below. They consist of one simple one-dimensional example and three two-dimensional examples each with solutions containing boundary layers for small $\varepsilon$. The solution to the last example additionally contains an interior layer.
\begin{ex}\label{ex:1dconfusion}
\begin{align}
- \varepsilon u''+ u' = 0 \quad \text{ in }(0,1), \qquad u(0) = 0, \qquad u(1)=1.
\end{align}
The analytical solution is given by
\begin{align}
u (x) =
\frac{\exp\left({-\frac{1}{\varepsilon}}\right)-\exp\left({\frac{x-1}{\varepsilon}}\right)}{\exp\left(-\frac{1}{\varepsilon}\right)-1}.
\end{align}
\end{ex}
\begin{ex}[Eriksson-Johnson model problem]\label{ex:eriksson-johnson}
\begin{align}
\frac{\partial u}{\partial x}-\varepsilon\left(\frac{\partial^2 u}{\partial x^2}+\frac{\partial^2 u}{\partial y^2}\right) = 0 \quad \text{ in } (0,1)^2,\\ u = 0 \text{ if } x=1, y=0,1,\qquad
u = \sin(\pi y) \text{ if } x=0.
\end{align}
The analytical solution is given by
\begin{align}
  u(x,y) = \frac{\exp(r_1(x-1))- \exp(r_2(x-1))}{\exp(r_1)-\exp(r_2)}\sin(\pi y),
\end{align}
where
\begin{align*}
  r_1 = \frac{1+\sqrt{1+4\pi^2\varepsilon^2}}{2\varepsilon}, && r_2 = \frac{1-\sqrt{1+4\pi^2\varepsilon^2}}{2\varepsilon}.
\end{align*}
\end{ex}
\begin{ex}[Boundary Layer in the Corner of the Domain]\label{ex:skew}
\begin{align}
\bm{b} \cdot \nabla u-\varepsilon \Delta u &= b_2h_1(x)+b_1h_2(y) \quad \text{ in } (0,1)^2,\\ u &= 0 \quad \text{ on }\partial \Omega,
\end{align}
with $\bm b = (b_1, b_2)^T$, $b_{1,2} >0$ and
\begin{align*}
  h_1(x) = x-\frac{1-\exp\left(\frac{b_1x}{\varepsilon}\right)}{1- \exp\left(\frac{b_1}{\varepsilon}\right)}, && h_2(y) = y-\frac{1-\exp\left(\frac{b_2y}{\varepsilon}\right)}{1- \exp\left(\frac{b_2}{\varepsilon}\right)}.
\end{align*}
The analytical solution is given by
\begin{align}
  u(x,y) = h_1(x)h_2(y).
\end{align}
Since $b_{1,2}>0$, the outflow boundary is defined by the two lines $x=1$ and $y=1$. For small $\varepsilon$ we can observe a boundary layer at the outflow boundary
and in particular near the corner $(x,y) = (1,1)$.
\end{ex}
\begin{ex}[Interior and Boundary Layer]\label{ex:interior_layer}
\begin{align}
  \begin{aligned}
\bm{b} \cdot \nabla u-\varepsilon \Delta u &=0 \quad \text{ in } (0,1)^2,\\  u &= 1 \text{ on }\partial \Omega\cap \{ x = 0\}, \\ u &= 0 \text{ on }\partial \Omega \setminus \{ x = 0\},
\end{aligned}
\end{align}
with $\bm b = (2, 1)^T$.

For this example, a boundary layer develops at
$\Gamma_+ \cap \{ y>0.5\}$ and an interior layer along the line $ y = 0.5 x$.
\end{ex}
\subsection{Convergence Tests}\label{sec:convergence}
\begin{figure}[t]
\centering
\includegraphics[width = 9.5cm]{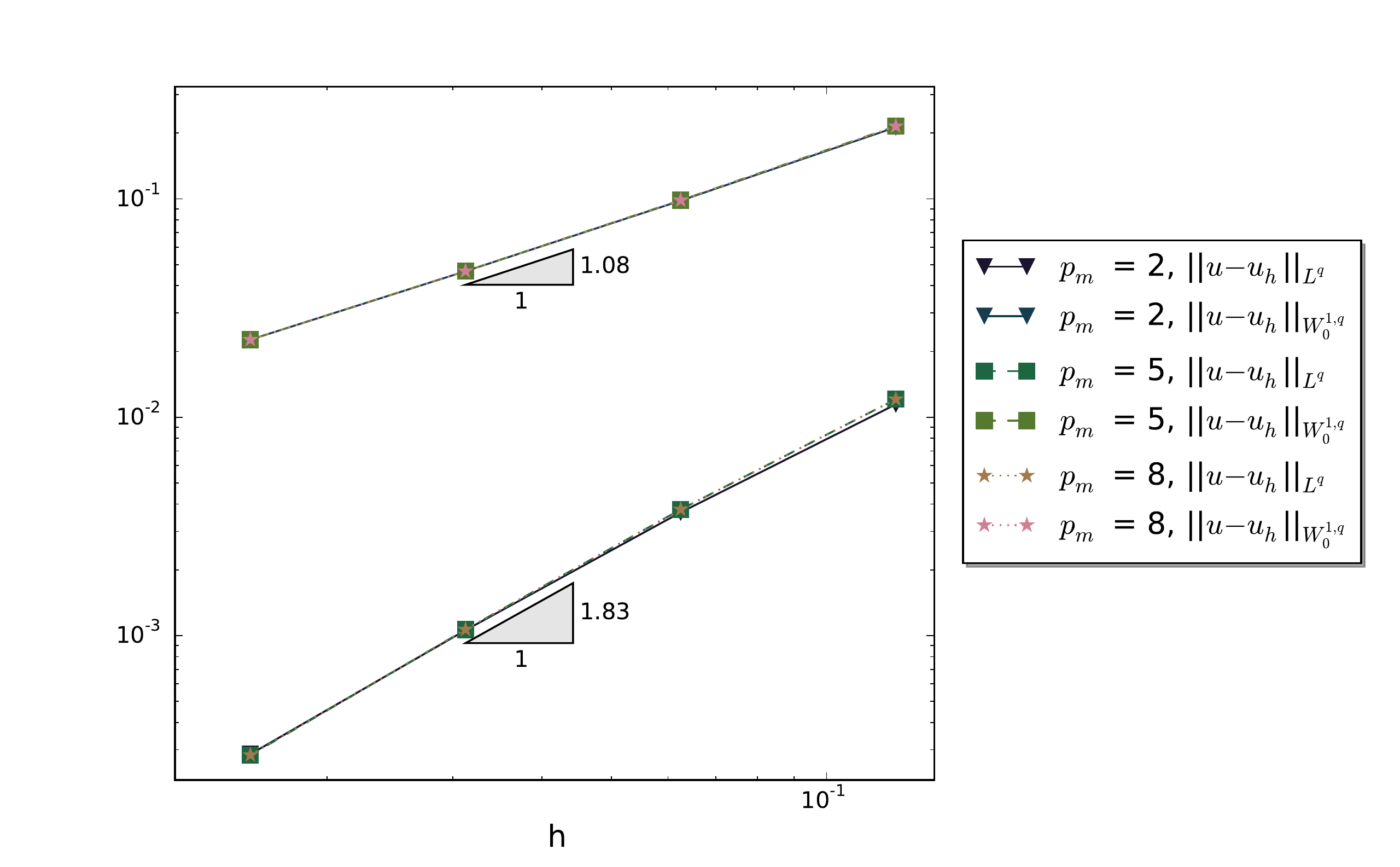}
\caption{Example \ref{ex:eriksson-johnson} with $\varepsilon = 1$ and $q=1.2$. Convergence
for $p_n = 1$ and $\varDelta p = 1,4,7$. }
\label{fig:ej_conv_enrichment}
\end{figure}

\begin{figure}[t]
\centering
\includegraphics[width = 6cm]{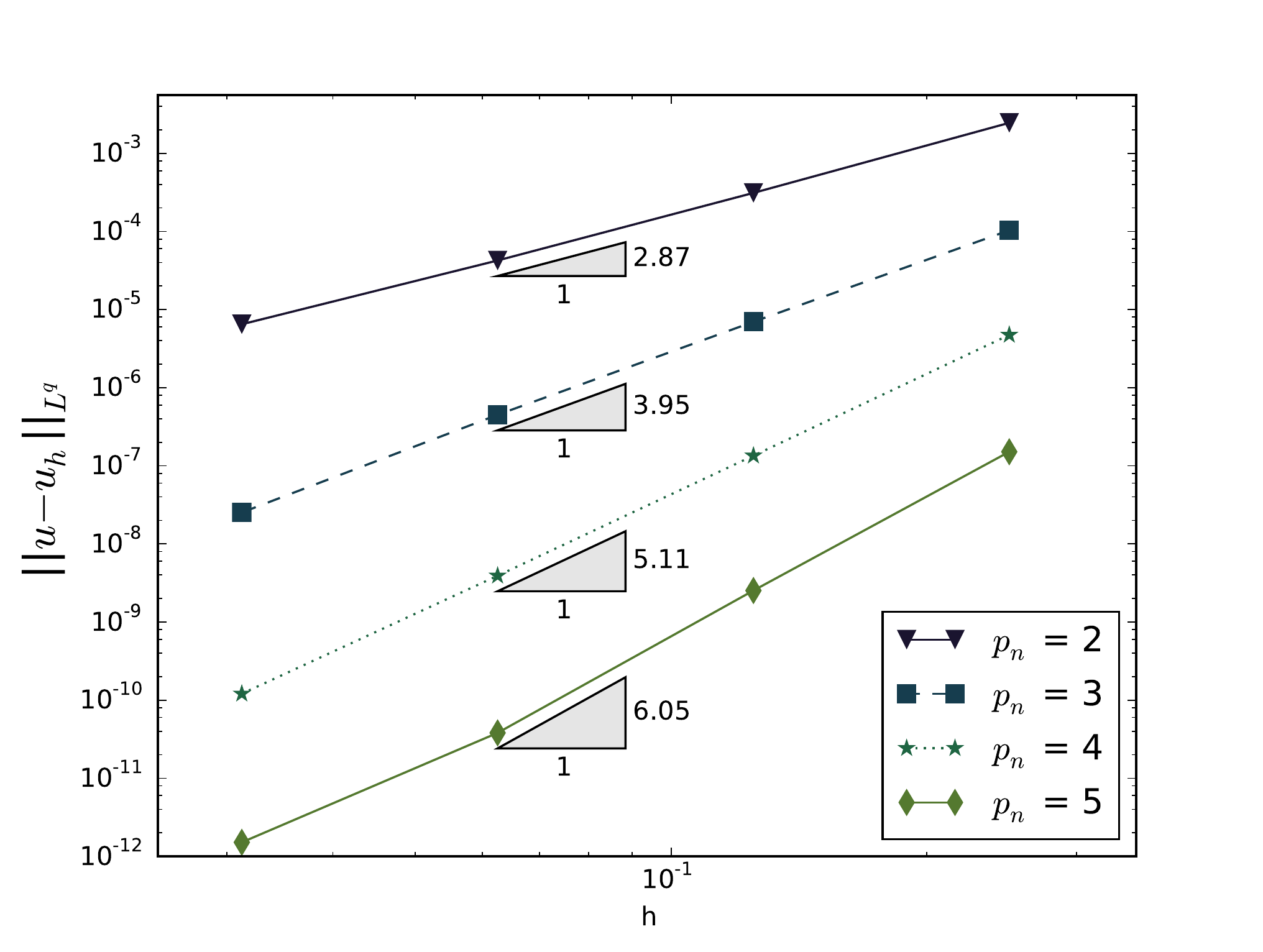}
\includegraphics[width = 6cm]{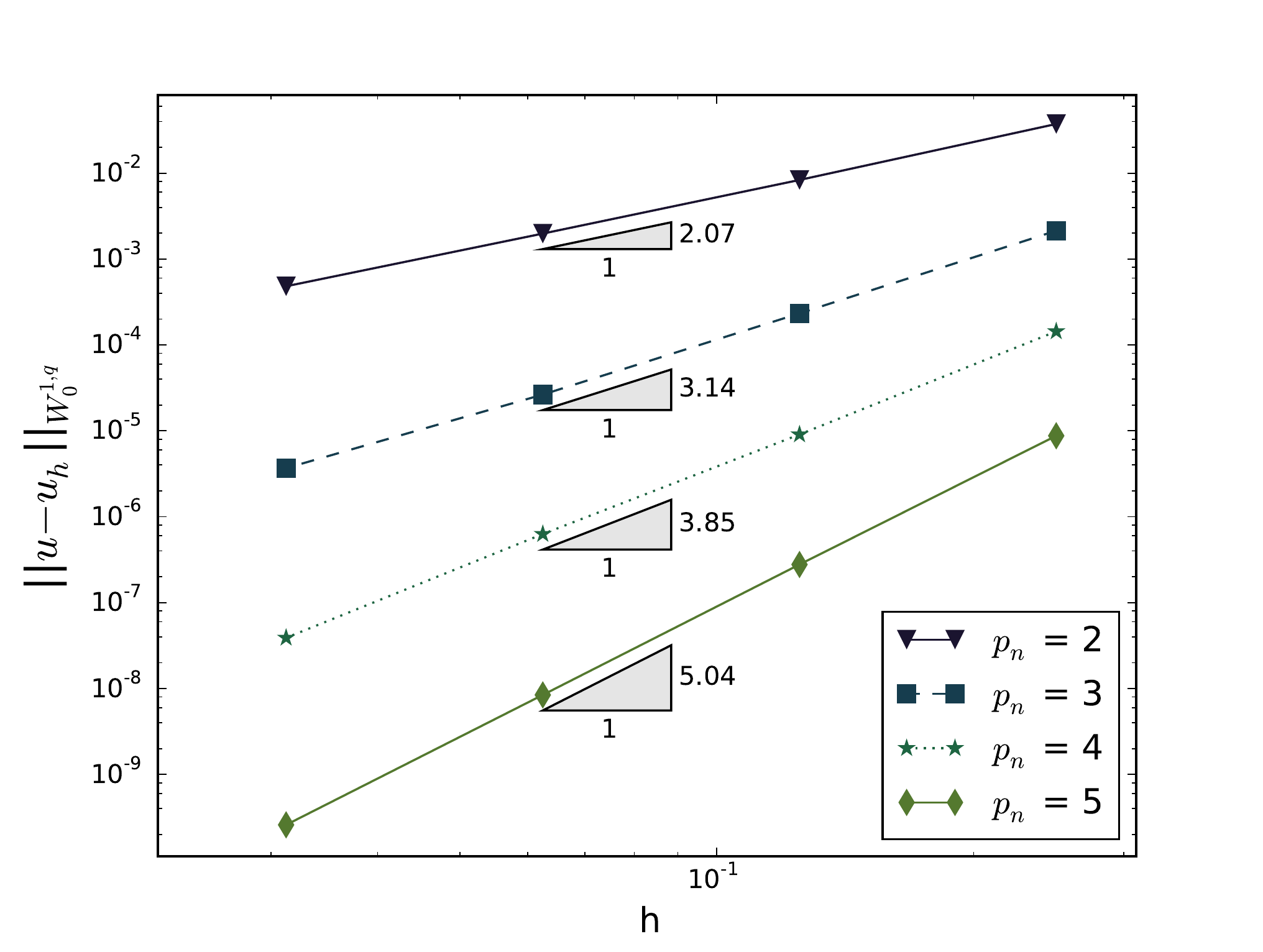}
\caption{Example \ref{ex:eriksson-johnson} with $\varepsilon = 1$ and $q=1.2$,
for $p_n = 2,3,4,5$ and $\varDelta p = 2$. Left: $L^q(\Omega)$-norm, Right: $\sob{1}{q}$-norm}
\label{fig:ej_conv}
\end{figure}

We start with investigating the convergence of the proposed method with weak boundary conditions on the inflow boundary
in the space $V_m$, i.e., formulation \eqref{eq:Beps} with $\tilde{\mathcal{B}}_{\varepsilon}$ from \eqref{eq:weak_inflow_bc}, and the norm \eqref{eq:norm_eps1q} with $\omega(x) \equiv 1$ and $\alpha = 1$. To this end,
we first consider the convergence of the method for Example \ref{ex:eriksson-johnson}
with $\varepsilon = 1$. We consider a uniform mesh with element size $h$. The
space $U_n$ consists of piecewise polynomials of degree $p_n$, while space $V_m$ consists of piecewise polynomials of degree $p_m = p_n +\varDelta p$ on the same mesh with
$\varDelta p \geq 1$. Figure \ref{fig:ej_conv_enrichment} shows that  for $q = 1.2$ the error is essentially independent
of the choice of $\varDelta p$ and hence there is no benefit to enriching the
space of the residual $r_m$ beyond $p_m = p_n +1$. Figure \ref{fig:ej_conv} shows
that we obtain the expected optimal convergence rates for different polynomial degrees
$p_n$ under uniform $h$-refinement.

\begin{figure}
\centering
\includegraphics[width = 6cm]{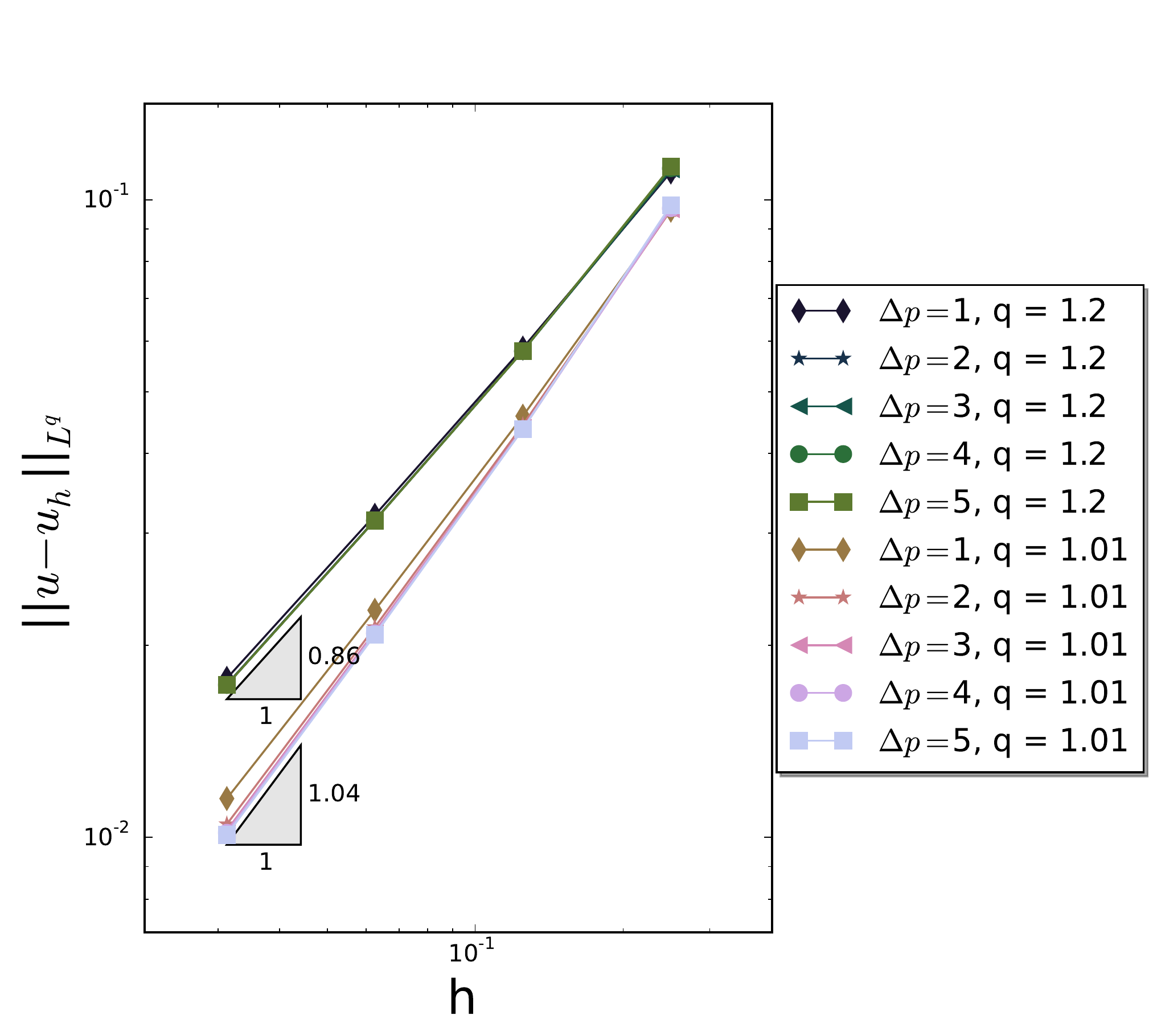}
\includegraphics[width = 6cm]{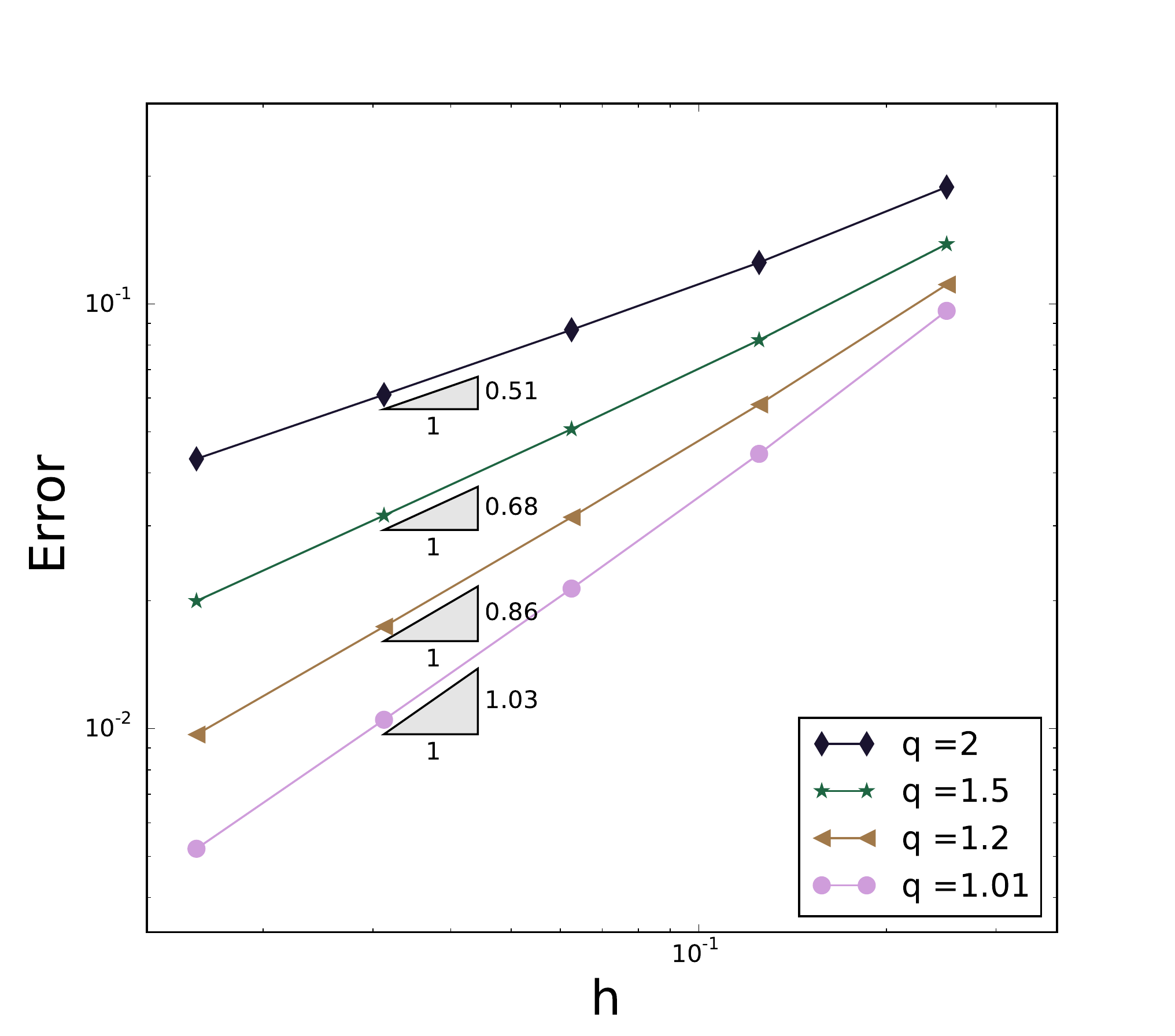}
\caption{$\Lp{q}$-error for Example \ref{ex:eriksson-johnson} with $\varepsilon = 10^{-4}$ and $p_n = 1$.}
\label{fig:ej_conv_eps4}
\end{figure}

Next, we consider the same example with the same finite element spaces for $\varepsilon = 10^{-4}$. Figure
\ref{fig:ej_conv_eps4} (left) shows that again both for $q = 1.01$ and $q = 1.2$ the error does not improve as $\varDelta p$ is
increased beyond $\varDelta p =2$; for $q=1.01$, however, we can  observe that the error is larger for $\varDelta p = 1$. Note that we did not observe this in the diffusive regime. The plot on the right shows that in the convection-dominated regime we obtain a convergence rate of approximately
$\mathcal{O}(h^{\frac{1}{q}})$ as $h$ tends to zero; this is consistent with the approximation error bound of the  piecewise linear interpolant of a jump discontinuity.
\subsection{Convection-Diffusion in 1D}\label{sec:1D_numerics}
In this section we consider different choices for $V$, $\|\cdot \|_V$ and $V_m$ and show how this affects the
approximation to the solution $u$ of Example \ref{ex:1dconfusion}. To this end, in Section \ref{sec:choice_of_norm},
we investigate how the choice of $\alpha$, $\omega(x)$ and the boundary conditions in $V_m$ affect
the approximation $u_n$. In Section \ref{sec:vanishing_osc} we consider one specific choice for the test norm and show that the oscillations vanish entirely as $q \rightarrow 1$. Next,
in Section \ref{sec:choice_Vm}, we investigate the choice of the finite element space
$V_m$ more closely. We have seen that increasing $\varDelta p$ has no significant effect on the convergence rates; however, we will see that this does indeed affect certain qualitative properties of the solution. Finally, we demonstrate that the method is robust in $\varepsilon$ in Section \ref{sec:eps_robust}

\subsubsection{Comparing the Choices for the Test Norm and the Boundary Conditions}\label{sec:choice_of_norm}
\begin{figure}
  \centering
  \begin{subfigure}{6cm}
    \includegraphics[width = 6cm, height = 4.2cm]{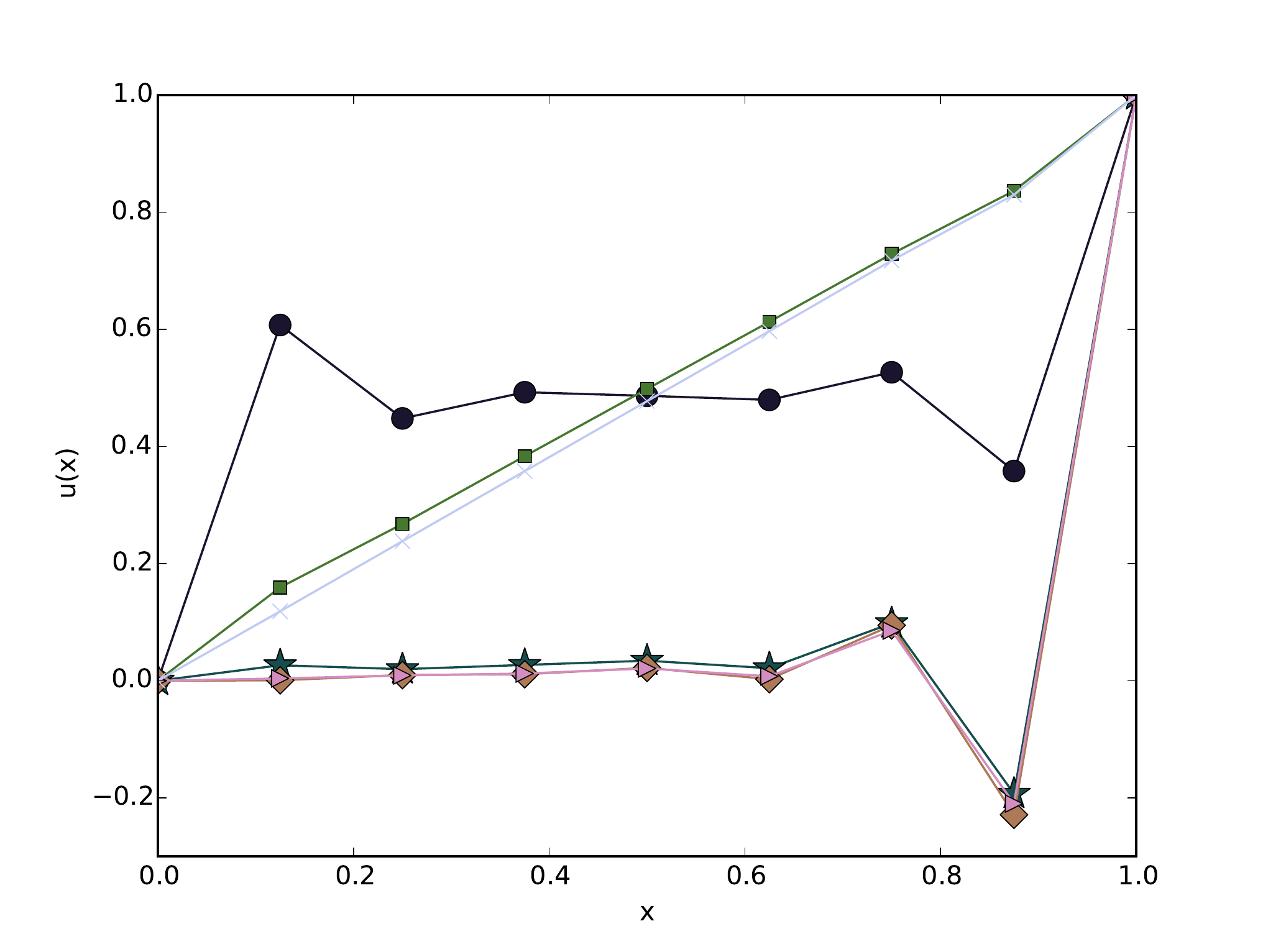}
    \caption{$\varepsilon = 10^{-3}$, $q=2$, $\alpha = 1$} \label{fig:norms_eps3_q2}
  \end{subfigure}
  \begin{subfigure}{6cm}
    \includegraphics[width = 6cm, height = 4.2cm]{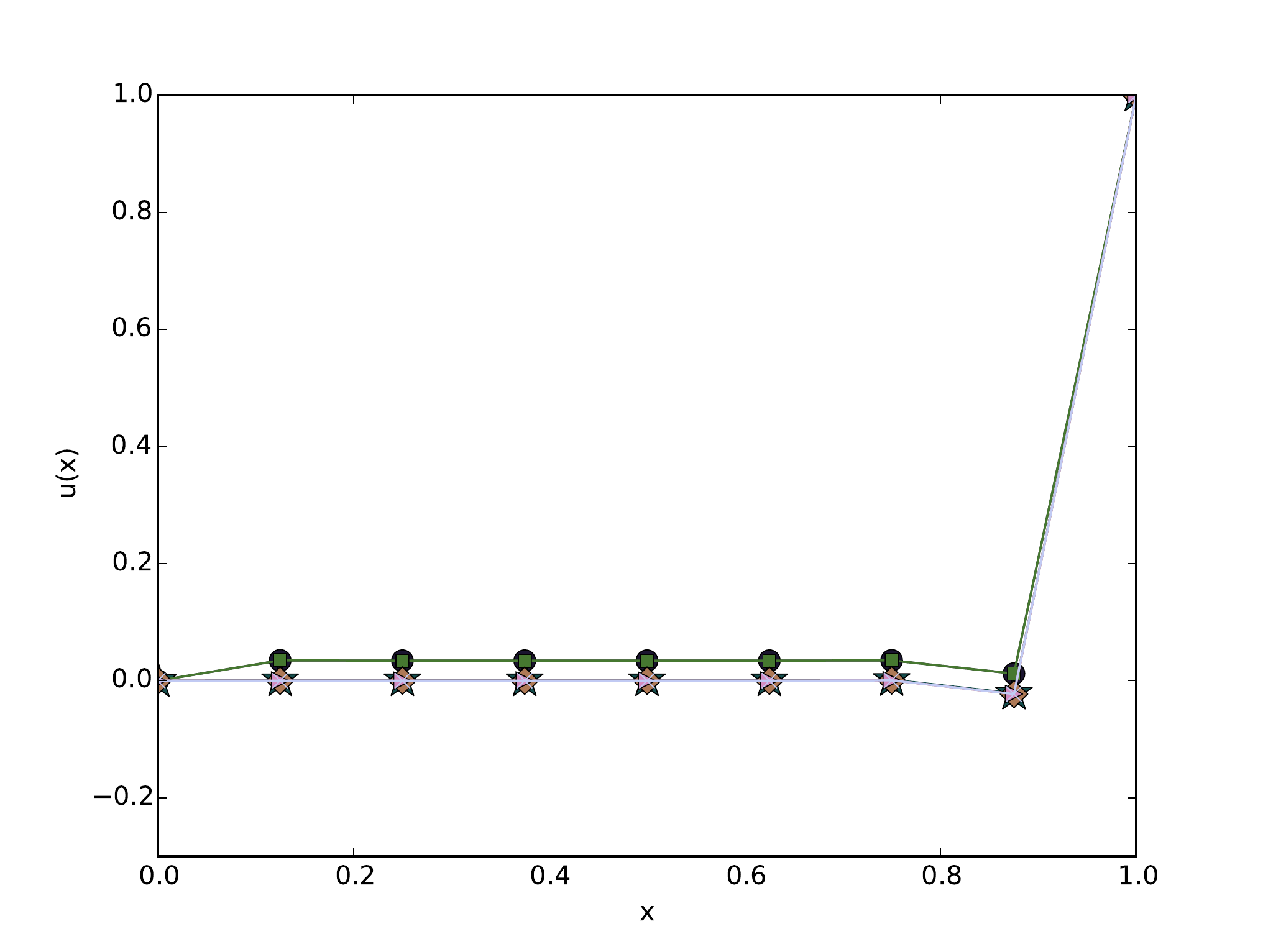}
    \caption{$\varepsilon = 10^{-3}$, $q=1.01$, $\alpha = 1$}\label{fig:norms_eps3_bc_q101}
  \end{subfigure}
  \begin{subfigure}{6cm}
    \includegraphics[width = 6cm, height = 4.2cm]{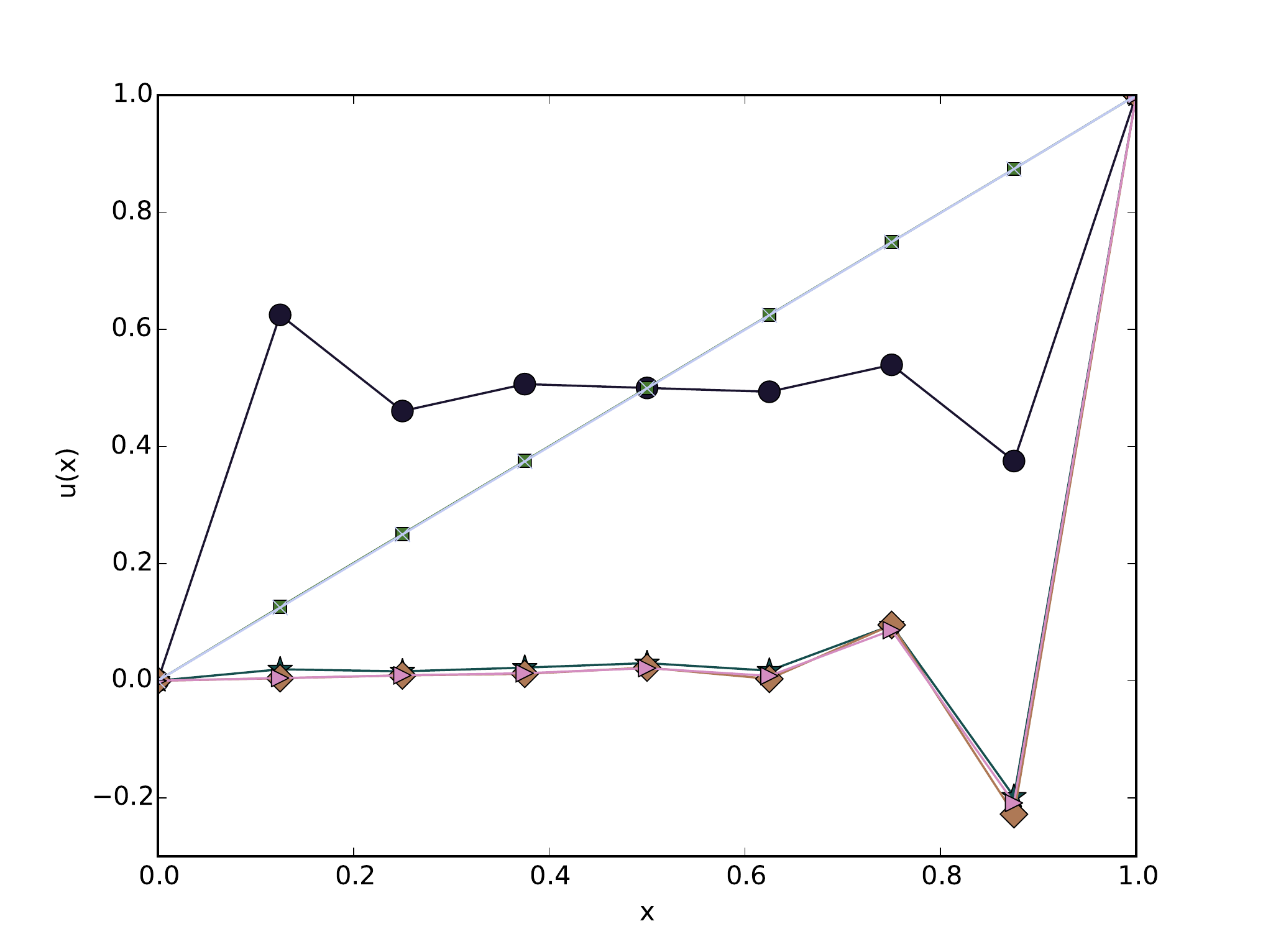}
    \caption{$\varepsilon = 10^{-6}$, $q=2$, $\alpha = 1$}\label{fig:norms_eps6_q2}
  \end{subfigure}
  \begin{subfigure}{6cm}
    \includegraphics[width = 6cm, height = 4.2cm]{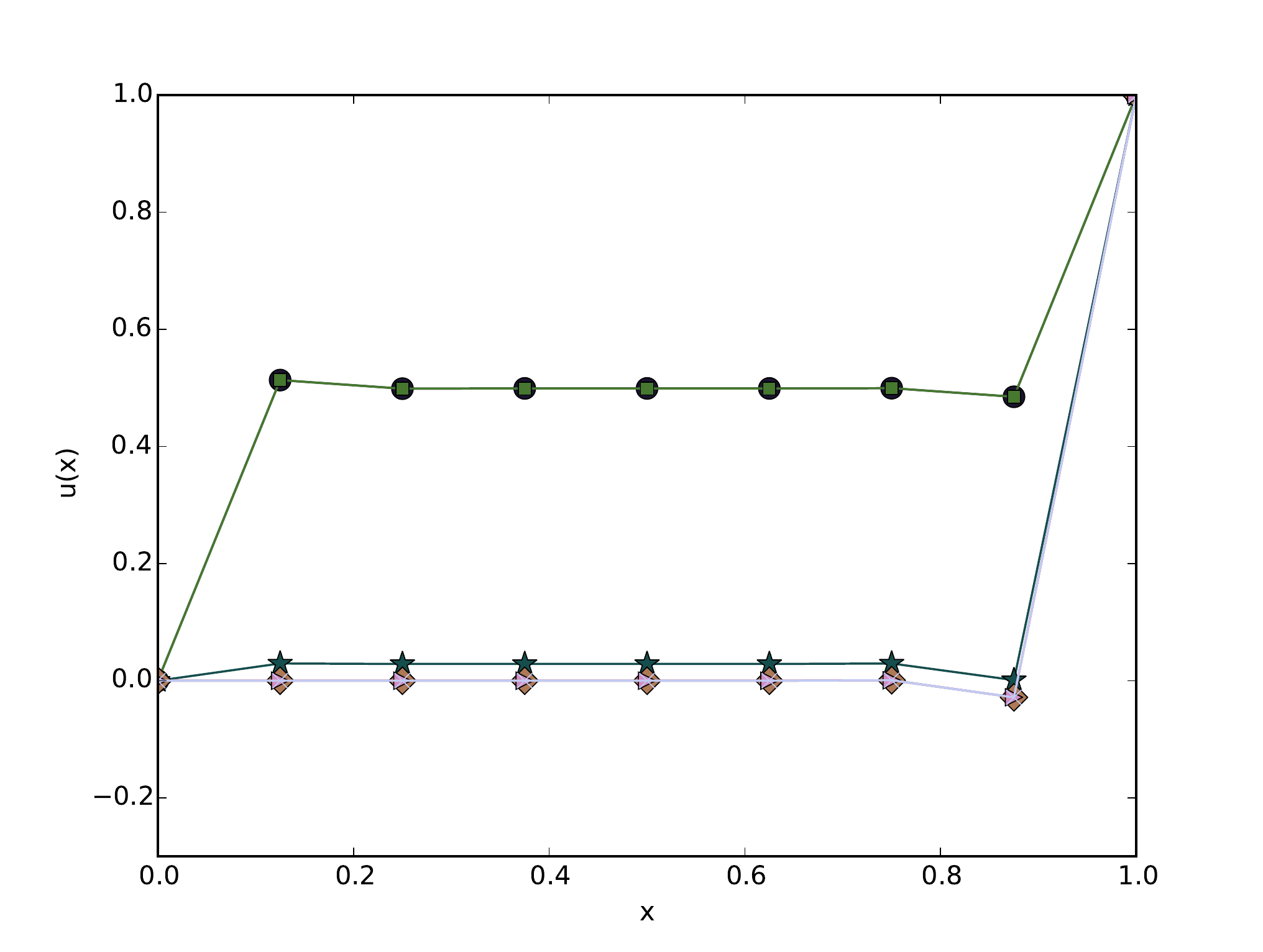}
    \caption{$\varepsilon = 10^{-6}$, $q=1.01$, $\alpha = 1$}\label{fig:norms_eps6_q101}
  \end{subfigure}
  \begin{subfigure}{6cm}
    \includegraphics[width = 6cm, height = 4.2cm]{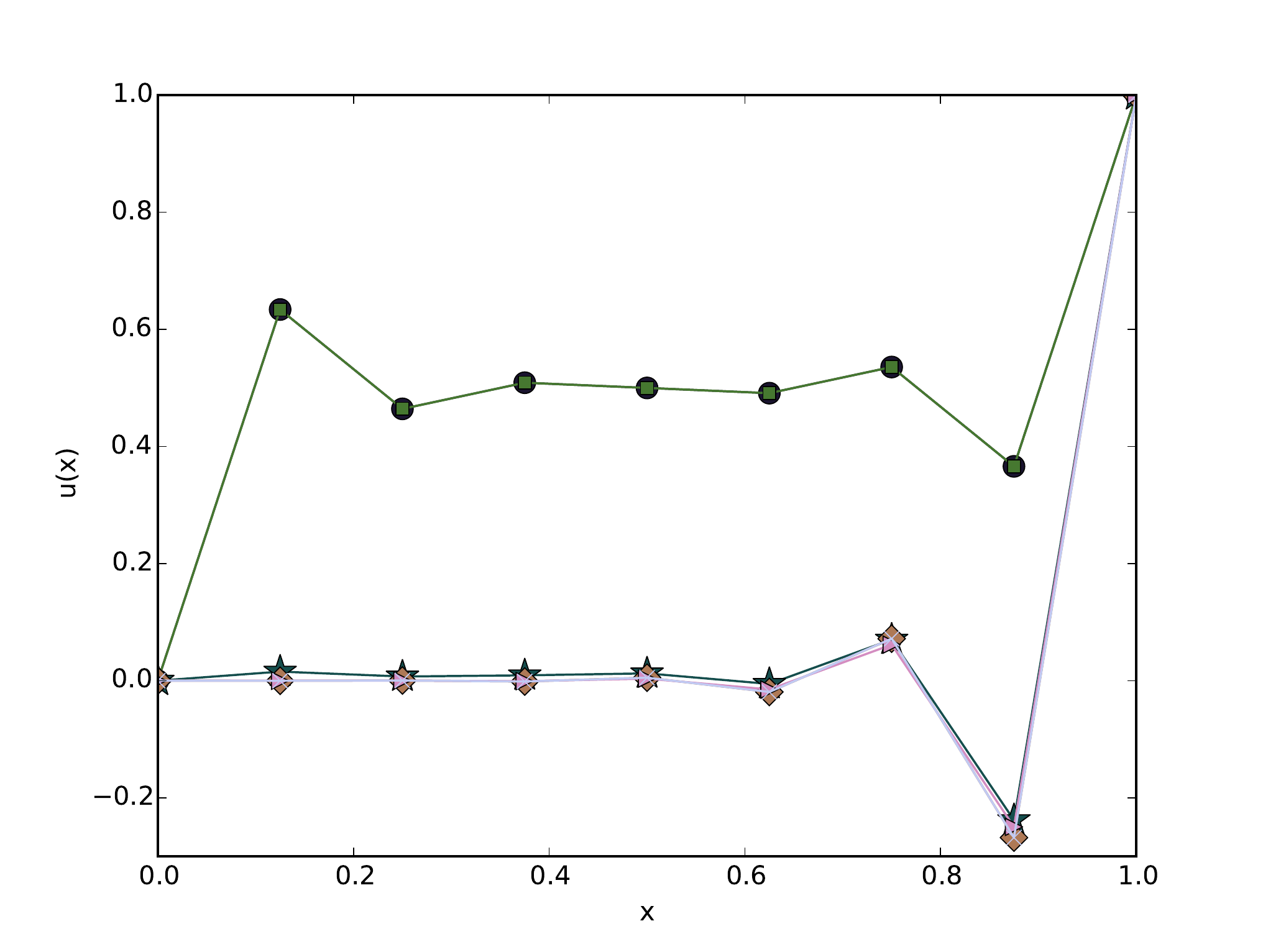}
    \caption{$\varepsilon = 10^{-6}$, $q=2$, $\alpha = 0$}\label{fig:norms_eps6_q2_alpha0}
  \end{subfigure}
  \begin{subfigure}{6cm}
    \includegraphics[width = 6cm, height = 4.2cm]{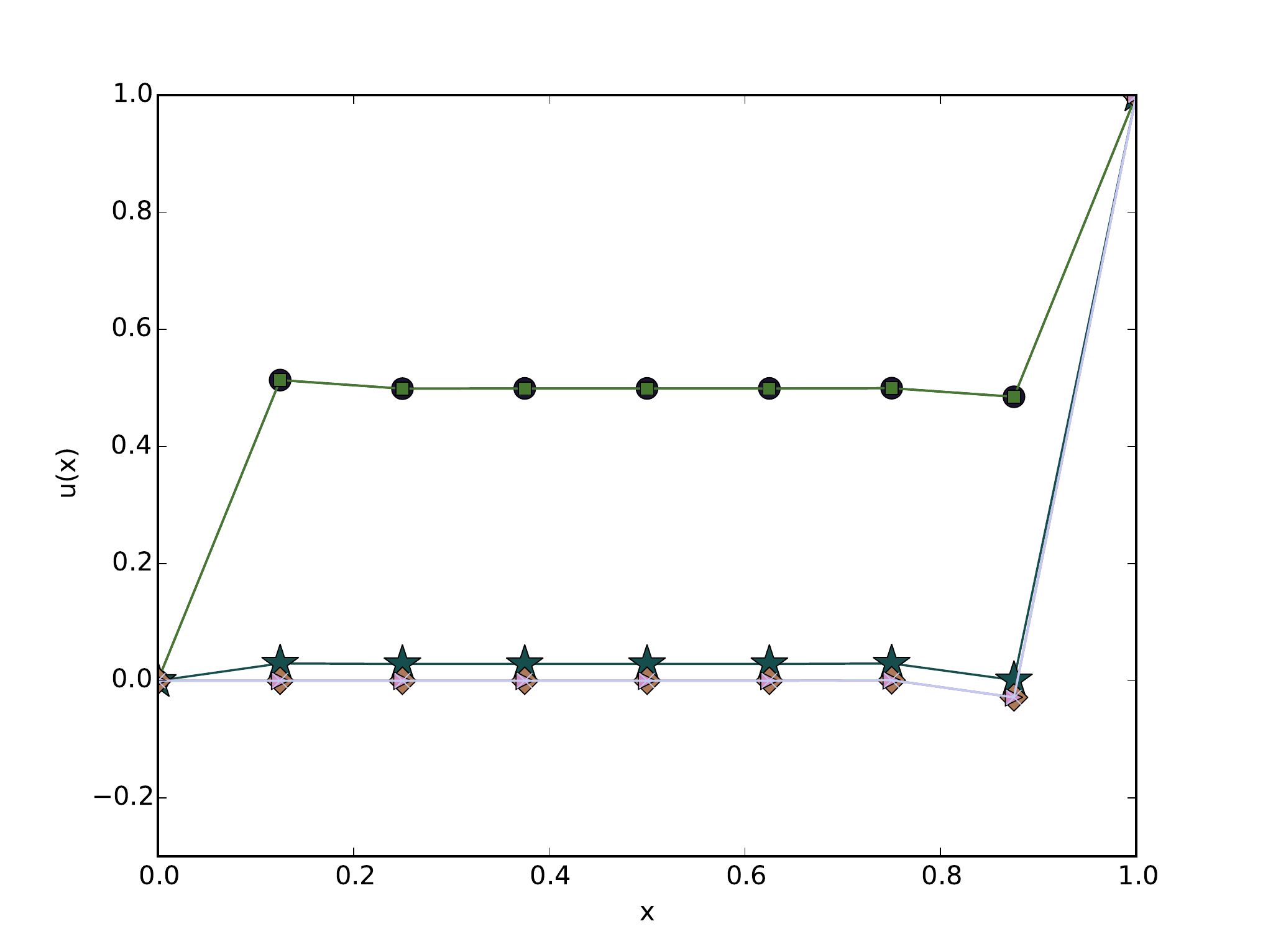}
    \caption{$\varepsilon = 10^{-6}$, $q=1.01$, $\alpha = 0$}\label{fig:norms_eps6_q101_alpha0}
  \end{subfigure}
  \begin{subfigure}{6cm}
    \centering
    \includegraphics[height = 1.8cm]{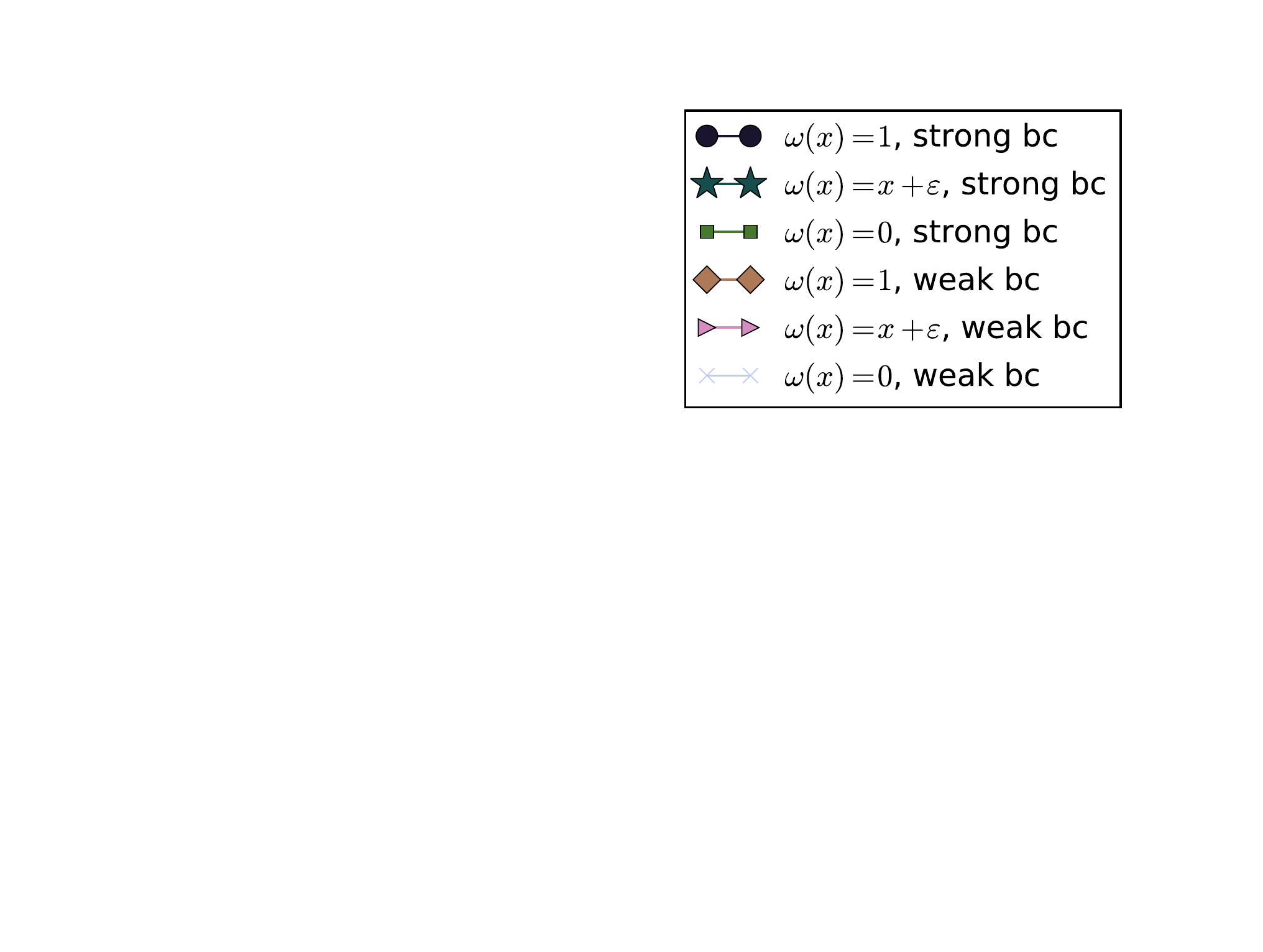}
  \end{subfigure}
  \caption{Solution $u_n$ for $q=2$ and $q=1.01$ using different norms, i.e., different weighting functions $\omega(x)$, and either Dirichlet boundary conditions for $r$ on the whole boundary (strong bc) or
  weak boundary conditions on $r$ on the inflow boundary and Dirichlet conditions on the outflow boundary (weak bc). }
  \label{fig:norms_1D}
\end{figure}

We start by comparing different versions of the non-linear Petrov-Galerkin method. In order to ensure that the space $V_m$ is sufficiently large we  use piecewise polynomials of degree $p_m = 10$ as the basis, while $U_n$ consists of piecewise linear polynomials. We split the interval $(0,1)$ uniformly into $8$ elements. We consider three different choices for the weighting function $\omega(x)$ as mentioned in Section \ref{sec:test_norms}, namely $\omega(x) = x+\varepsilon$, $\omega(x) \equiv 1$ and $\omega(x) \equiv 0$. We combine this with two different choices of boundary conditions on $r$, i.e., zero Dirichlet boundary conditions on the whole boundary and weak boundary conditions on the inflow boundary as described in Section
\ref{sec:weak_inflow_bc}. This creates six test cases; we first consider these six
test cases for $\varepsilon = 10^{-3}$ with $\alpha = 1$ and $q = 2, 1.01$. Secondly,
we consider all six test cases with  $\varepsilon = 10^{-6}$, $\alpha = 0,1$ and $q = 2, 1.01$
The solution $u_n$ for each of the cases is shown in Figure \ref{fig:norms_1D}.

We can see in Figures \ref{fig:norms_eps3_q2}, \ref{fig:norms_eps6_q2} and \ref{fig:norms_eps6_q2_alpha0} that for $q=2$ and Dirichlet boundary conditions on the whole boundary for $r$, the approximation $u_n$ only resembles the analytical solution if $\omega(x) = x+\varepsilon$.
Introducing weak boundary conditions on $r$ on the inflow boundary resolves this issue for $\omega(x) \equiv 1$ but the approximation with $\omega(x) \equiv 0$ only shows improvement if $\alpha = 0$.
These observations are the same for $\varepsilon = 10^{-3}$ and $\varepsilon = 10^{-6}$.
This indicates that we can circumvent constructing a non-constant function $\omega(x)$ -- which can be challenging for complicated geometries -- by introducing weak boundary conditions on $r$ on the inflow boundary.
These observations are consistent with the results for a very similar problem studied in  \cite{Chan2014}.
Furthermore, note that $\omega(x) =0$ and $\alpha = 1$ most closely resembles the method introduced in \cite{Guermond2004}.
In \cite{Guermond2004} it is demonstrated for a different (two-dimensional) example that the approximation for $q=2$ can be very inaccurate.

If we now consider $q=1.01$, cf., Figures \ref{fig:norms_eps3_bc_q101}, \ref{fig:norms_eps6_q101} and \ref{fig:norms_eps6_q101_alpha0}, we observe a different behaviour of the approximation. If $\varepsilon = 10^{-3}$, the approximations are all very similar and much closer to the analytical solution with close to no undershoots; we will investigate
the phenomenon of vanishing undershoots more closely both in one dimension and two dimensions later on. If
$\varepsilon = 10^{-6}$, the approximation with Dirichlet boundary conditions on the whole boundary and
$\omega(x) \equiv 0$ or $\omega(x) \equiv 1$ again leads to a very poor approximation to the analytical solution. This again points to robustness issues that for $q=1.01$ are only visible for much smaller $\varepsilon$ than for $q=2$. An interesting observation is that for $\alpha = 1$,
 the combination of weak inflow boundary conditions on $r$ and $\omega(x) \equiv 0$ is not distinguishable from other choices of $\omega(x)$ for $q = 1.01$ whereas, for $q = 2$ it is significantly different. The improvement of the approximation as $q \rightarrow 1$ for $\omega(x) \equiv 0$ and $\alpha = 1$ resembles the results presented in \cite{Guermond2004}.

 In conclusion, choosing the boundary conditions as described in Section \ref{sec:weak_inflow_bc} allows us to avoid constructing a non-constant weighting function $\omega(x)$. Close to $q=1$, the results are very similar for $\omega(x)\equiv 1$ and $\omega(x) \equiv 0$, whereas close to $q=2$, $\omega(x) \equiv 1$ clearly yields a significantly better approximation of the solution.

 The choice of $\alpha$ does not seem to have a big impact on the approximation except in the case $\omega(x)\equiv 0$ and $q = 2$, where choosing $\alpha = 0$ improves the approximation compared to $\alpha =1$.
 If $\omega(x) \equiv 1$, then both for $\alpha =1$ and $\alpha =0$, the dominant term in the norm on $V$ is essentially $\|\nabla u\|_{\Lp{q'}}$ since $\|\bm b\cdot \nabla v\|_{\Lp{q'}}$ and $\| \nabla v\|_{\Lp{q'}}$ are the same up to a constant in one dimension.
 If $\omega(x)\equiv 0$, then $\alpha = 0$ implies that $V$ is endowed with the $\sobzero{1}{q'}$-norm , whereas for $\alpha = 1$ the $\Lp{q'}$-term is the dominant term in the norm on $V$ due to the $\varepsilon$-weighting of the gradient. Since
 $\alpha = 0$ yields a better approximation in this case, this suggests that for $q=2$ it is favourable to choose a stronger norm on $V$ that is either dominated by the norm of the gradient or the streamline term $\|\bm b\cdot \nabla v\|_{\Lp{q'}}$.

\subsubsection{Vanishing oscillations as $q \rightarrow 1$}\label{sec:vanishing_osc}
\begin{figure}[!t]
\centering
\includegraphics[width = 12cm]{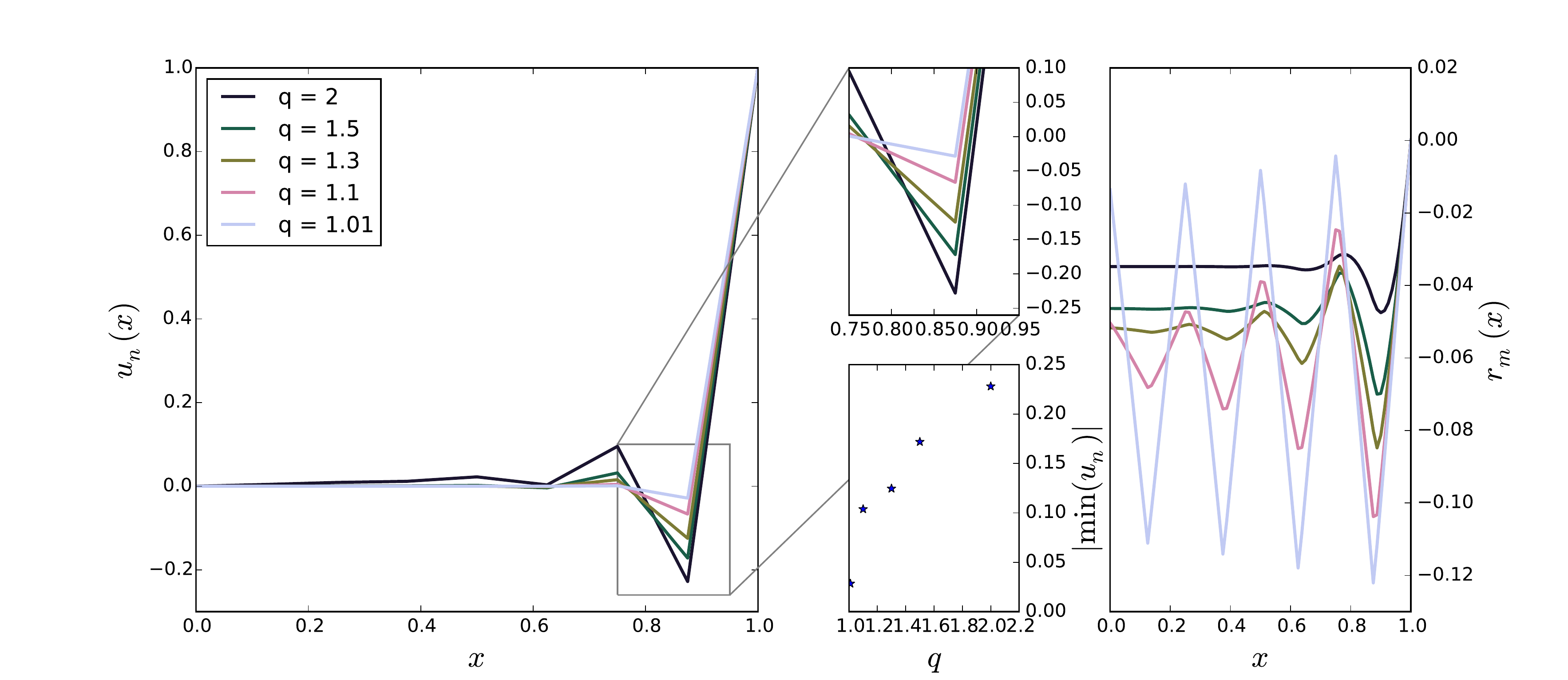}
\caption{Example \ref{ex:1dconfusion} with $\varepsilon = 10^{-5}$:
Numerical approximations $u_n$ (left) and $r_m$ (right) with a uniform mesh consisting of 8 intervals
using piecewise linear polynomials for $u_n$ and polynomials of degree $p_m = 10$ for $r_m$
for varying values of $q$. Top center: Zoom-in showing over- and undershoots.
 Bottom center: $|\min (u_n)|$ vs $q$.}
\label{fig:osc_u}
\end{figure}

The examples in Figure \ref{fig:norms_1D} already illustrate that the
undershoot in the approximate solution nearly vanishes for $q = 1.01$.
We will now investigate in more detail how the undershoot depends on $q$
if we choose $\varepsilon = 10^{-5}$, $\omega(x) \equiv 1$, $\alpha = 1$ and impose weak inflow boundary conditions on $r_m$.
To this end, we choose a large space $V_m$ with polynomial degree $p_m = 10$
and piecewise linear polynomials for $U_n$. We split the interval into
$8$ elements and compute the approximate solutions $(r_m,u_n)$ for
several choices of $q$. Figure \ref{fig:osc_u} shows $u_n$ on the left and $r_m$ on the right. A close-up
of the undershoot and a plot $\min(u_n)$ vs. $q$ is shown in the centre of the figure.
We can see that the undershoot decreases monotonically as $q$ approaches $1$.
In contrast to this, we can see that the oscillations in the residual $r_m$
increase as $q' \rightarrow \infty$ which suggests that a large space $V_m$ may
be necessary. This is the next aspect of the method we will investigate.

\subsubsection{The Choice of the Space $V_m$} \label{sec:choice_Vm}
In \cite{Muga2017}, it has been shown that the inexact method is well-posed if a
Fortin projector, cf., \eqref{eq:Fortin_operator} exists. It is easy to see
that $\dim (V_m) \geq \dim(U_n)$ is a necessary condition. Finding a sufficient condition or
in other words finding a compatible pair $(U_n, V_m)$ is highly non-trivial. In \cite{MugTylZee2018},
certain special cases for the convection-reaction equation are considered. The observations in Section \ref{sec:convergence} suggest that for $\varDelta p \geq 2$, the Fortin condition is typically satisfied.
We now investigate how the choice of $V_m$ affects the undershoot for
$q$ close to $1$. To this end, we consider two different strategies of enlarging
$V_m$: global $p$-enrichment and uniform $h$-refinement. We again consider $\varepsilon = 10^{-5}$, $\omega(x) \equiv 1$, $\alpha = 1$ and impose weak inflow boundary conditions on $r_m$.

\begin{figure}[!t]
\centering
\includegraphics[width = 12cm]{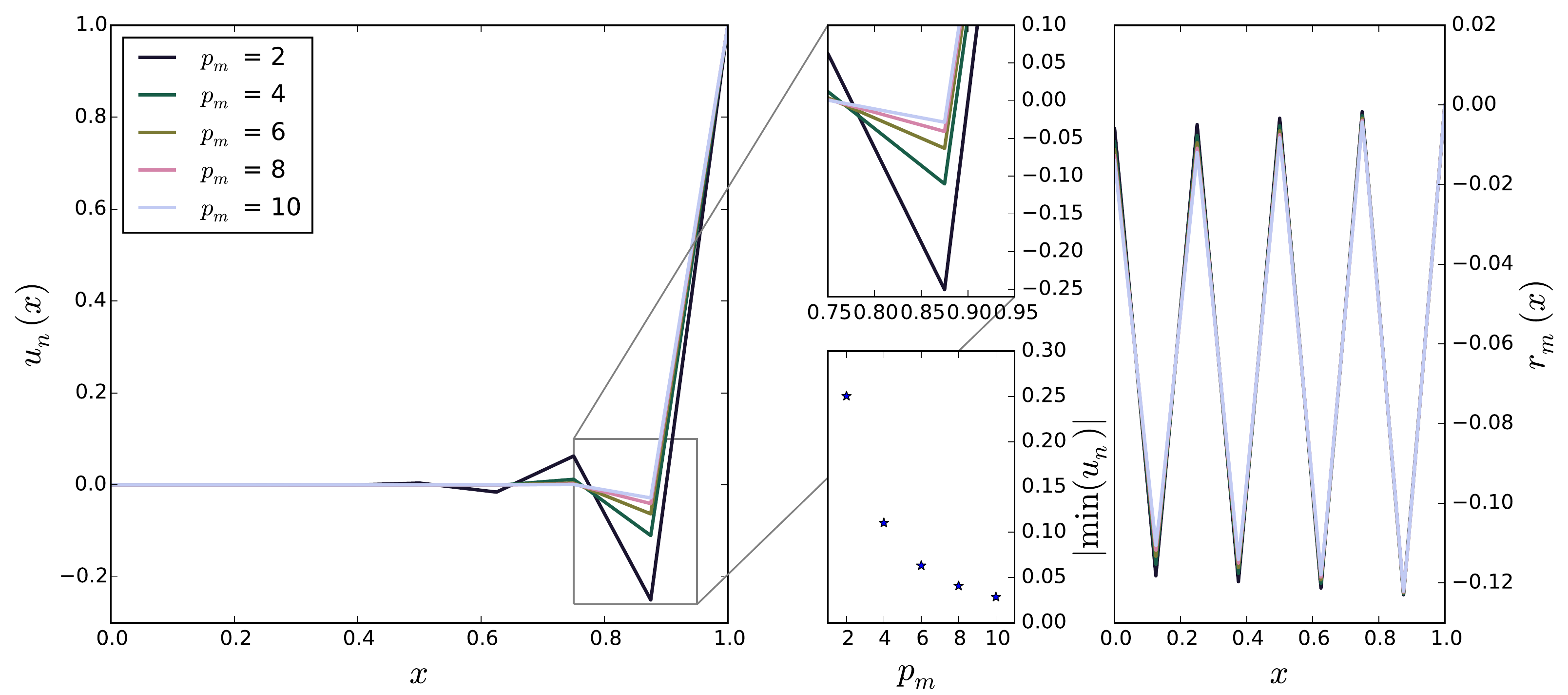}
\caption{Example \ref{ex:1dconfusion} with $\varepsilon = 10^{-6}$ and $q = 1.01$, $q' = 101$: Numerical approximations $u_n$ (left) and $r_m$ (right) with a uniform mesh consisting of 8 intervals using piecewise linear polynomials for $u_n$ and varying polynomial degrees $p_m$ for $r_m$. Top center: Zoom-in showing over- and undershoots. Bottom center: $|\min (u_n)|$ vs $p_m$. Right: Projection of $r_m$ onto the space of piecewise linear functions with a uniform mesh consisting of 16 intervals using piecewise linear polynomials for $u_n$ and varying polynomial degrees $p_m$ for $r_m$}
\label{fig:p_u}
\end{figure}

\textbf{p-enrichment:}
In terms of implementation, the simplest choice for the space $V_m$
such that $\mathrm{dim}(V_m) > \mathrm{dim}(U_n)$ is a finite element
space over the same grid with the polynomial degree increased globally
by some integer $\varDelta p$. Figure \ref{fig:p_u} shows how the overshoot reduces as we increase $\varDelta p$.
The difference in $r_m$ for varying $p_m$ is barely visible even though we have used a finer mesh for plotting in order to capture the
behaviour of the higher modes and yet it has a significant effect on the approximation $u_n$.

\begin{figure}[!t]
\centering
\includegraphics[width = 12cm]{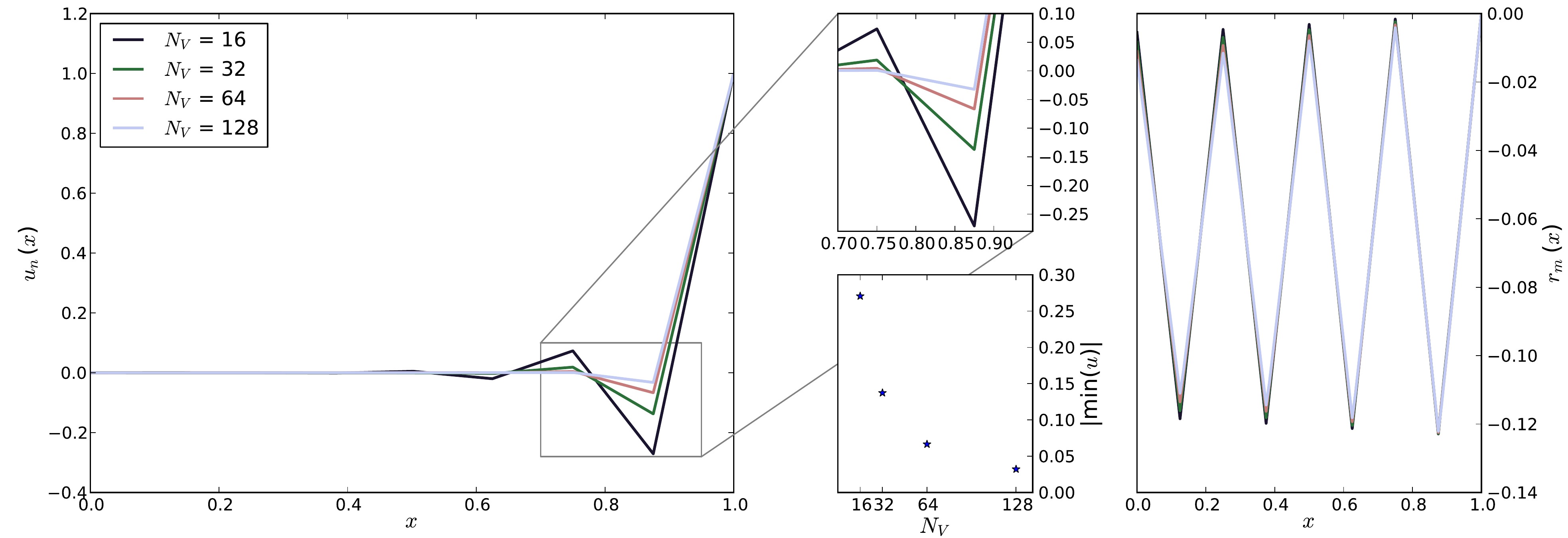}
\caption{Example \ref{ex:1dconfusion} with $q = 1.01$, $q' = 101$: Numerical
approximations $u_n$ (left) and $r_m$ (right) with a uniform mesh consisting of 16 intervals using piecewise linear polynomials for $u_n$ and polynomials of degree $p_m = 10$ for $r_m$
for varying values of $\varepsilon$.
Top center: Zoom-in showing over- and undershoots.
Bottom center: $|\min (u_n)|$ vs $\varepsilon$.}
\label{fig:N_u}
\end{figure}

\textbf{h-refinement:}
This time both $U_n$ and $V_m$ consist of piecewise linear polynomials.
To ensure that $\mathrm{dim}(V_m) > \mathrm{dim}(U_n)$, we choose a refinement of
the underlying mesh of the space $U_n$ to construct the space $V_m$. This is the only numerical
experiment in this article that cannot be implemented in FEniCS; instead, e.g., the C++ library Hermes2D \cite{Solin2014} can be used.
 Figure \ref{fig:N_u} shows how the overshoot reduces as we refine the mesh for $V_m$;
the difference in $r_m$ is again barely visible despite the significant effect on the approximation $u_n$.

\subsubsection{Robustness in $\varepsilon$}\label{sec:eps_robust}
\begin{figure}[!t]
\centering
\includegraphics[width = 12cm]{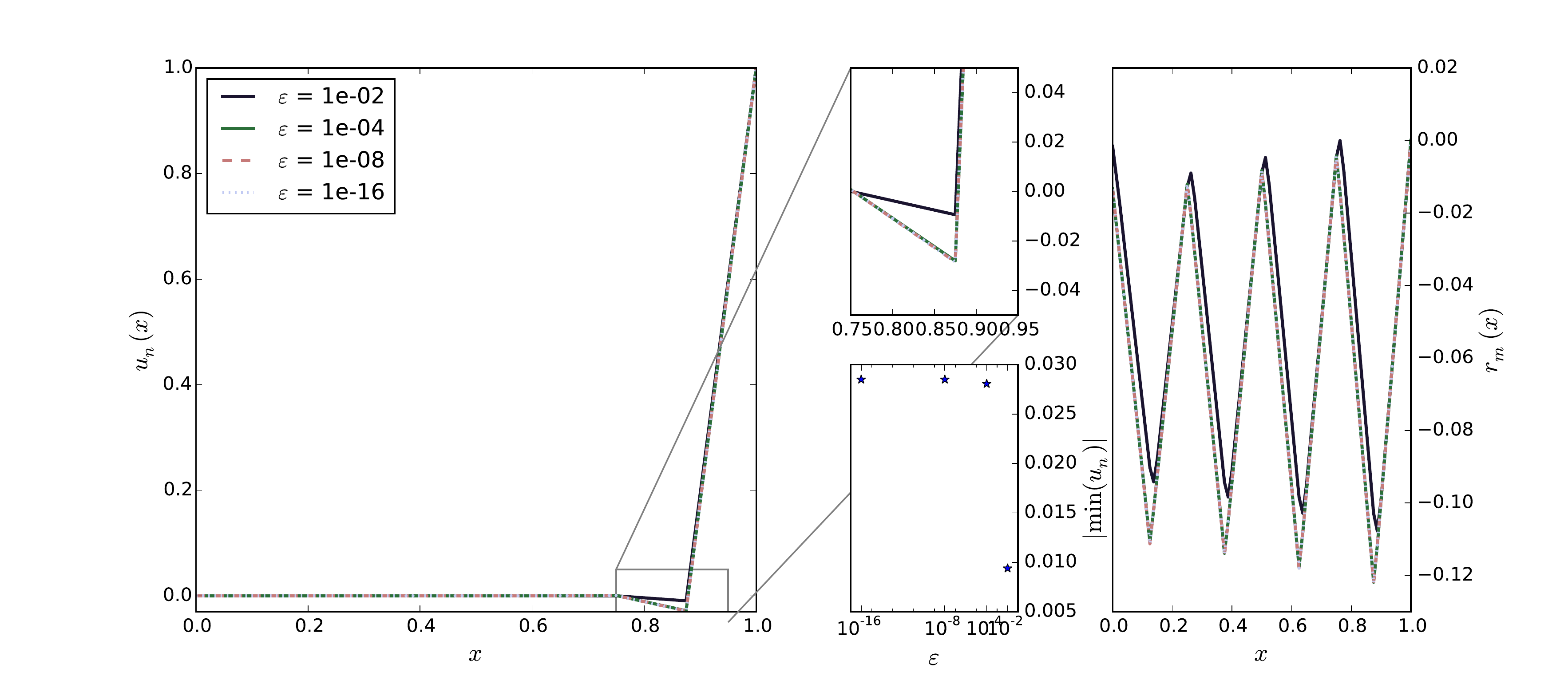}
\caption{Example \ref{ex:1dconfusion} with $q = 1.01$, $q' = 101$: Numerical
approximations $u_n$ (left) and $r_m$ (right) with a uniform mesh consisting of 16 intervals using piecewise linear polynomials for $u_n$ and polynomials of degree $p_m = 10$ for $r_m$
for varying values of $\varepsilon$.
Top center: Zoom-in showing over- and undershoots.
Bottom center: $|\min (u_n)|$ vs $\varepsilon$.}
\label{fig:eps_u}
\end{figure}

As a final experiment in one dimension, we study how $\varepsilon$ affects the
undershoot in the approximation. We continue with the same setting as in the previous section, but keep $V_m$ fixed as the space of piecewise polynomials of degree $p_m = 10$ on the same mesh as used for $U_n$ and vary $\varepsilon$ instead.
Figure \ref{fig:eps_u} shows that, although the undershoot is not the same for
all $\varepsilon$ leading to an under resolved layer, the undershoot seems to
be approximately the same for all $\varepsilon \leq \varepsilon_0$ for some
$\varepsilon_0$. This can be traced back to the inflow boundary conditions
on $r_m$ ---  which are scaled with $\varepsilon$ --- if we look at $r_m$ on the right in
Figure \ref{fig:eps_u}.  The method therefore seems to be robust in
$\varepsilon$.

\subsection{Vanishing Oscillations in Two Dimensions}\label{sec:numerics_2D}

In this section we explore whether the under- and overshoots in the approximation still
vanish as $q \rightarrow 1$ if we apply the method to two-dimensional problems. To this end, we consider the Examples \ref{ex:eriksson-johnson}, \ref{ex:skew} and \ref{ex:interior_layer} on different meshes for $\varepsilon = 10^{-6}$. We consider $\omega(x) \equiv 1$, $\alpha = 1$ and impose weak boundary conditions on $r_m$ as described in Section \ref{sec:weak_inflow_bc}.
In one case we will also consider $\alpha = 0$ to compare the approximations for both choices of $\alpha$.
We will see that the overshoots disappear on certain meshes, but remain present on others. We will furthermore demonstrate that the approximation qualitatively behaves like the $L^q(\Omega)$-best approximation of the analytical solution and thus we can apply the observations in \cite{Houston2019a} to (a) predict on which meshes the overshoot will disappear as
$q \rightarrow 1$ and (b) design meshes that have this property in specific situations.

\subsubsection{The Eriksson-Johnson Model Problem}
\begin{figure}[t!]
  \centering
  \begin{subfigure}{2.9cm}
    \includegraphics[width = 2.9cm]{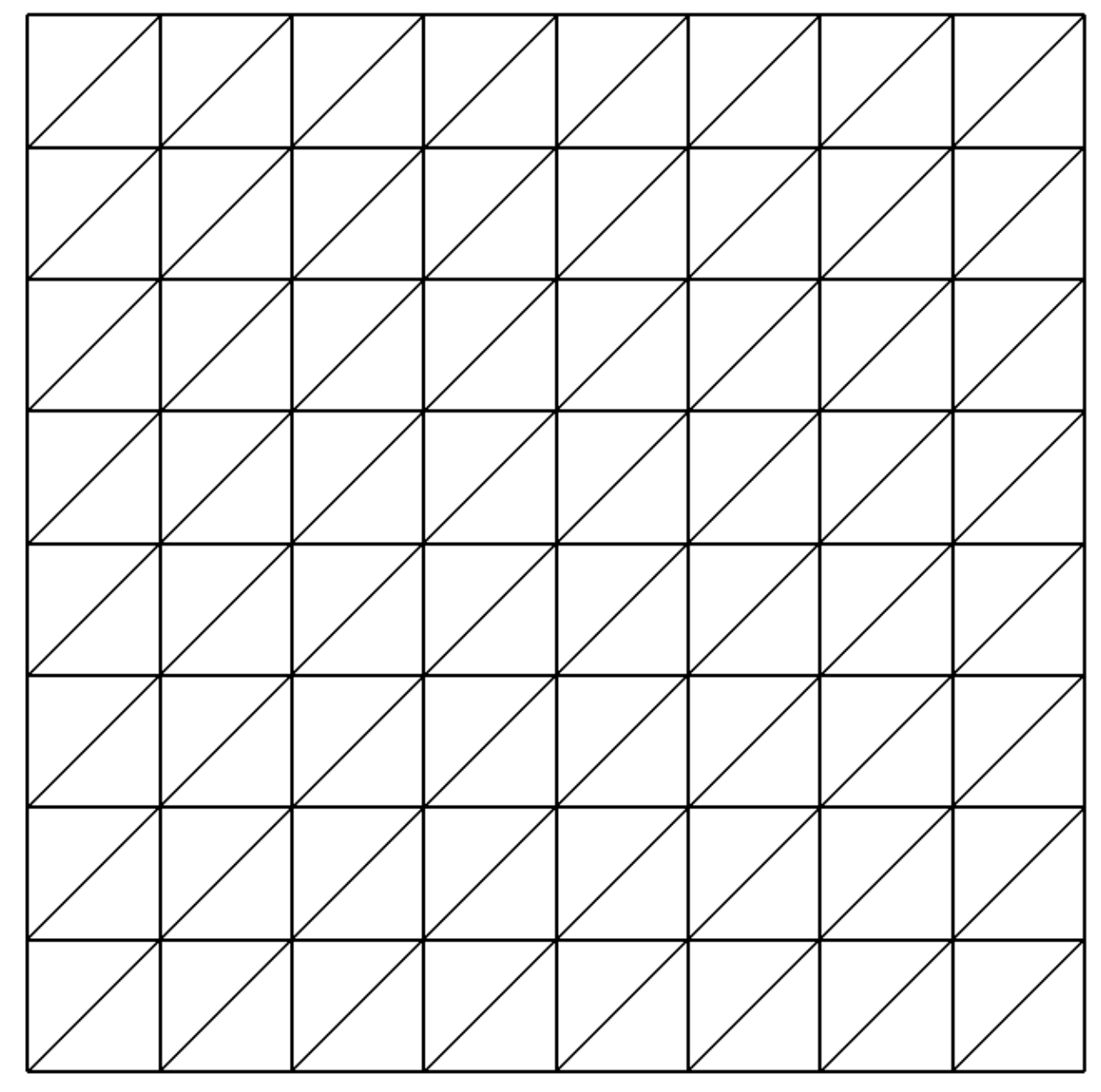}
    \caption{Mesh 1}
  \end{subfigure}
  \begin{subfigure}{2.9cm}
    \includegraphics[width = 2.9cm]{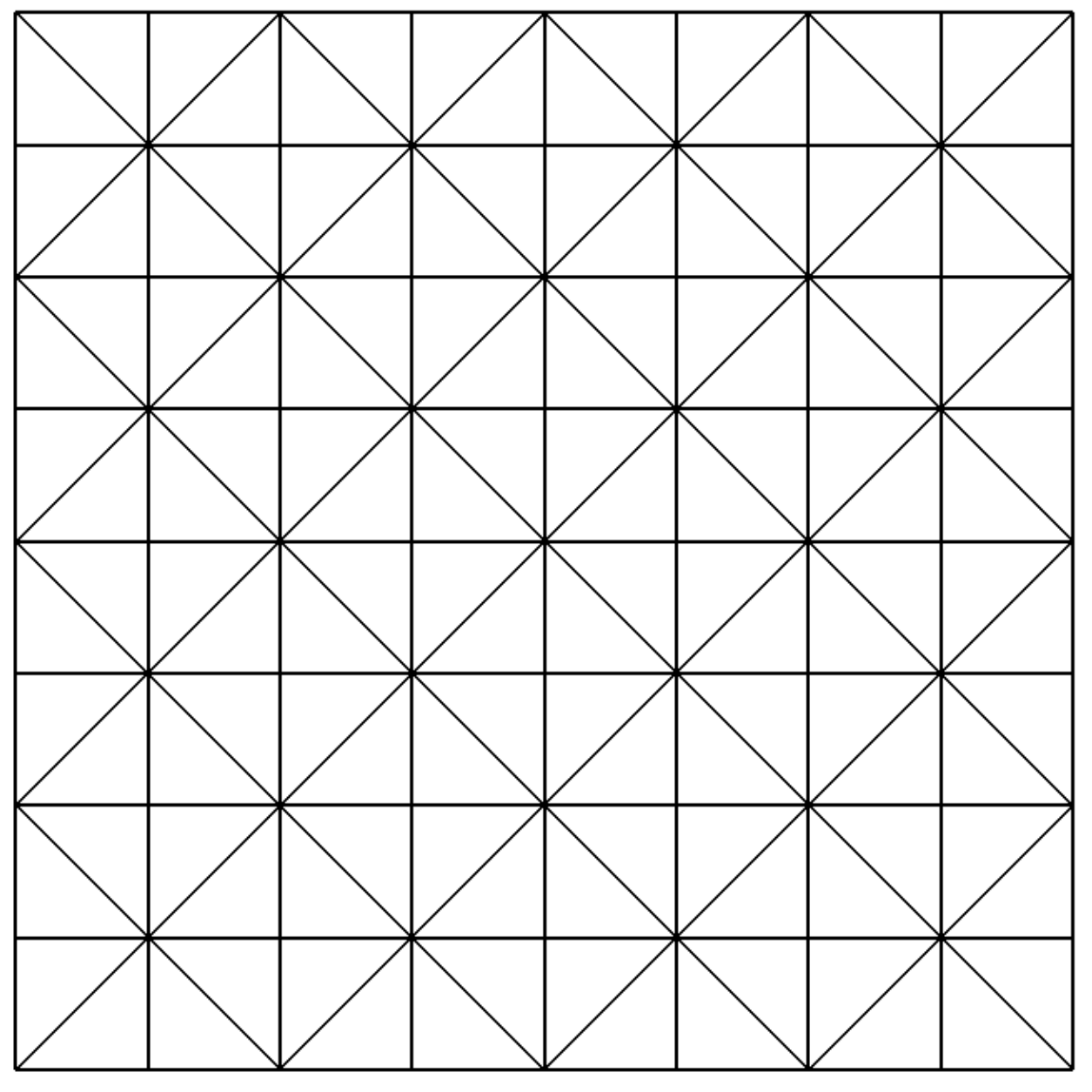}
    \caption{Mesh 2}
  \end{subfigure}
  \begin{subfigure}{2.9cm}
    \includegraphics[width = 2.9cm]{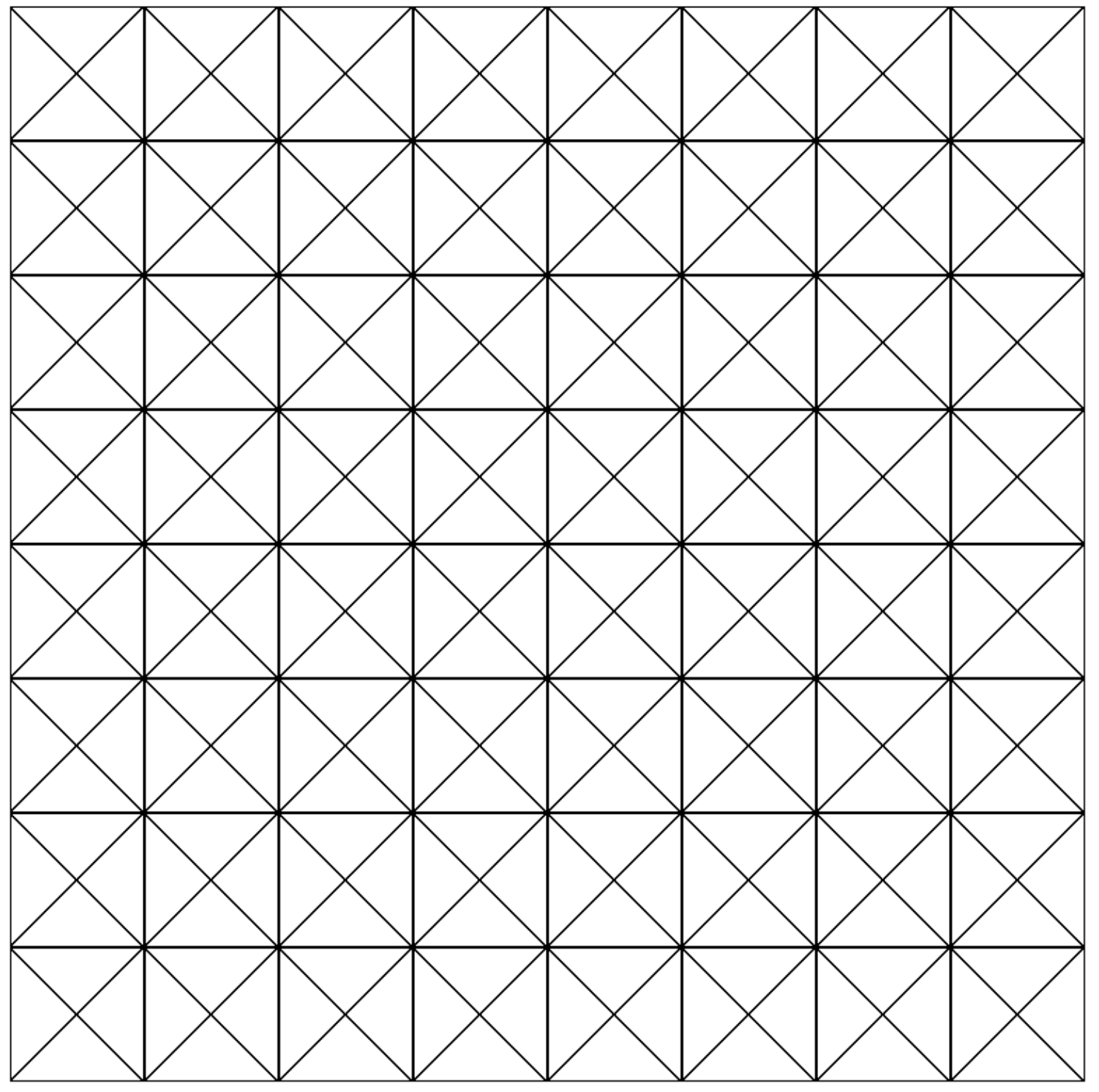}
    \caption{Mesh 3}
  \end{subfigure}
  \begin{subfigure}{2.9cm}
    \includegraphics[width = 2.9cm]{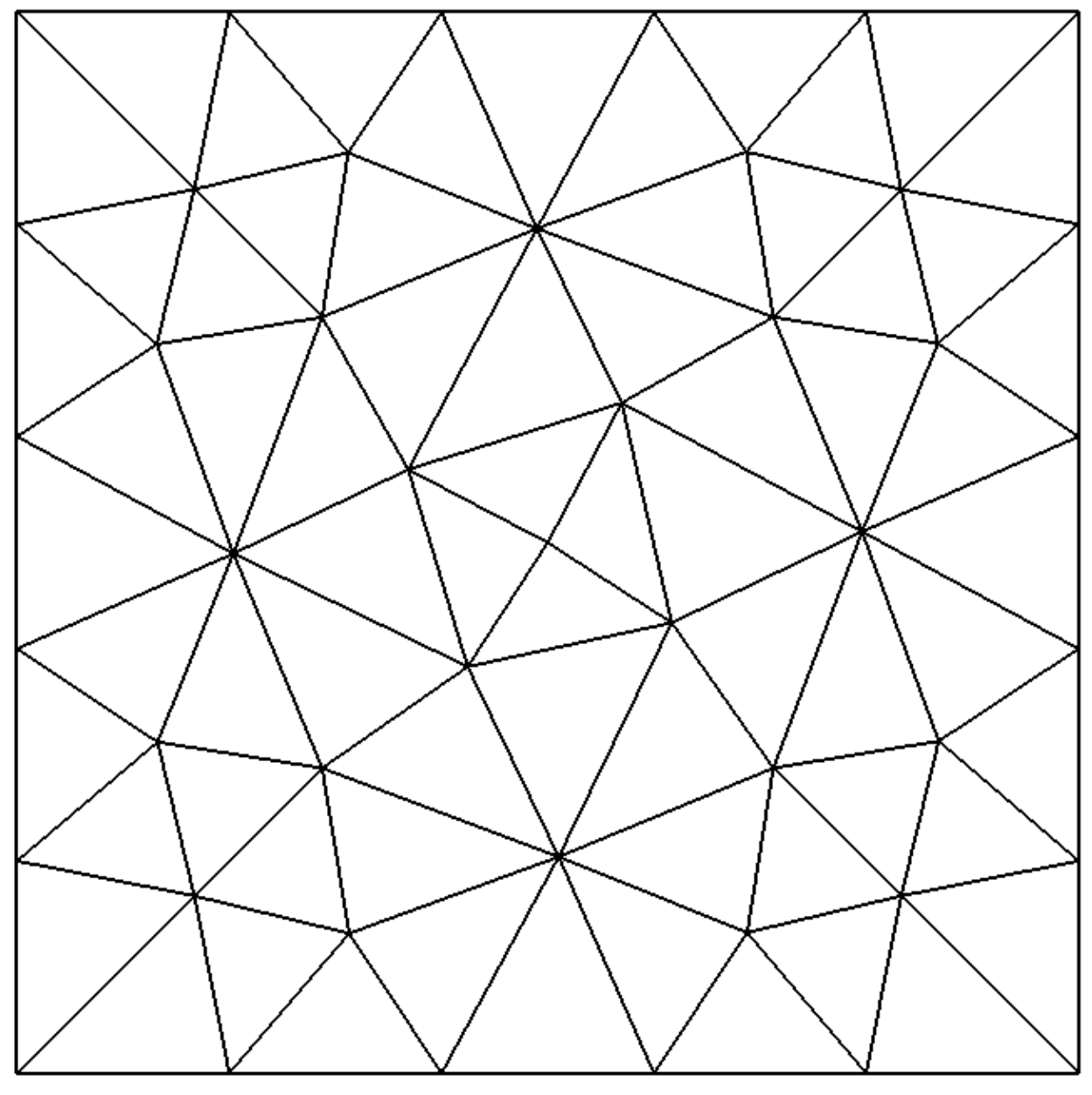}
    \caption{Mesh 4}
  \end{subfigure}
  \caption{Four Different Meshes on the Unit Square}
  \label{fig:meshes}
\end{figure}

We consider the Eriksson-Johnson Model Problem (Example \ref{ex:eriksson-johnson}) on the four different meshes depicted in Figure \ref{fig:meshes}. Note that these are the same meshes that are investigated in \cite{Houston2019a}.

\begin{figure}[!t]
  \centering
  \begin{subfigure}[t]{3.9cm}
    \centering
    \includegraphics[width = 3.9cm]{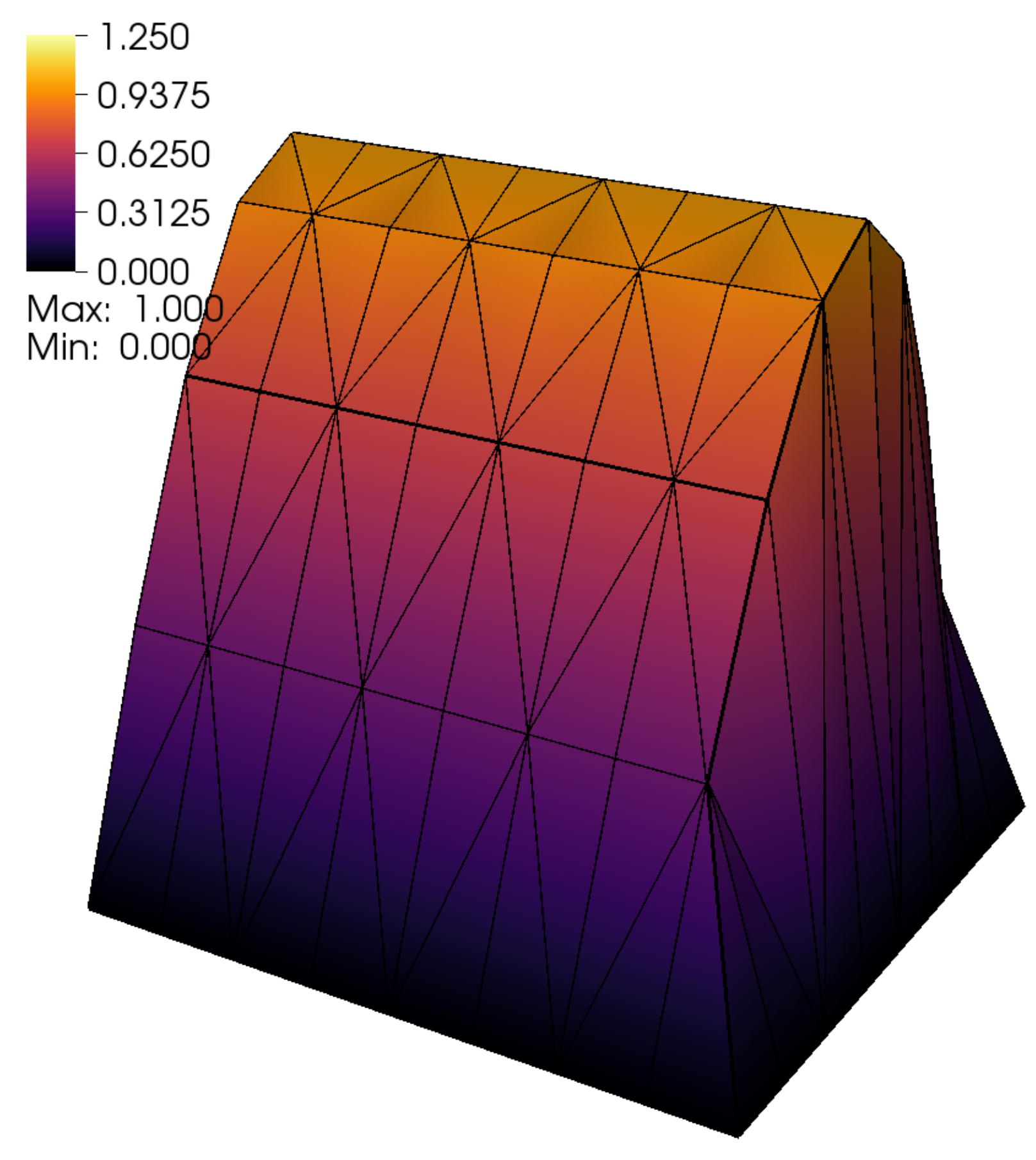}
    \caption{Interpolant of the analytical solution}
  \end{subfigure}
  \begin{subfigure}[t]{3.9cm}
    \centering
    \includegraphics[width = 3.9cm]{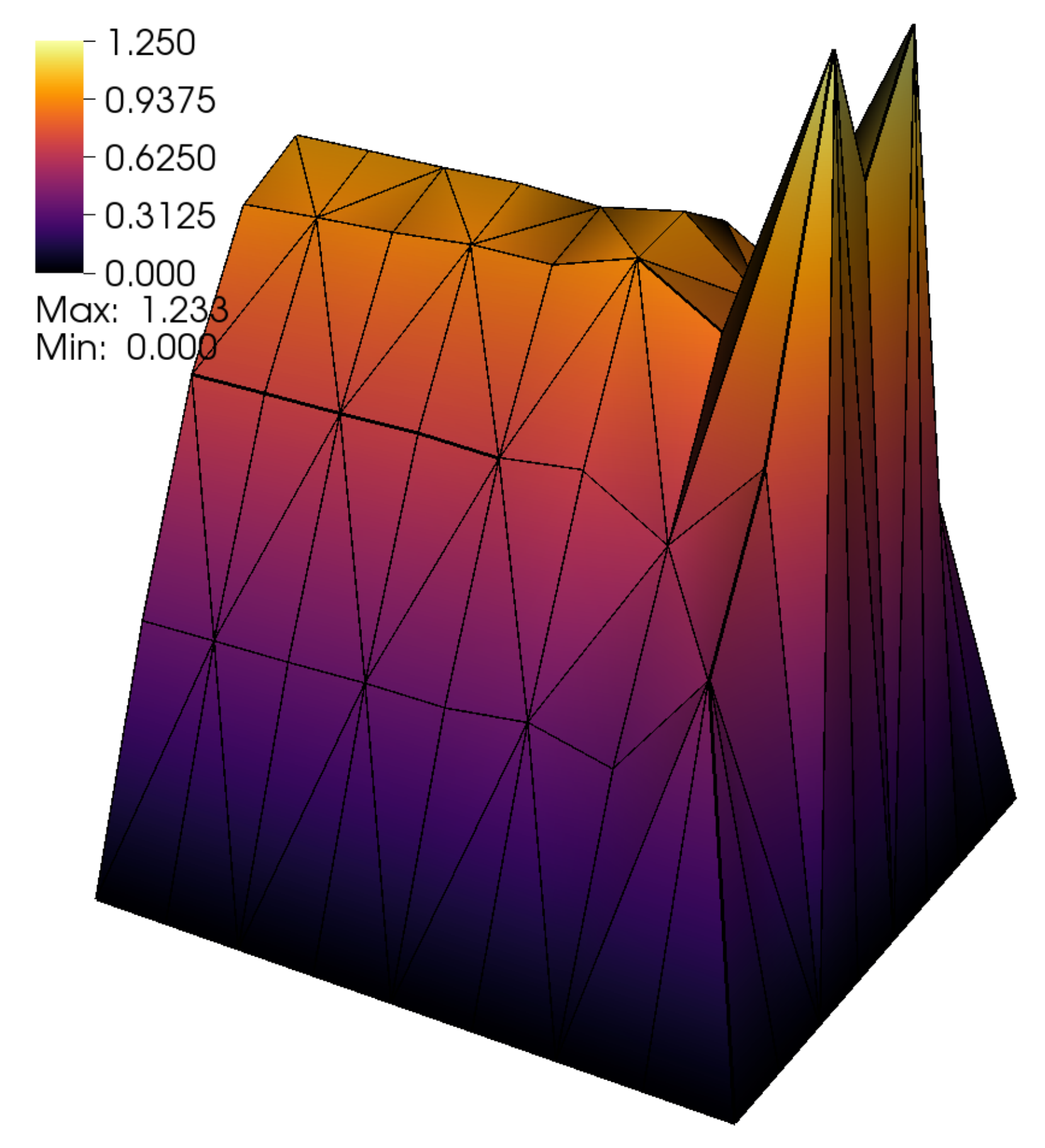}
    \caption{Approximation with $q = 2$}
  \end{subfigure}
  \begin{subfigure}[t]{3.9cm}
    \centering
    \includegraphics[width = 3.9cm]{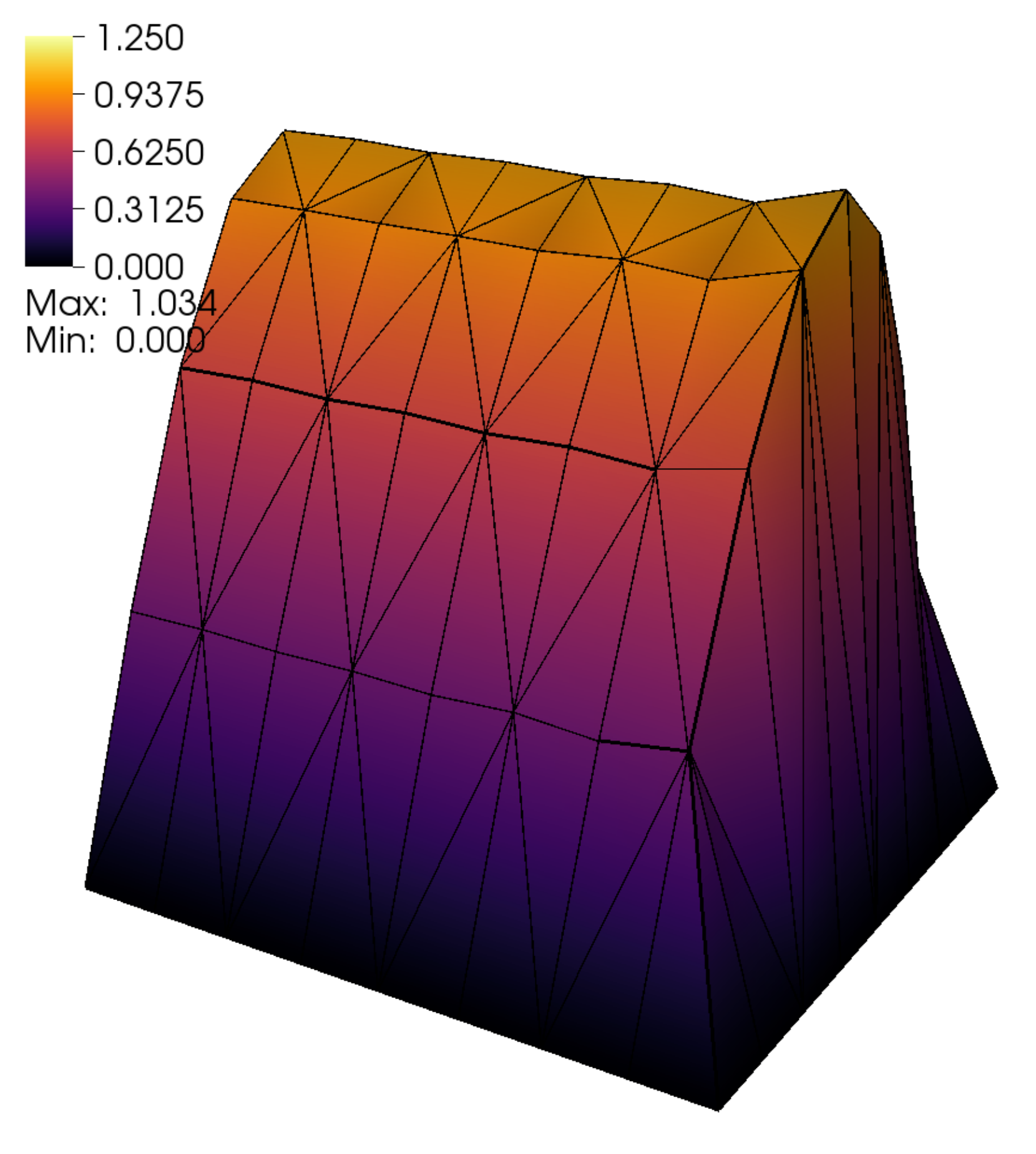}
    \caption{Approximation with $q = 1.01$}
  \end{subfigure}
  \caption{Approximation of the analytical solution of Example \ref{ex:eriksson-johnson} (Eriksson-Johnson model problem) on Mesh 2 for $\varepsilon = 10^{-6}$.}\label{fig:ej_left_right}
\end{figure}

\begin{figure}[!t]
  \centering
  \begin{subfigure}[t]{3.9cm}
    \centering
    \includegraphics[width = 3.9cm]{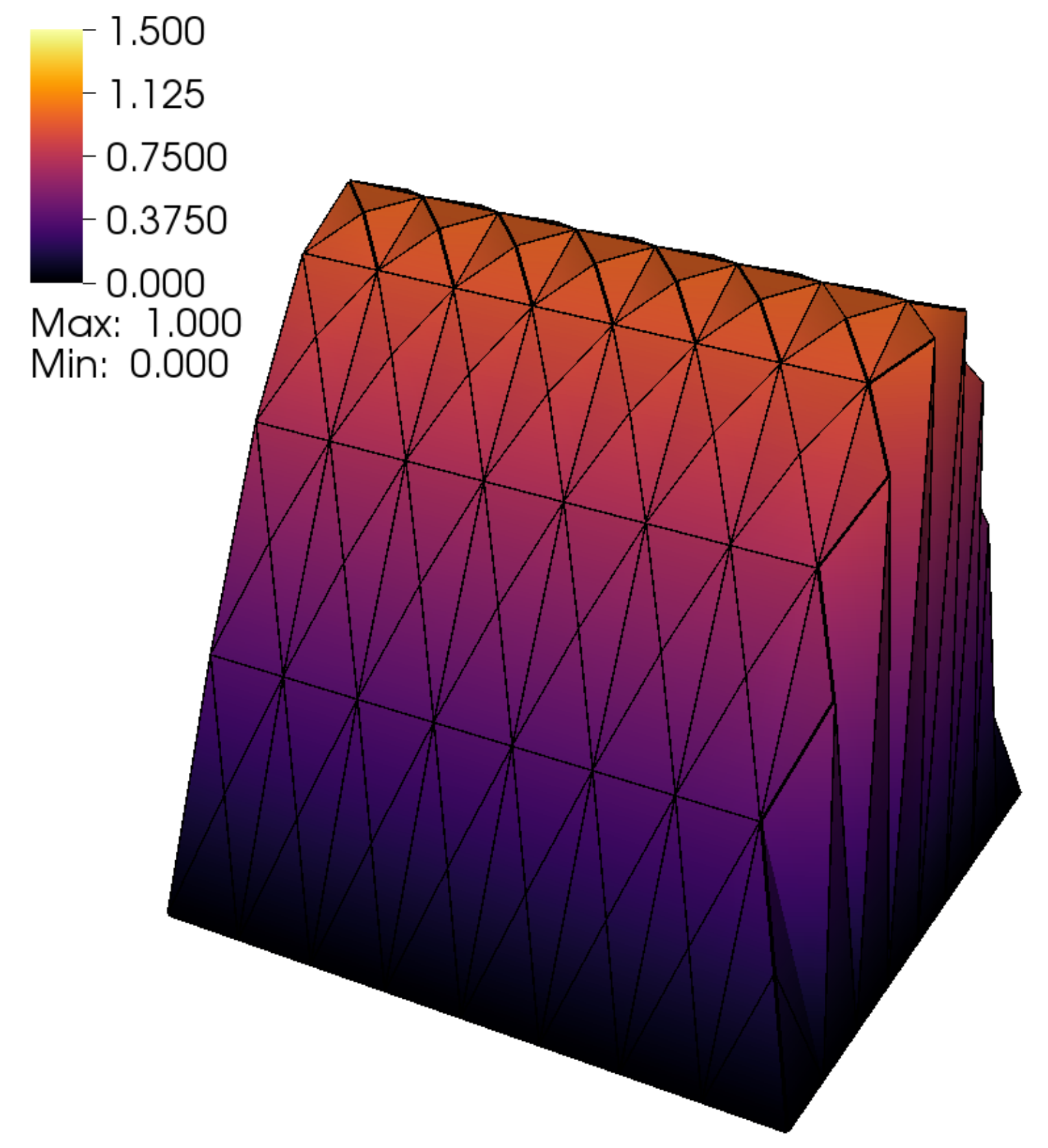}
    \caption{Interpolant of the analytical solution}
  \end{subfigure}
  \begin{subfigure}[t]{3.9cm}
    \centering
    \includegraphics[width = 3.9cm]{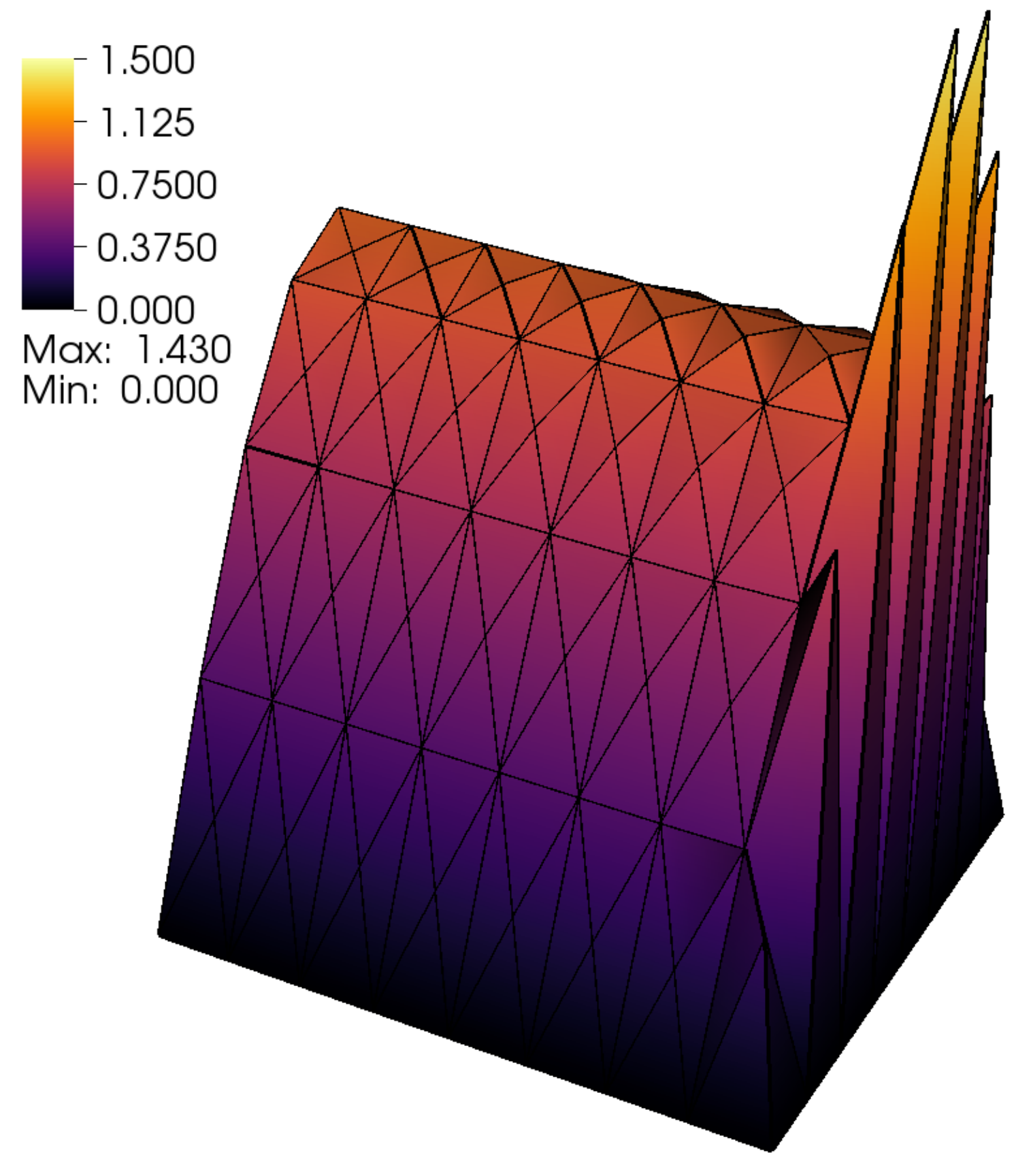}
    \caption{Approximation with $q = 2$}
  \end{subfigure}
  \begin{subfigure}[t]{3.9cm}
    \centering
    \includegraphics[width = 3.9cm]{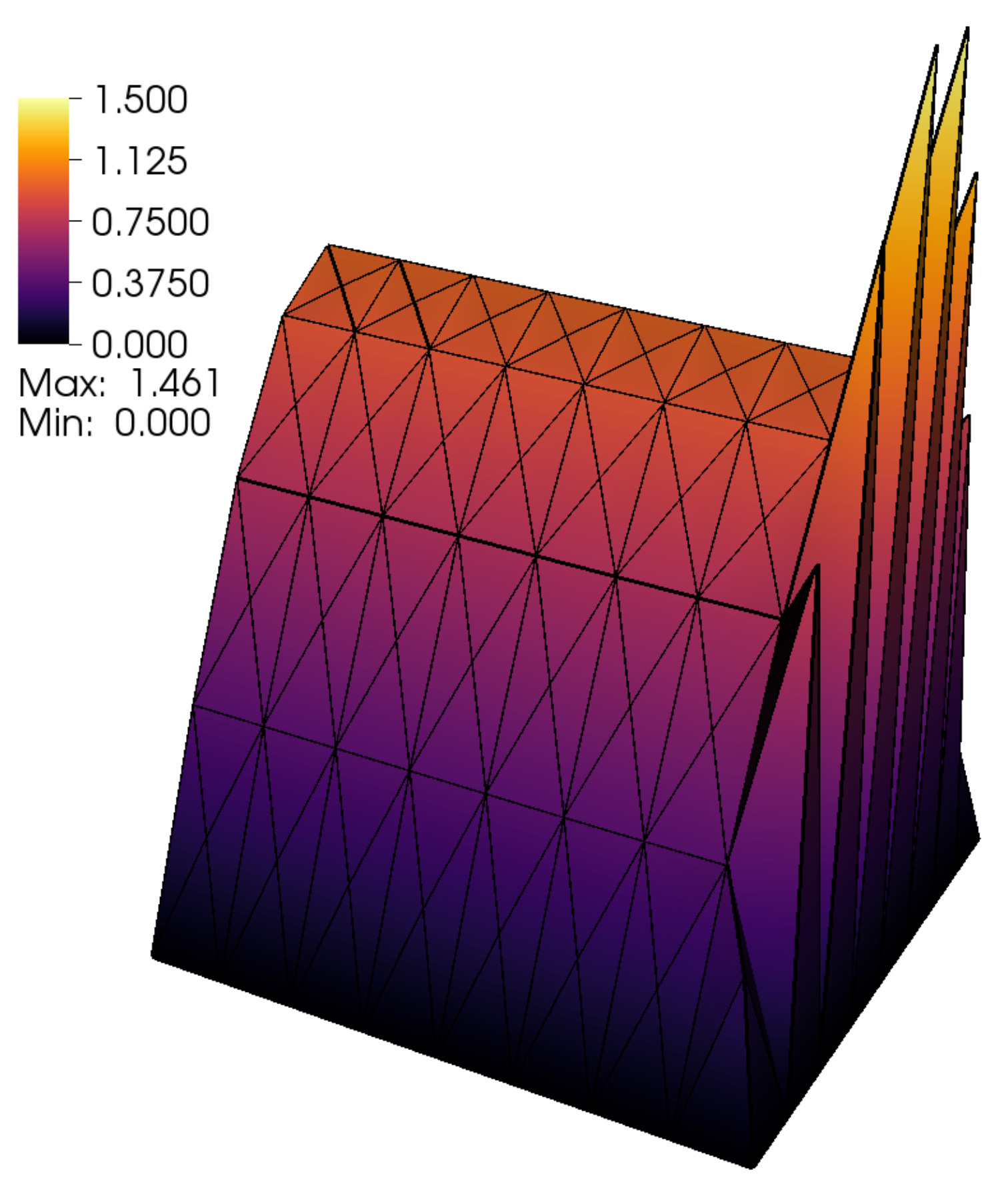}
    \caption{Approximation with $q = 1.01$}
  \end{subfigure}
    \caption{Approximation of the analytical solution of Example \ref{ex:eriksson-johnson} (Eriksson-Johnson model problem) on Mesh 3 for $\varepsilon = 10^{-6}$.}\label{fig:ej_crisscross}
\end{figure}

First, we compare the approximations on Mesh 2 and Mesh 3 for $q=2$ and $q=1.01$ with the interpolant of the analytical solution in $U_n$. Figure \ref{fig:ej_left_right} shows the approximations on Mesh 2; for $q = 2$ we  clearly observe overshoots near the boundary layer (center) compared to the interpolant of the analytical solution (left). On the other hand,
for $q = 1.01 $ we can see that the approximation approaches the interpolant of the analytical solution.
Figure \ref{fig:ej_crisscross} shows the approximations on Mesh 3; in this case we can see that the overshoot is nearly the same for $q = 2$ and $q = 1.01$. From \cite{Houston2019a} we note that this is the qualitative behaviour that the $L^q$-best approximation exhibits.

\begin{figure}[!t]
\centering
\includegraphics[width = 9.5cm]{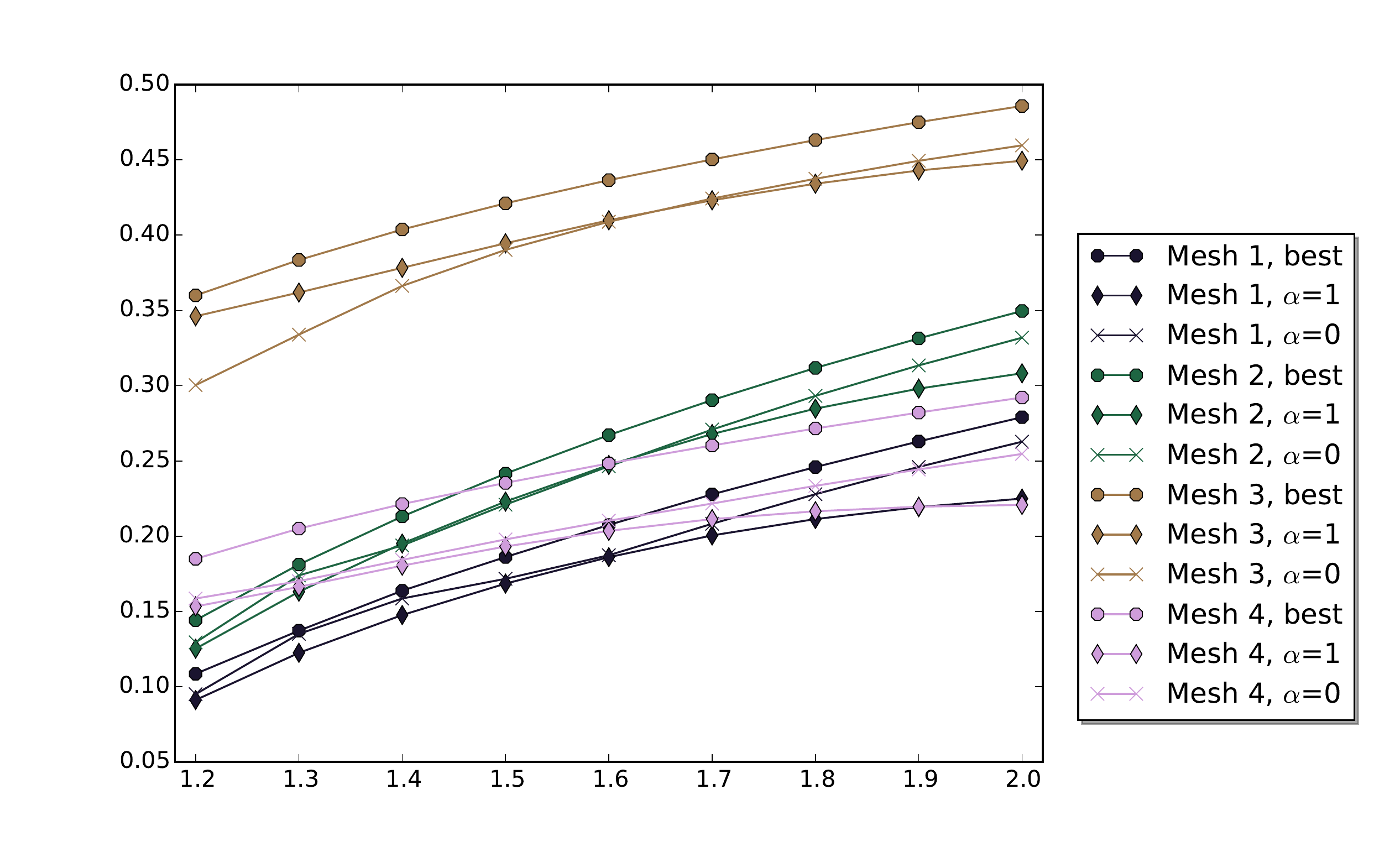}
\caption{Example \ref{ex:eriksson-johnson}: $\max(u_n -\hat{u})$, where $\hat{u}$ is the interpolant of the analytical solution $u$ in the space $U_n$, for the $L^q(\Omega)$-best approximation and the finite element approximations with $\alpha=1$ and $\alpha =0$.}
\label{fig:overshoot_2D}
\end{figure}

With this in mind, we have numerically computed the $L^q(\Omega)$-best approximation of the analytical solution for several choices of $q$ and  have also computed the finite element approximations with $\omega(x) =1$ and both, $\alpha=1$ and $\alpha = 0$, with weak boundary conditions imposed on $r_m$ on the inflow boundary for the same values of $q$.
Figure \ref{fig:overshoot_2D} shows the maximal error between the interpolant of the analytical solution and the approximation for the $L^q(\Omega)$-best approximation and the error between the interpolant and the numerical approximations for $\alpha = 1$ and $\alpha =0$. We can clearly see that the overshoot observed for our proposed scheme is very similar to the overshoot in the $L^q(\Omega)$-best approximation. Therefore, understanding for which meshes discontinuities/layers are captured as sharply as the grid allows by the $L^q(\Omega)$-best approximation can serve as an indicator regarding whether the over- and undershoots dissappear as $q \rightarrow 1$.  This is  the subject of \cite{Houston2019a}, where the $L^q$-best approximation of discontinuous functions is analysed in detail in certain cases.

\subsubsection{Boundary Layer in a Corner of the Domain}
\begin{figure}[!t]
  \centering
\begin{subfigure}[t]{6cm}
  \centering
  \includegraphics[width = 4.3cm]{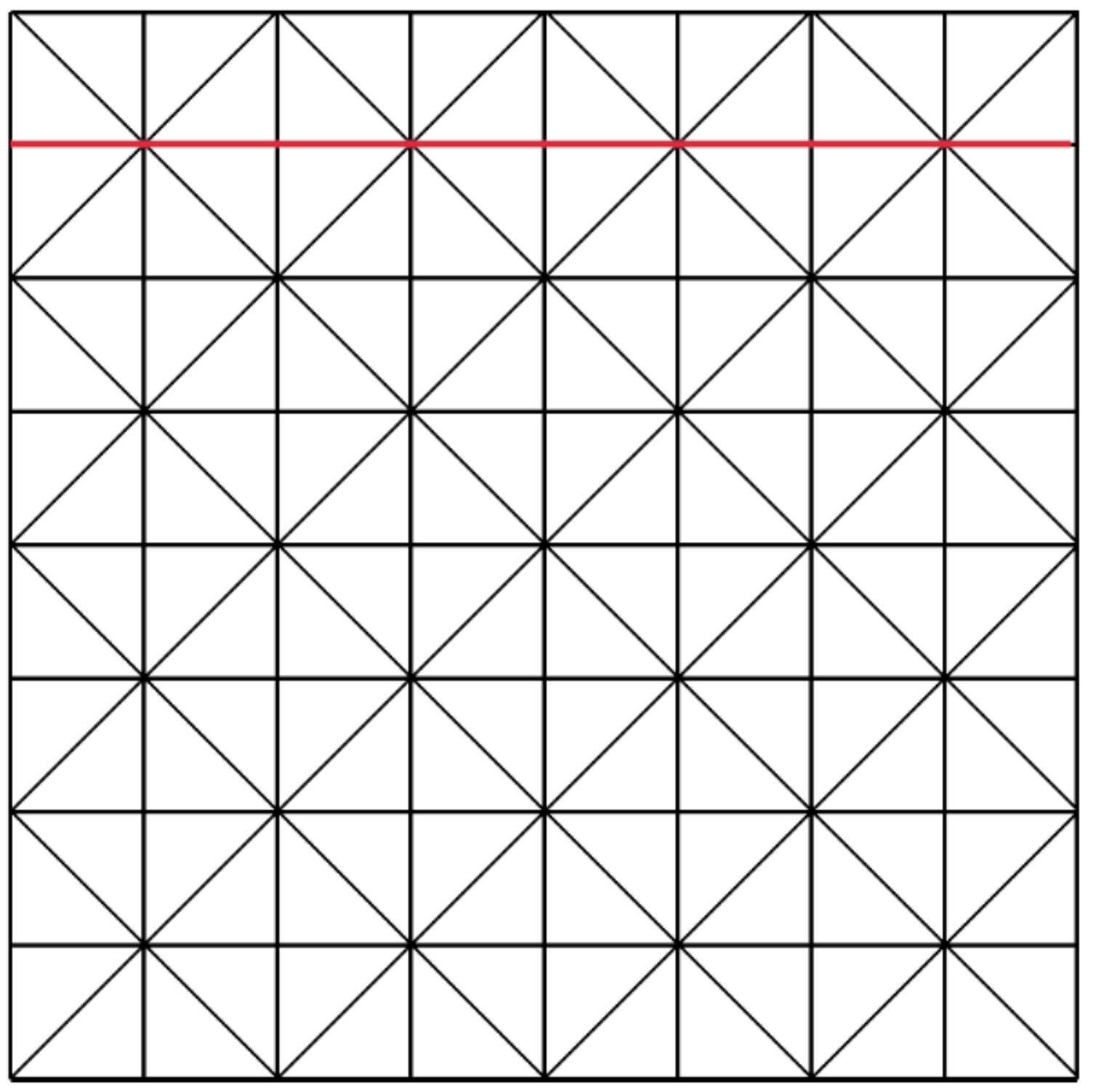}
  \caption{Mesh 2}\label{fig:mesh2_redline}
\end{subfigure}
\begin{subfigure}[t]{6cm}
  \centering
  \includegraphics[width = 4.3cm]{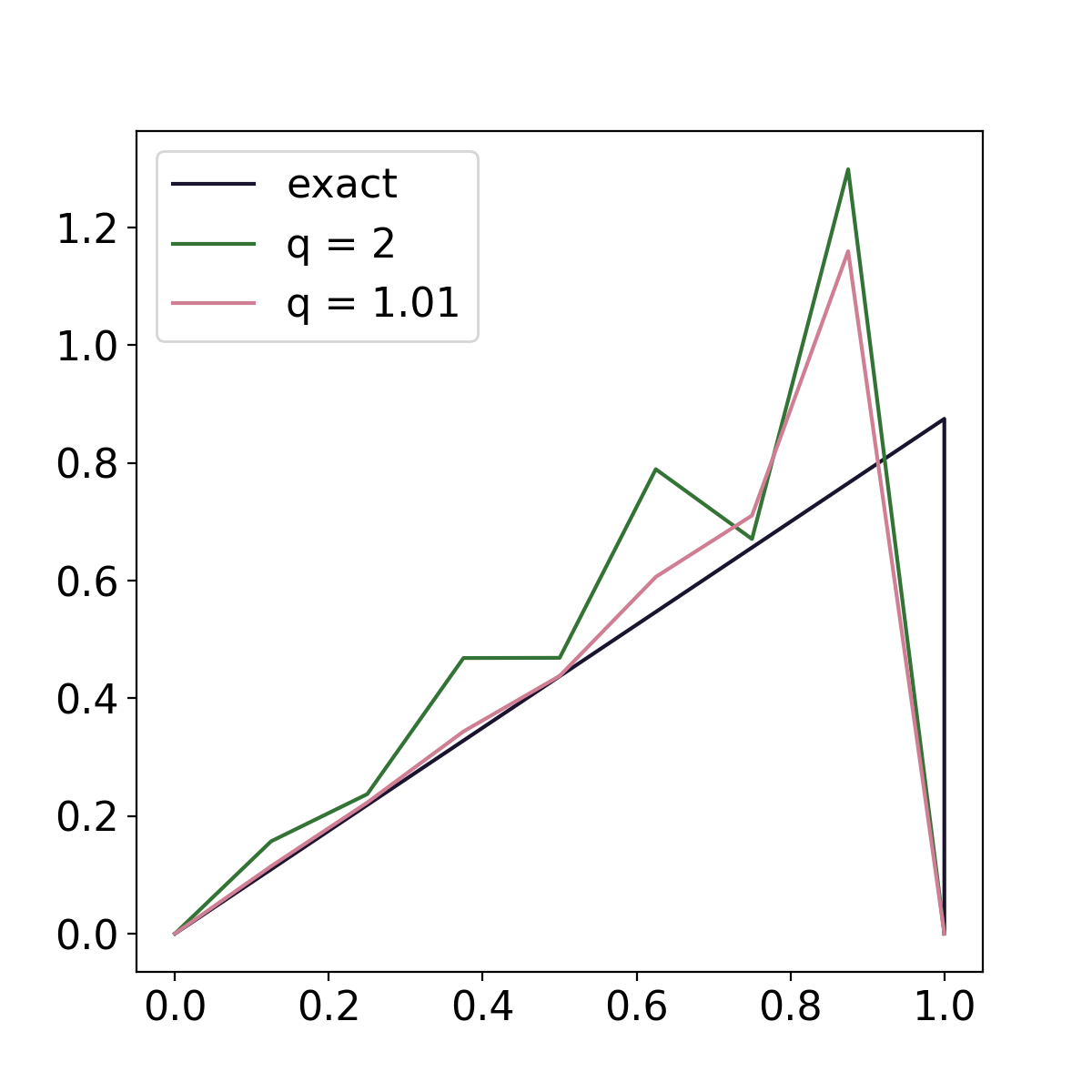}
  \caption{Exact solution $u$ and approximation $u_n$ along red line as indicated in Figure \ref{fig:mesh2_redline} for $q = 2$ and $q = 1.01$}
\end{subfigure}\\
\begin{subfigure}[t]{6cm}
  \centering
  \includegraphics[width = 4.3cm]{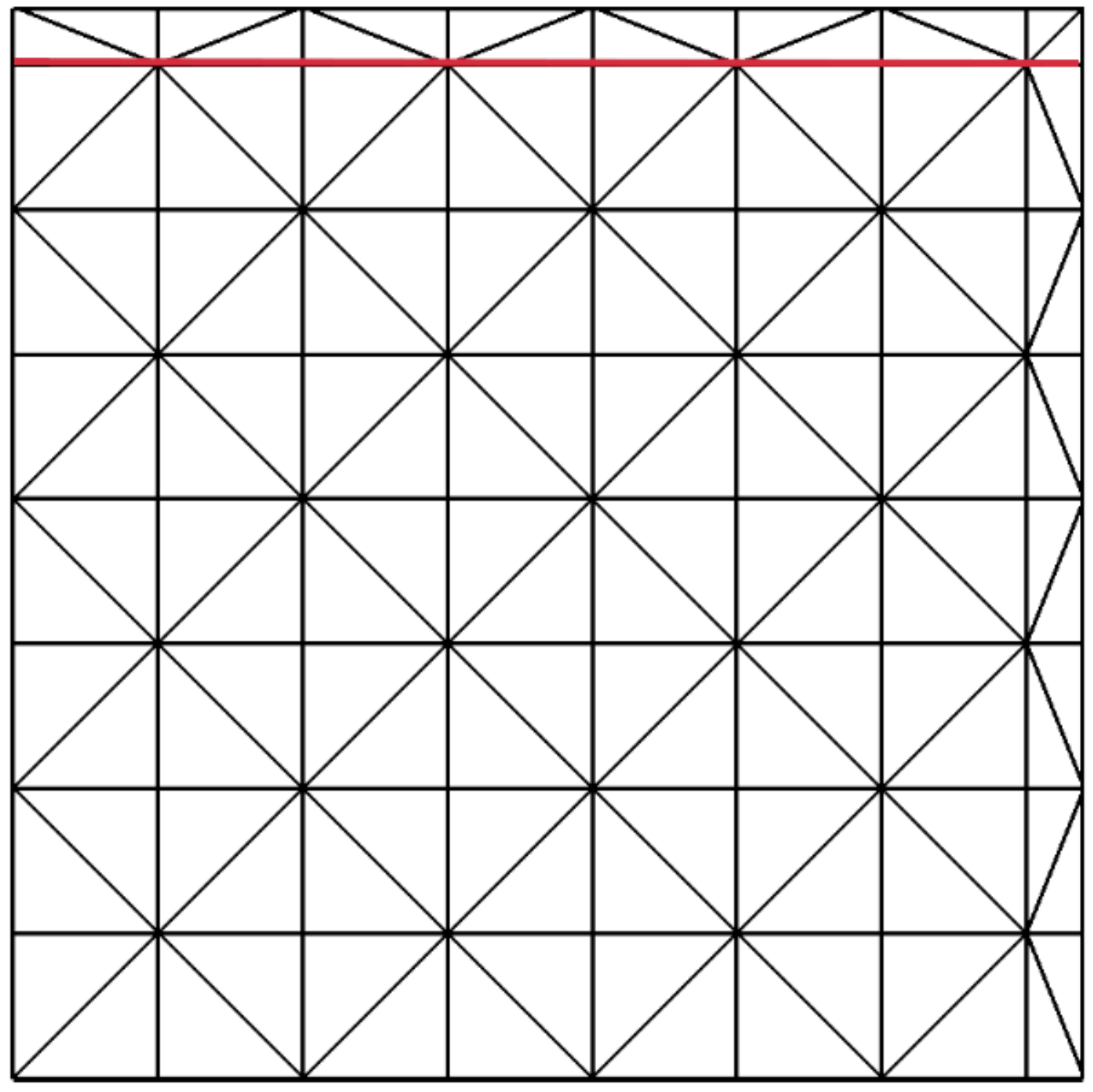}
  \caption{Modified version of Mesh 2}\label{fig:mesh2_modified}
\end{subfigure}
\begin{subfigure}[t]{6cm}
  \centering
  \includegraphics[width = 4.3cm]{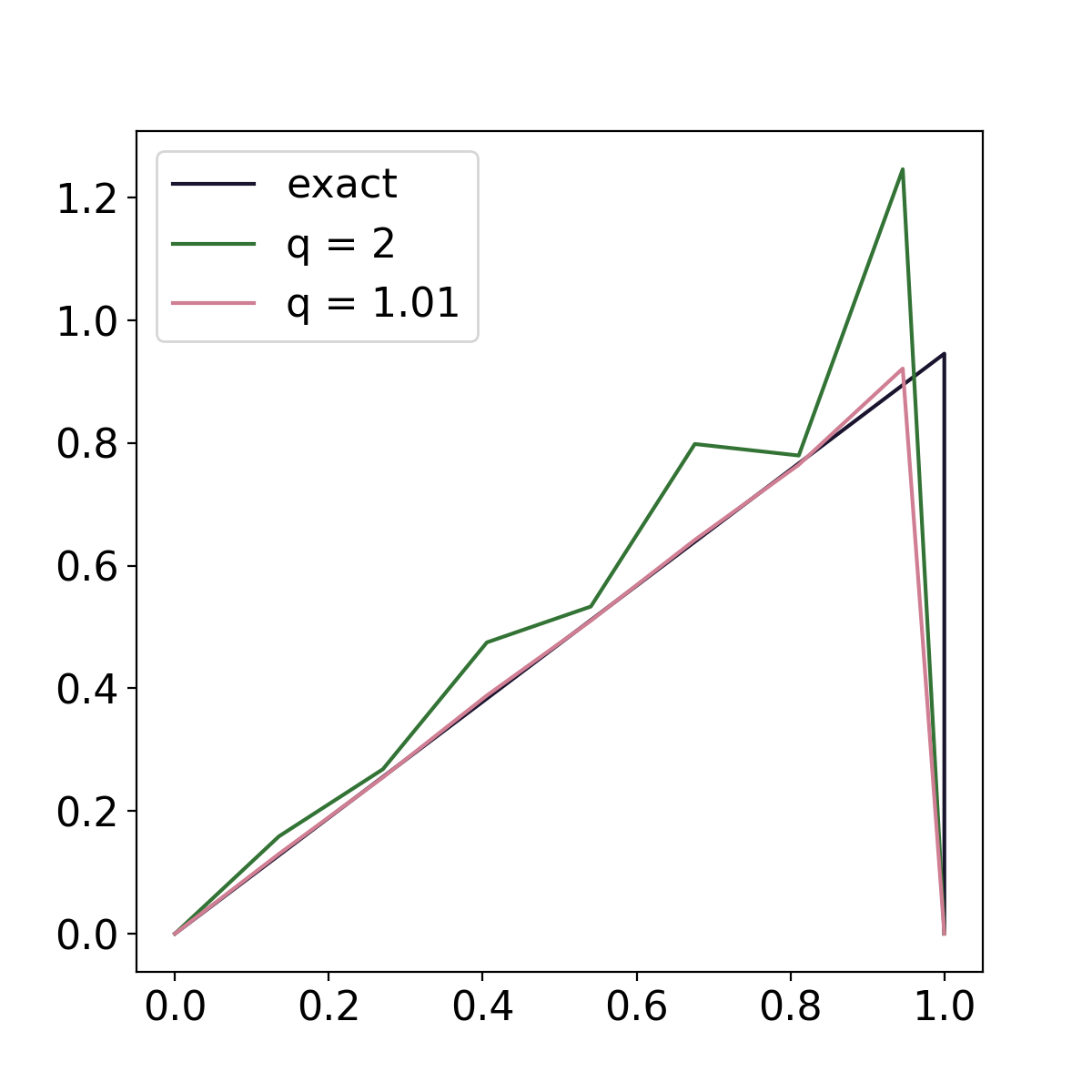}
  \caption{Exact solution $u$ and approximation $u_n$ along red line as indicated in Figure \ref{fig:mesh2_modified} for $q = 2$ and $q = 1.01$}
\end{subfigure}
\caption{Example \ref{ex:skew} with $\bm b = (2,1)^T$: Mesh and approximation for $\varepsilon = 10^{-6}$ with $p_n = 1$ and $p_m = 8$.}\label{fig:skew}
\end{figure}

In this section we  consider Example \ref{ex:skew} with $\bm b = (2,1)^T$. We will see that on Mesh 2 in this case the overshoot does not disappear as $q\rightarrow 1$  for our method. This is due to the layer not only appearing along an edge of the unit square but also near the corner
$(x,y) = (1,1)$. With this in mind, we suggest an alternative mesh that is only a slight modification of Mesh 2 for which the overshoot disappears.
Figure \ref{fig:skew} shows the two versions of Mesh 2 that we have used and the
corresponding approximations for $q=2$ and $q=1.01$ along the line through the
interior nodes closest to the boundary $y=1$. We can see that on Mesh 2, the overshoot
near the corner does not disappear as $q \rightarrow 1$, whereas it does reduce significantly away from the corner. For the modified version of Mesh 2, the overshoot disappears everywhere as $q \rightarrow 1$.

\begin{figure}[!t]
  \centering
  \begin{tikzpicture}[thick,scale=0.9, every node/.style={scale=0.9}]
    \fill[green!30!white] (0,0) -- (4, 0)-- (4,4) -- (0,4) -- (0,0);
    \fill[blue!50!white]  (0,0) -- (2, 0)-- (2,2) -- (0,2) -- (0,0) ;
\draw [thick] (0,0) rectangle (4,4);
\draw [thick](0,0) -- (4,4);
\draw [thick](4,0) -- (0,4);
\draw [thick] (0,2) -- (4,2);
\draw [thick] (2,0) -- (2,4);
\draw [thick,->] (4.2,2)--(4.8,2);
\end{tikzpicture}
\begin{tikzpicture}[thick,scale=0.9, every node/.style={scale=0.9}]
  \fill[green!30!white] (0,0) -- (4, 0)-- (4,4) -- (0,4) -- (0,0);
  \fill[blue!50!white]  (0,0) -- (2.857, 0)-- (2.857,2.857) -- (0,2.857) -- (0,0) ;
  \draw [thick] (0,0) rectangle (4,4);
  \draw [thick](0,0) -- (4,4);
  \draw [thick](4,0) -- (2.857,2.857);
  \draw [thick] (2.857,2.857) -- (0,4);
  \draw [thick] (0,2.857) -- (4,2.857);
  \draw [thick] (2.857,0) -- (2.857,4);
\end{tikzpicture}
\caption{Modification of Mesh 2: Elements connected to the boundary are marked in green, elements separated from the boundary in blue.} \label{fig:mesh_modification}
\end{figure}
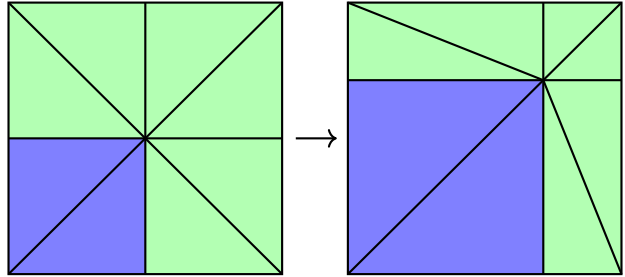

To understand the improvement when using this modified mesh, we focus on the interior node closest to the corner $(x,y) = (1,1)$. Figure
\ref{fig:mesh_modification} shows all elements connected to this node. The elements with either a node or an edge on the boundary are marked in green; the elements separated from the boundary are blue. From \cite{Houston2019a} we can infer that the
$L^1$-best approximation does not exhibit overshoots if the volume of the green area is smaller or
equal to the volume of the blue area. This is obviously violated on Mesh 2 with the green area being three times as large as the blue area. The modified version of
the mesh is  designed to satisfy this condition by changing the mesh as indicated
in Figure \ref{fig:mesh_modification}.

\subsubsection{Example with an Interior Layer}
\begin{figure}[!t]
  \centering
  \begin{subfigure}[t]{3.9cm}
  \includegraphics[width = 3.9cm]{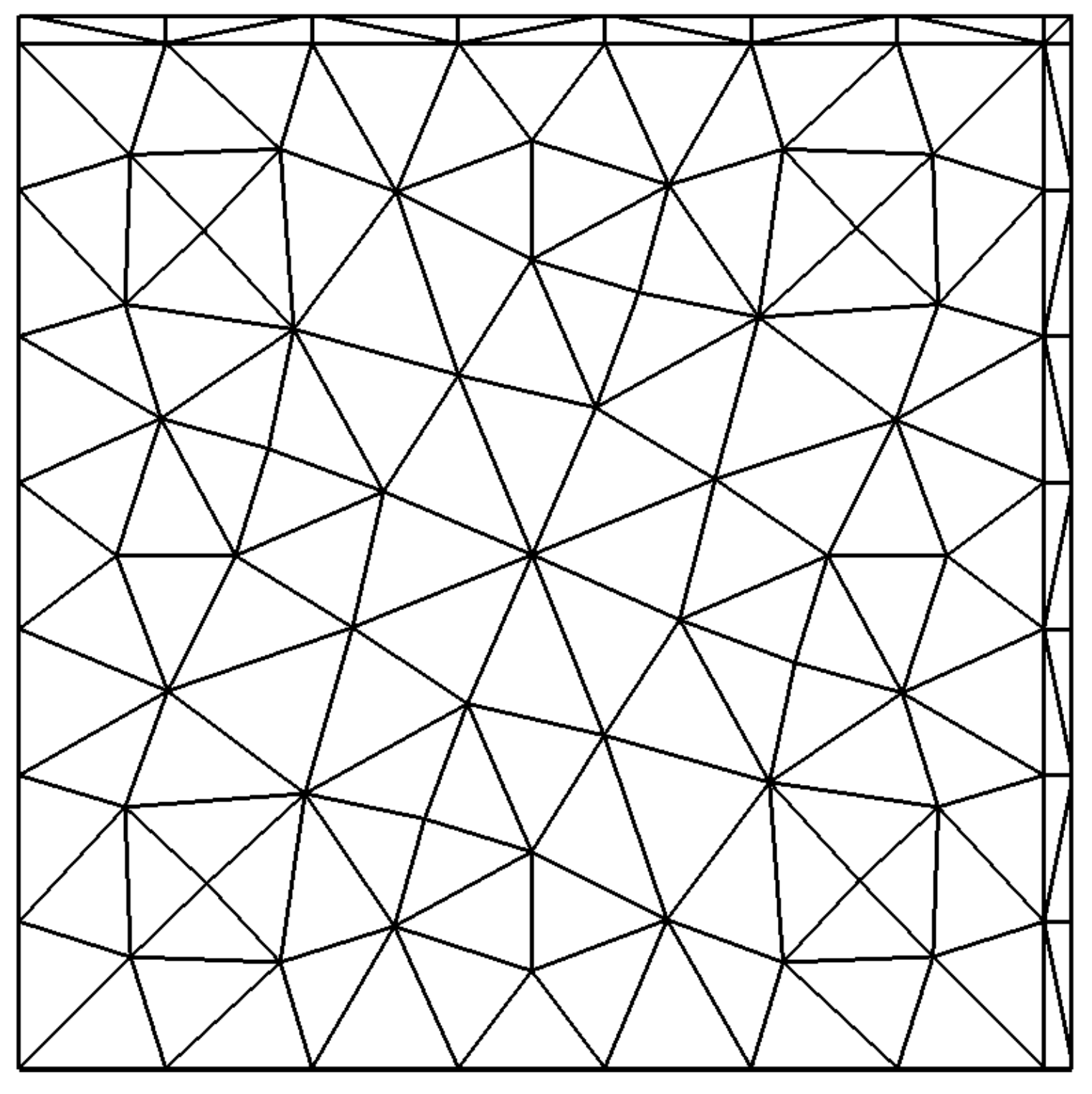}
  \caption{Mesh A}
  \end{subfigure}
  \begin{subfigure}[t]{3.9cm}
  \includegraphics[width = 3.9cm]{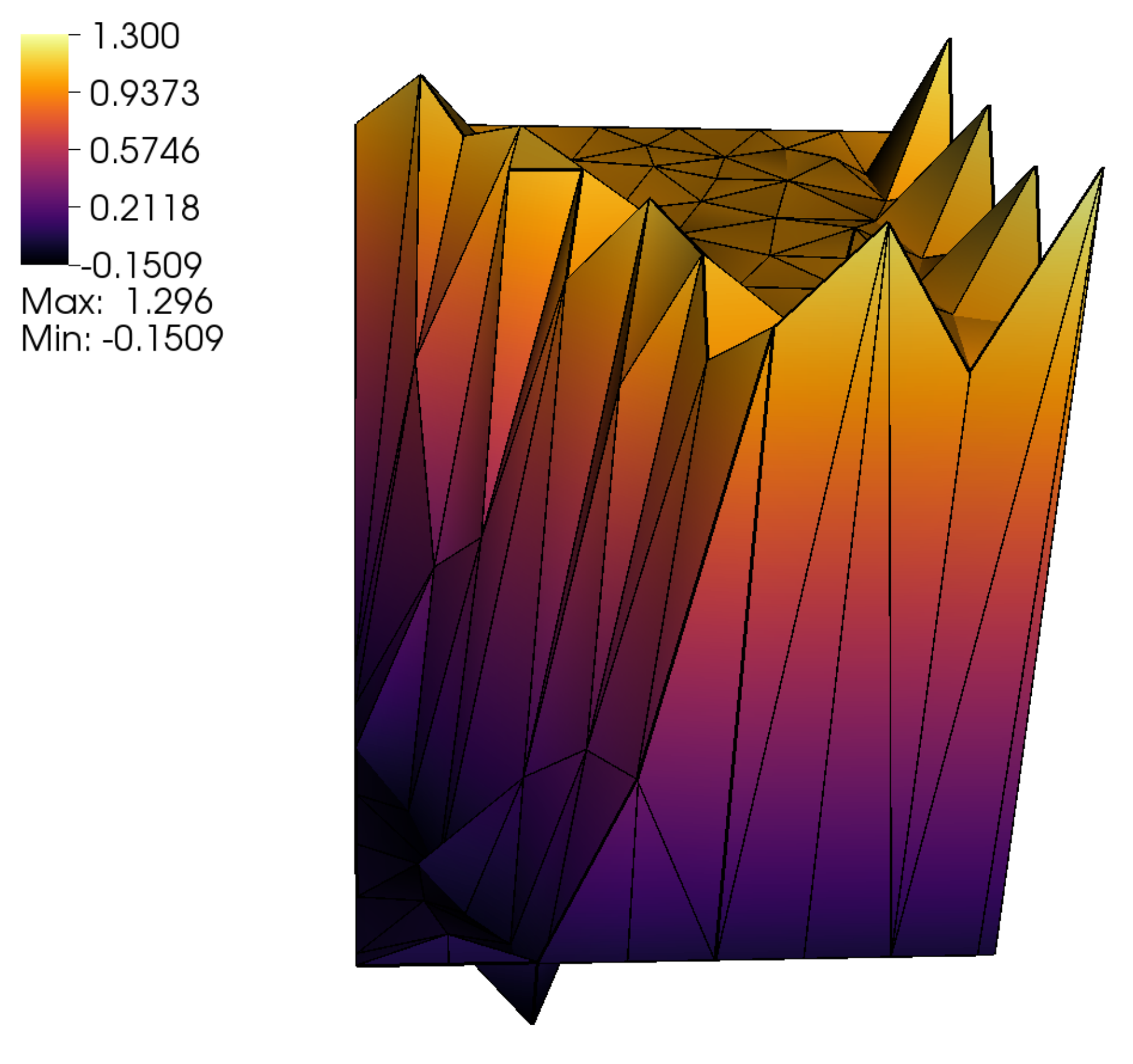}
  \caption{Approximation with $q = 2$ on Mesh A.}
  \end{subfigure}
  \begin{subfigure}[t]{3.9cm}
  \includegraphics[width = 3.9cm]{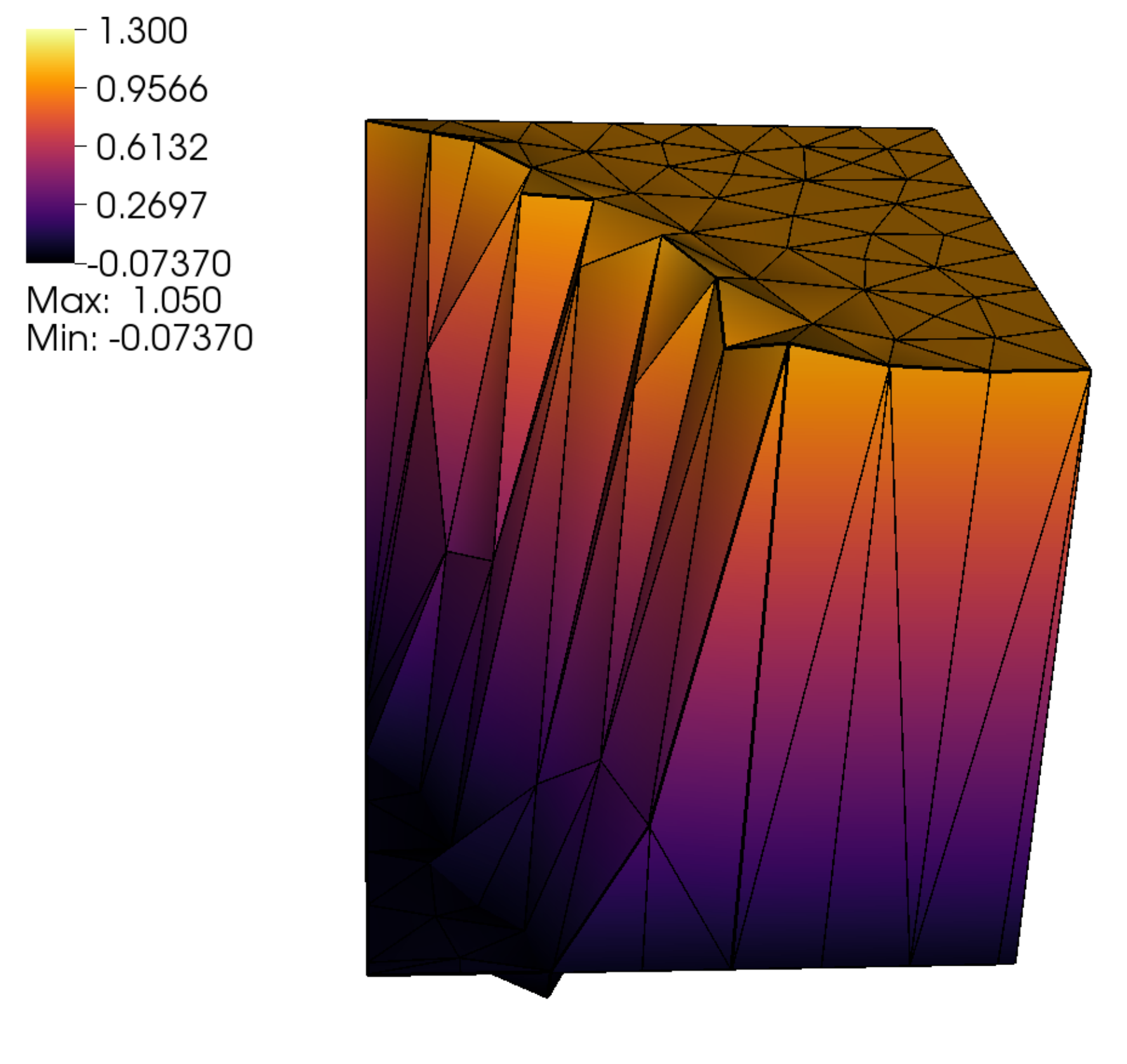}
  \caption{Approximation with $q = 1.1$ on Mesh A.}
  \end{subfigure}\\
  \begin{subfigure}[t]{3.9cm}
  \includegraphics[width = 3.9cm]{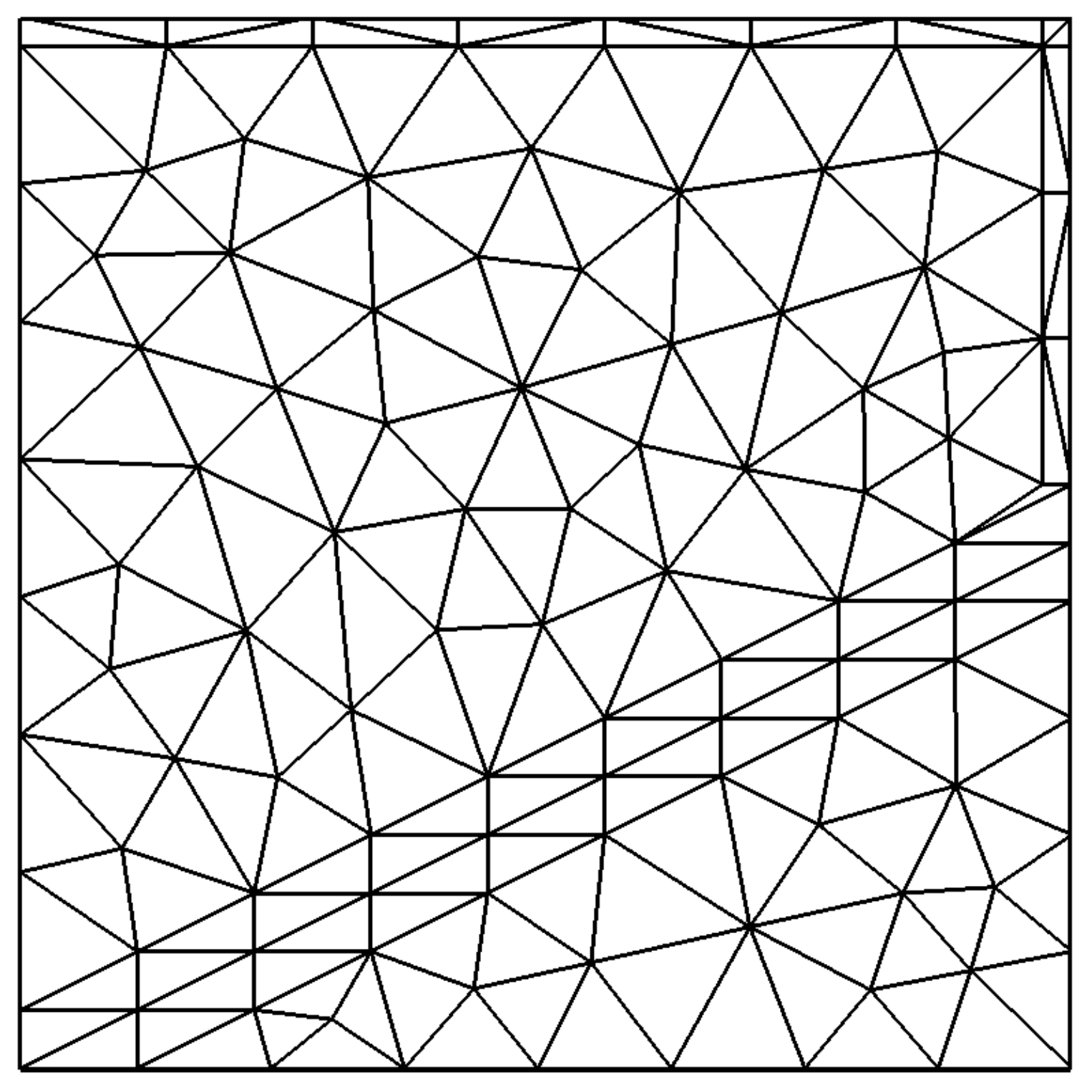}
  \caption{Mesh B}
  \end{subfigure}
  \begin{subfigure}[t]{3.9cm}
  \includegraphics[width = 3.9cm]{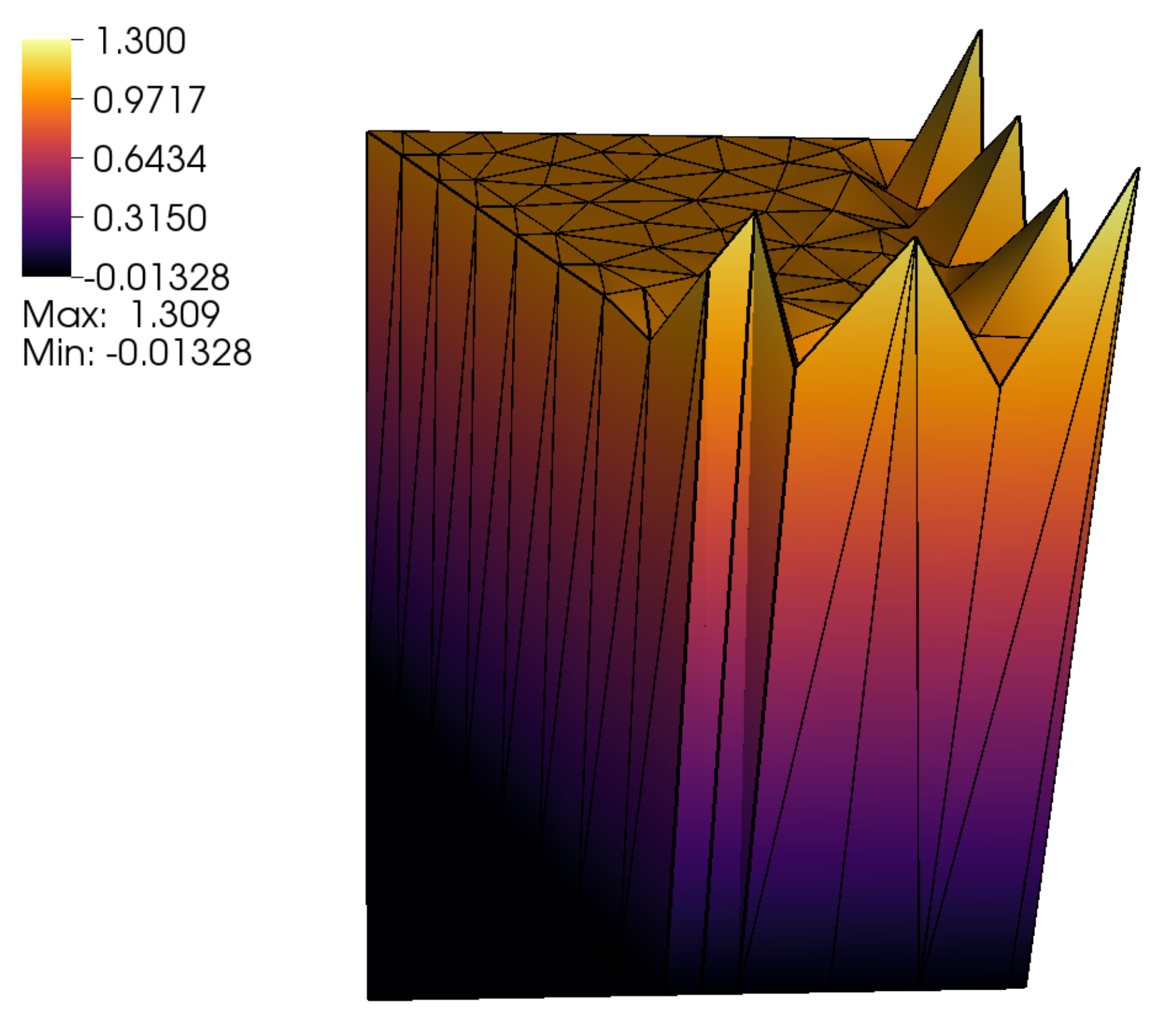}
  \caption{Approximation with $q = 2$ on Mesh B.}
  \end{subfigure}
  \begin{subfigure}[t]{3.9cm}
  \includegraphics[width = 3.9cm]{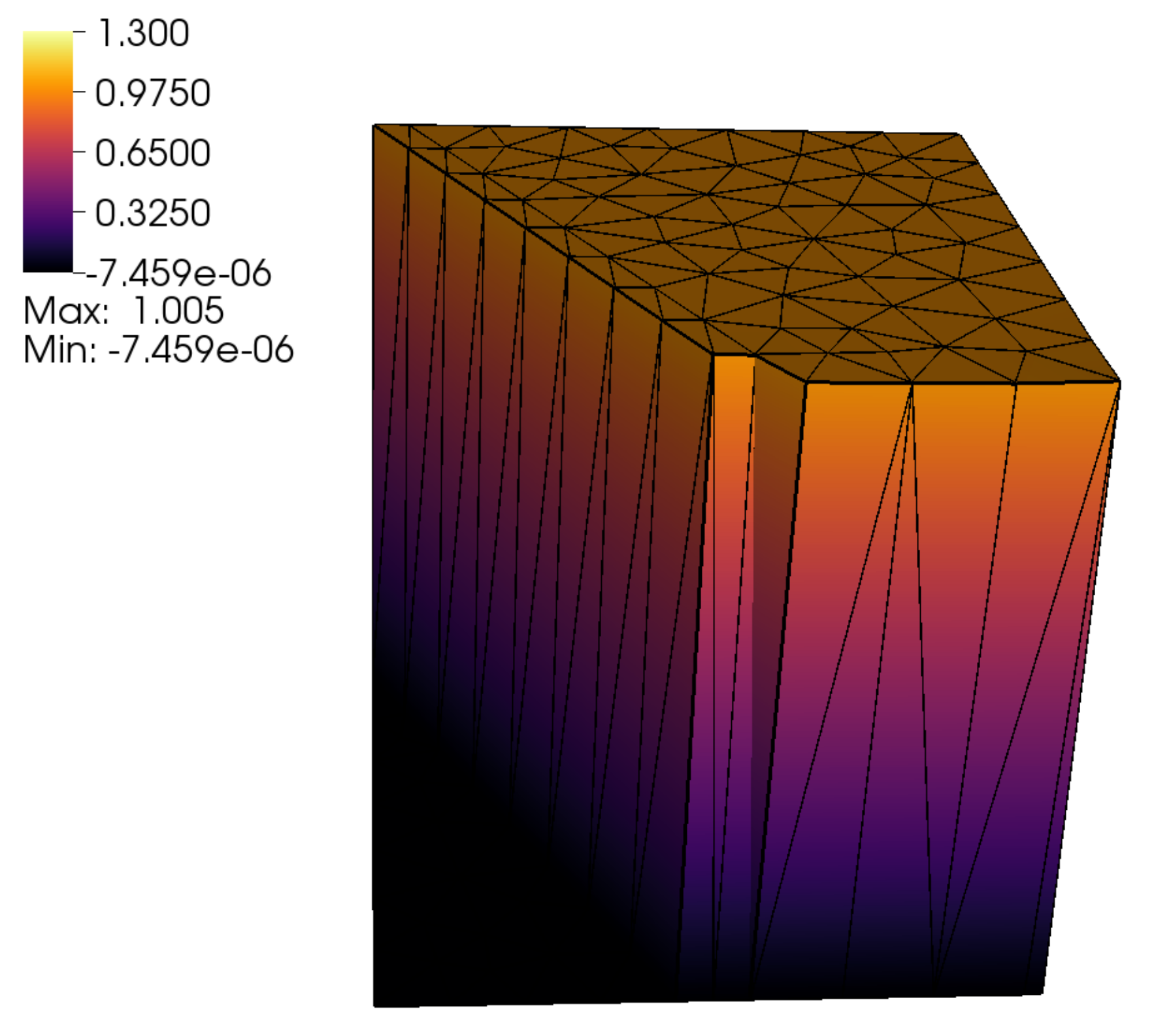}
  \caption{Approximation with $q = 1.1$ on Mesh B.}
  \end{subfigure}
    \caption{Example \ref{ex:interior_layer} with $\varepsilon = 10^{-6}$.}\label{fig:interior_layer}
\end{figure}

In this section we consider Example \ref{ex:interior_layer} and demonstrate that a mesh can be constructed such that the over- and undershoots both near the boundary layer and the interior layer can be eliminated as $q\rightarrow 1$.
To this end,
Figure \ref{fig:interior_layer} shows the approximation to Example  \ref{ex:interior_layer} for $q = 2$ and $q =1.1$ on two different meshes. The first mesh is designed such that the overshoot at the boundary layer disappears as $q \rightarrow 1$; this is already clearly visible for $q = 1.1$. Along the interior layer, the overshoot also reduces significantly, but some over- and undershoots are still present for $q=1.1$. The second mesh is additionally designed to align with the interior layer in such a way that we can expect the over- and undershoots to disappear as $q\rightarrow 1$.
We can see, that this mesh even nearly eliminates the overshoot along the interior layer for $q = 2$ and for $q =1.1$ all over- and undershoots along both layers have essentially vanished.

\begin{figure}[!t]
  \centering
  \begin{subfigure}[t]{4cm}
  \includegraphics[width = 4cm]{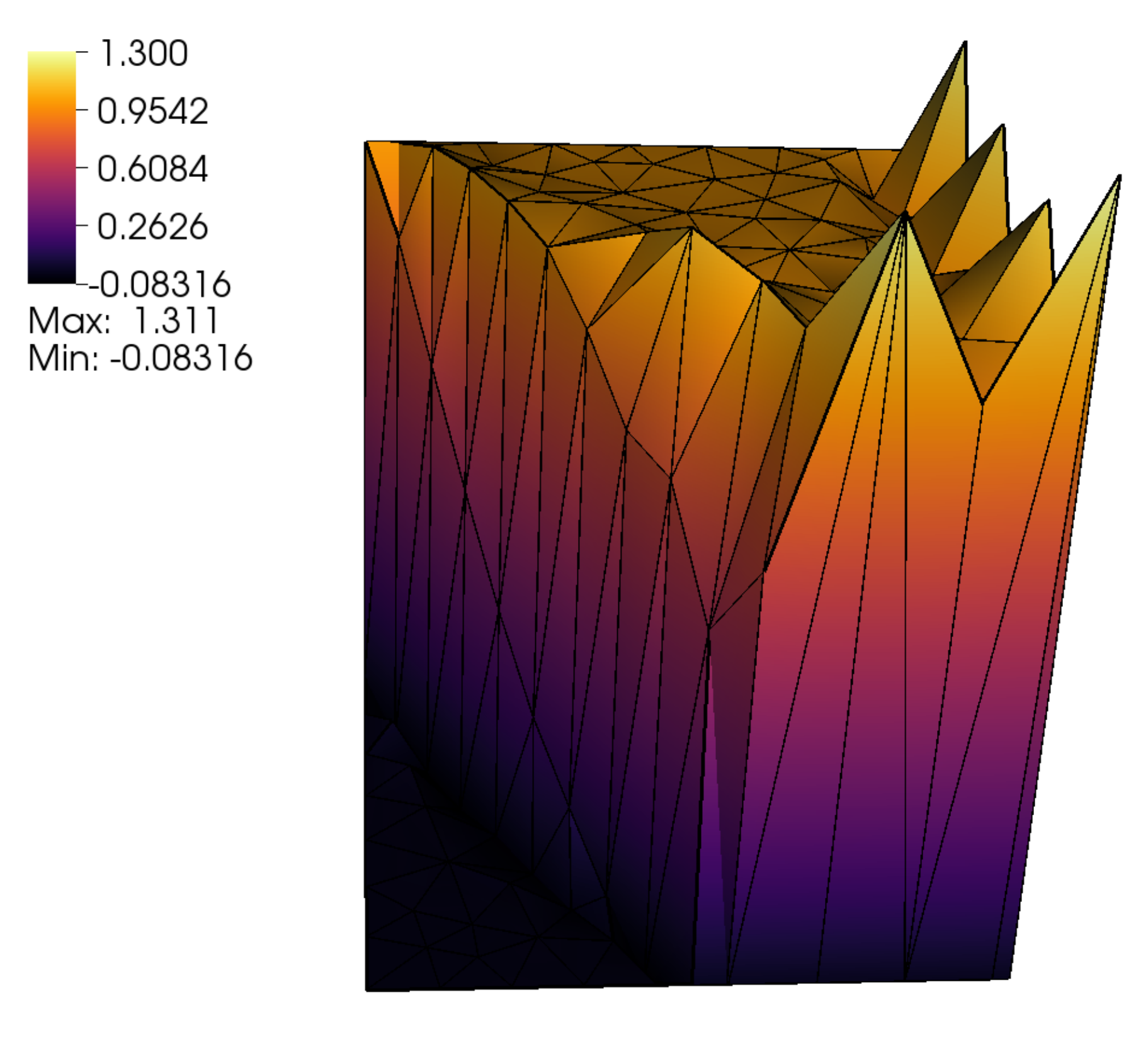}
  \caption{Approximation with $q = 2$ and $\bm b = (2, 1.2)^T$ on Mesh B.}
  \end{subfigure}
  \begin{subfigure}[t]{4cm}
  \includegraphics[width = 4cm]{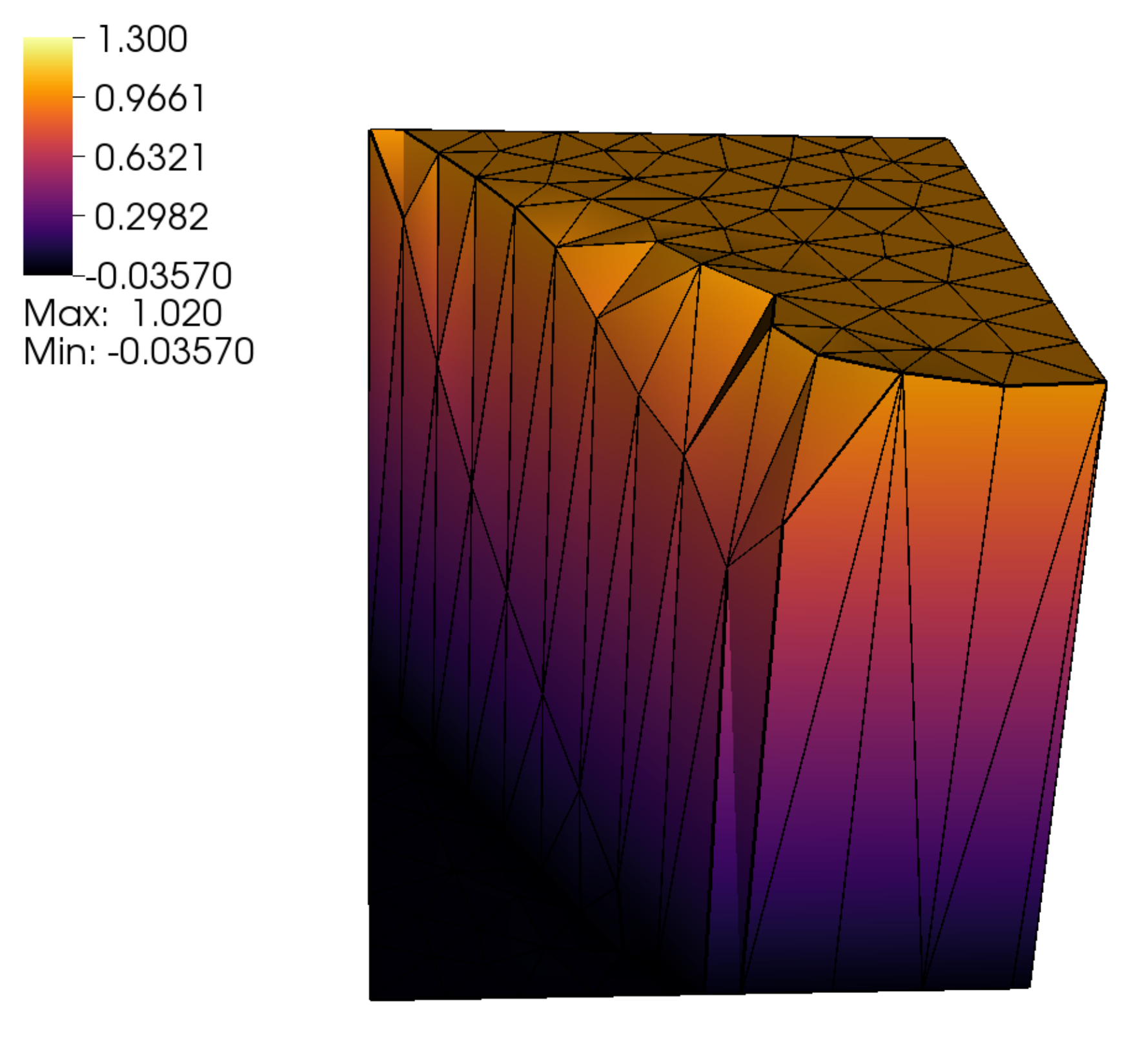}
  \caption{Approximation with $q = 1.1$ and $\bm b = (2, 1.2)^T$ on Mesh B.}
  \end{subfigure}

  \begin{subfigure}[t]{4cm}
  \includegraphics[width = 4cm]{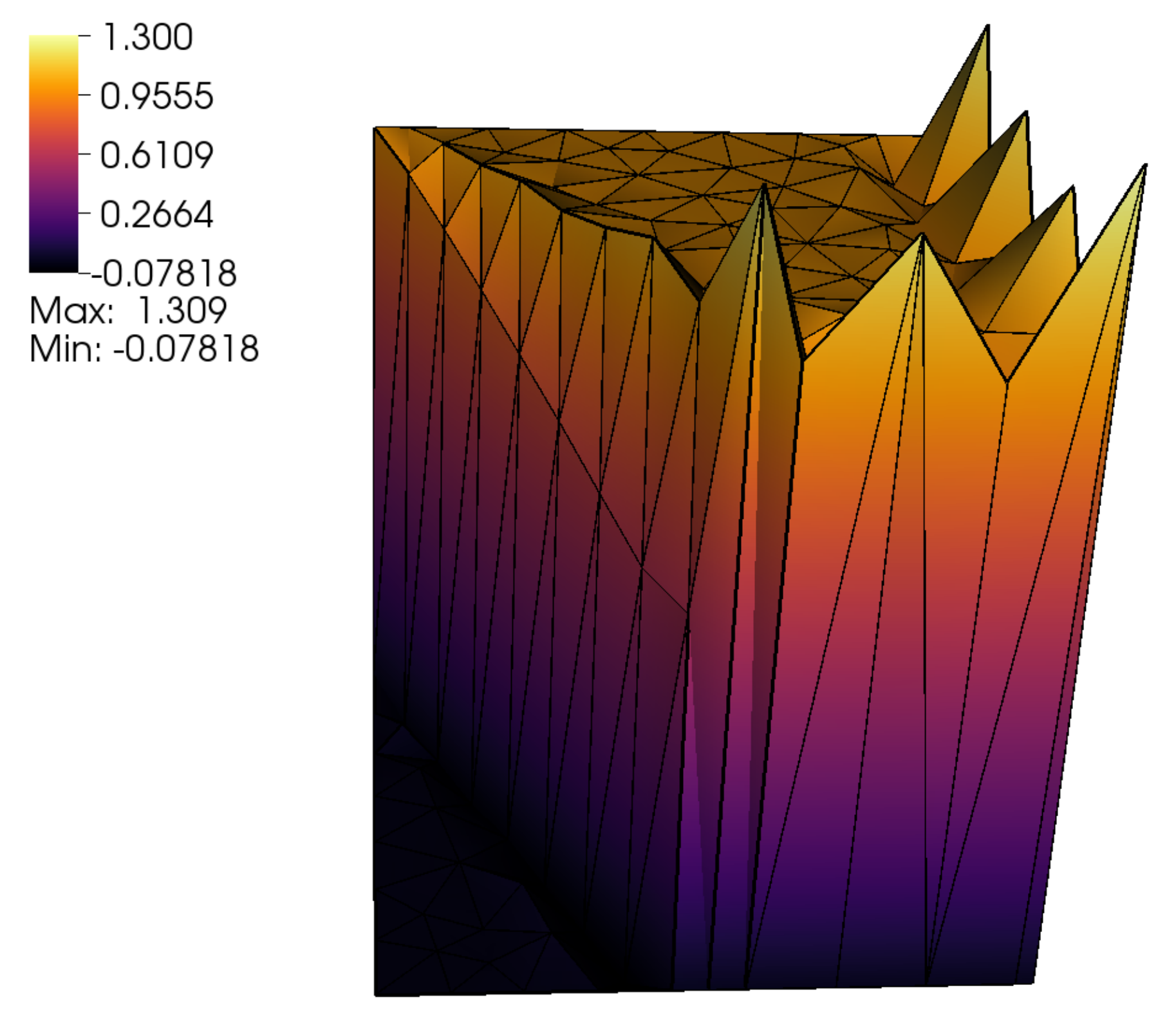}
  \caption{Approximation with $q = 2$ and $\bm b = (2, 1.06)^T$ on Mesh B.}
  \end{subfigure}
  \begin{subfigure}[t]{4cm}
  \includegraphics[width = 4cm]{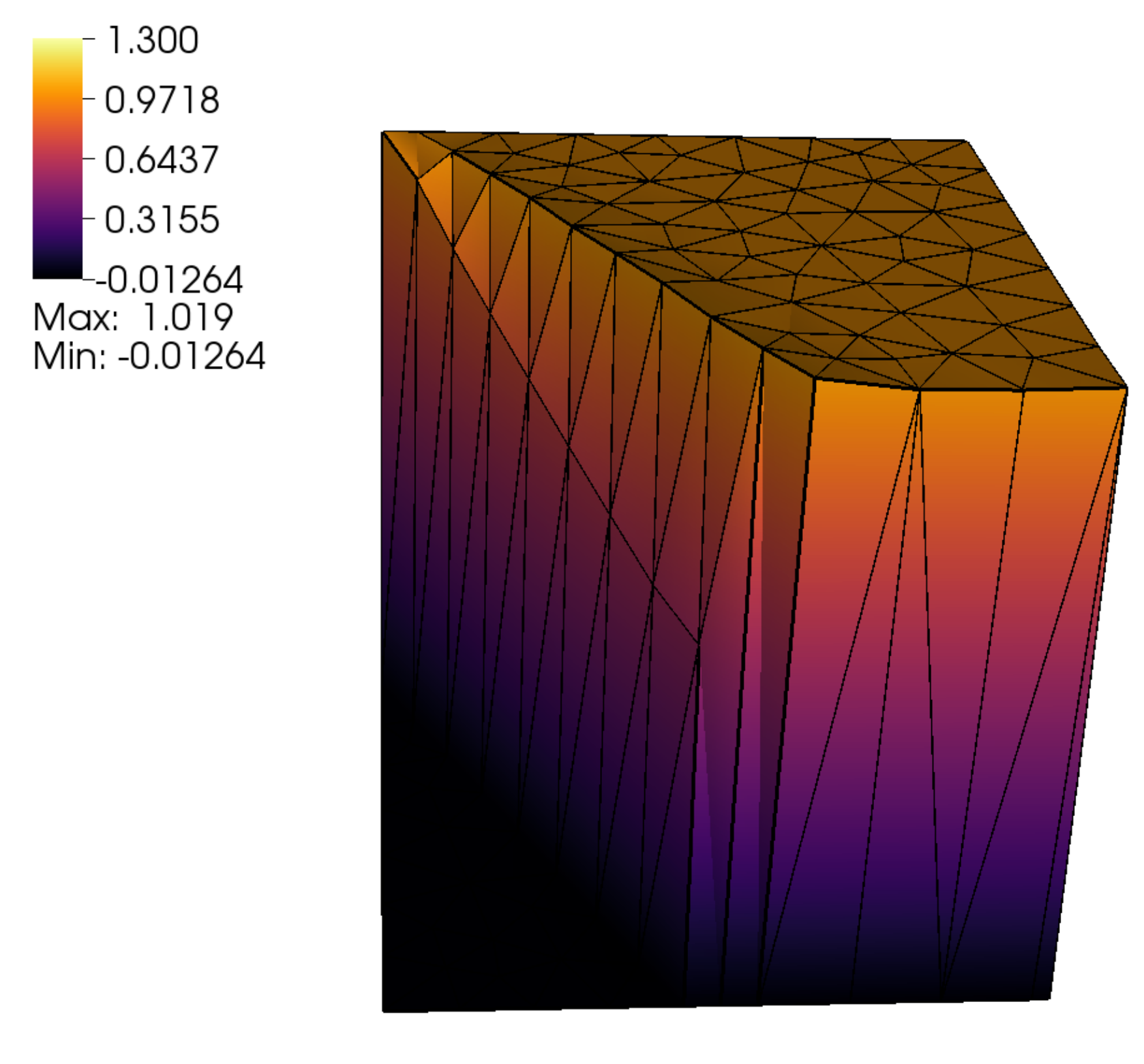}
  \caption{Approximation with $q = 1.1$ and $\bm b = (2, 1.06)^T$ on Mesh B.}
  \end{subfigure}
\caption{Example \ref{ex:interior_layer} with $\varepsilon = 10^{-6}$, $\bm b = (2, 1.2)^T$(top) and$\bm b = (2,1.06)^T$ (bottom).}\label{fig:interior_layer_sensitivity}
\end{figure}

Figure \ref{fig:interior_layer_sensitivity} shows the approximations for $q = 2$ and $q = 1.1$ if $\bm b$ is not perfectly aligned with the mesh, i.e., $\bm b = (2, 1.2)^T$ and $\bm b = (2, 1.06)^T$.
The difference between these two choices for $\bm b$ is that in the latter case the layer is still contained in between the two lines parallel to the line $(x,y) = t(2,1)$, whereas in the former case it is not.
In both cases, we clearly observe overshoots along the interior layer for $q = 2$ which are reduced for $q = 1.1$.
Both, for $q = 2$ and $q=1.1$ the overshoot reduces slowly the closer $\bm b$ is to $(2,1)^T$.  For $\bm b = (2, 1.06)^T$ and $q = 1.1$ the overshoot is barely visible in Figure \ref{fig:interior_layer_sensitivity}, but comparing the maximum and minimum values of the approximations reveals that the over- and undershoots are similar to the case $\bm b = (2, 1.2)$ in magnitude.

\section{Conclusions and Future Directions}\label{sec:conclusion}
In this article, we generalized the framework employed for DPG methods to more general Banach spaces. This framework in principle allows us to choose the norm in which the solution is approximated, provided  a suitable well-posed variational formulation with  robust constants and a test-norm that is computable is given. Another challenge in designing such a method is that a Fortin-condition must be satisfied for the discrete spaces $U_n$ and $V_m$. The recent work \cite{Demkowicz2019} describes a double adaptivity approach for the DPG method to circumvent this challenge in the context of DPG methods. Due to the close relationship between DPG methods and our approach this has potential to be generalized to the Banach space framework.

Furthermore, we have demonstrated how to
design a finite element method  for the convection-diffusion-reaction problem that in the convection-dominated case yields solutions that qualitatively behave like the $L^q(\Omega)$-best approximation of the analytical solution. This means that on meshes where the $L^1(\Omega)$-best approximation does not contain any over- and undershoots, these
can also be avoided by taking $q\rightarrow 1$ in our proposed method.
Indeed, the final two examples presented in this article demonstrate that it is possible to adjust meshes to specific situations in order to eliminate over- and undershoots.

 An important implication of the connection to the $\Lp{q}$-best approximation is that we can use the insights from \cite{Houston2019a} to design suitable meshes. Note that even in one dimension the over- and undershoots will not disappear for certain non-uniform meshes. Roughly speaking, meshes where the elements near the discontinuity or layer present in the analytical solution  are smaller than at a distance from those features are more favourable. Additionally, it is desirable that all interior nodes that are closest to such a feature are aligned roughly in parallel with the feature; for further details, we refer to \cite{Houston2019a}. This is similar to the observations used to design the modified version of the mesh in the final two examples. This suggests that if the mesh is refined in a suitable way near the layer, we can eliminate over- and undershoots. This observation could be used to design a mesh refinement strategy that modifies the mesh in such a way that the $L^1(\Omega)$-best approximation of the layer or discontinuity does not contain any over- or undershoots. Note, however, that even though we modified the mesh in the final example in a way that resembles refinement towards the boundary layer, the boundary layer clearly remains under resolved even on the modified mesh. It is therefore not necessary to fully resolve the layer to eliminate the oscillations in this case.

Another challenge of our proposed method is computational feasibility. Possible strategies to addressing this include reducing the degrees of freedom by choosing the test space $V_m$ adaptively, introducing broken test spaces $V_m$, similar to DPG methods, to allow for more effective parallelization, and developing a more efficient non-linear solver.
Despite the challenges outlined above, the numerical results for our proposed method show great potential for the application to non-linear problems whose analytical solutions contain sharp layers or discontinuities. It is well-known that higher-order monotonicity preserving methods can be expected to be non-linear; thus, the non-linear nature of our approach  by no means diminishes its potential.

\section*{Acknowledgements}
The authors would like to thank Gabriel Barrenechea, Jesse Chan, Leszek Demkowicz, Jay Gopalakrishnan, Ignacio Muga and Martin Stynes for fruitful discussions on the topic of this article.

\bibliography{main}

\end{document}